\documentclass[10pt]{article}

\usepackage{amsmath,amsthm,epsfig,amsfonts,amssymb}

\newtheorem{Theorem}{Theorem}

\newtheorem{Lemma}[Theorem]{Lemma}
\newtheorem{Proposition}[Theorem]{Proposition}
\newtheorem{Corollary}[Theorem]{Corollary}

\newcommand{\R}{R}
\newcommand{\Z}{Z}
\newcommand{\N}{N}
\newcommand{\I}{I}
\newcommand{\equalD}{\stackrel{\mathcal{D}}{=}}      

\numberwithin{equation}{section}

\begin{document}

\title{A phase diagram for a stochastic reaction diffusion system.}
\author{Carl Mueller, University of Rochester\\
and \\
Roger Tribe, Warwick University.}
\maketitle

\begin{abstract}
In this paper a stochastic reaction diffusion system is considered, 
which models the spread of a finite population reacting with a 
non-renewable resource in the presence of individual based noise. 
A two-parameter phase diagram is established to describe the large time evolution, 
distinguishing between certain death or possible life of the population.  
\end{abstract}


\section{Introduction} \label{introduction} \label{s1}
\subsection{Statement of results} \label{statement} \label{s1.1}
In this paper we study the stochastic reaction diffusion system on $\R^d$, for $d \leq 3$,
\begin{equation}  \label{spde}
\left\{ 
\begin{array}{l}
\partial_t u =  \Delta u + \beta uv -\gamma u + \sqrt{u} \, \dot{W}, \\ 
\partial_t v = -uv.
\end{array}
\right.
\end{equation}
One interpretation for this equation is that the solution $u \geq 0$ describes
the distribution of a population on $\R^d$ at times $t \geq 0$
and that $v \geq 0$ is the density of a nutrient gradually used up by the population. 
The noise $W$ is a space-time white noise, and the multiplicative factor 
$\sqrt{u}$ models `individual based' noise. This noise arises in situations where each
individual in the population contributes independently to the noise, so that
the variance of the noise is proportional to the population density.
In dimensions $d \leq 3$, linear scaling reduces the number of possible parameters to two. 
We have chosen to take the diffusion coefficient, the rate at which
the nutrient is used, and the noise coefficient all to be one, 
leaving two parameters $\beta,\gamma \geq 0$, which we think of as a 
reaction rate and a death rate.

The main results concern the long time behaviour of solutions.
Without noise there are travelling wave solutions and it is expected that
sufficient noise can destroy these waves. For our description of the 
long time behaviour, we take an initial condition $u_0=\mu$ in $\mathcal{M}(\R^d)$, the space of
finite measures in $\R^d$, and an initial 
nutrient level $v_0(x)$ taking a constant value which, by 
scaling again, we have taken to be $1$.
We say that {\it certain death} occurs for the parameters values $\beta,\gamma$ if 
\[
\mbox{$P \left[ u_t =0 \; \mbox{for large} \; 
t \right]=1$ for any initial condition
 $\mu \in \mathcal{M}(\R^d)$ and $v_0 = 1$.} 
\]
We say that {\it possible life} occurs for the parameter values $\beta,\gamma$ if 
\[
\mbox{$P \left[u_t \not =0 \; \mbox{for all} \; 
t \right]>0$ for any non-zero initial condition
 $\mu \in \mathcal{M}(\R^d)$ and $v_0=1$.} 
\]
The main results of the paper describe a phase diagram for the
parameter values $\beta, \gamma$ at which certain death or possible life occur.
\begin{Theorem} \label{mainresults} 
Consider solutions to (\ref{spde}) with initial condition $v_0 = 1$.
Then for any values of the parameters $\beta,\gamma \geq 0$ either 
possible life or certain death occurs. 
\begin{itemize}
\item When \textbf{d=3}, there exists a non-decreasing function $ \beta \to \Psi(\beta) \in (0,\beta)$ so that
when $0 \leq \gamma < \Psi(\beta)$ possible life occurs and when 
$\gamma > \Psi(\beta)$ certain death occurs. 
Moreover 
\[
0 < \liminf_{\beta \to 0} \beta^{-2} \Psi(\beta) \leq \limsup_{\beta \to 0} \beta^{-2} \Psi(\beta) < \infty,
\]
\item When \textbf{d=2} there exists a critical curve $\beta \to \Psi(\beta) \in [0,\beta)$ as in dimension $d=3$ 
(but allowing the possibility that $\Psi(\beta) = 0$ for small $\beta$).  
\item When \textbf{d=1} certain death occurs for any $\gamma, \beta \geq 0 $.
\end{itemize}
\end{Theorem}

\noindent \textbf{Remarks}

\noindent \textbf{1.} The estimates in the paper imply certain bounds on the critical curve. 
In section \ref{s6.3} we show that $\Psi(\beta) / \beta \to 1$ as $\beta \to \infty$. This asymptotic 
is not sharp, and the correct large $\beta$ asymptotic is, we believe, 
$\beta - \Psi(\beta) \sim \beta^{-2/(6-d)}$. The remarks in section \ref{s6.3} indicate that a certain 
improved spatial Markov property is required to establish these asymptotics, and the details are left 
to a subsequent paper. 
Many other questions remain, such as regularity of $\Psi$ and the
behaviour on the critical curve. 
Most importantly perhaps, we do not know the small $\beta$ behaviour in $d=2$, though
we conjecture (see section \ref{s7.3}) that there is certain death for $\gamma=0$ and small enough $\beta$.

\noindent \textbf{2.} 
The possible life behaviour is not an example of global coexistence of two species, 
for locally either the population $u_t$ becomes extinct or the nutrient 
levels converge to zero (see Lemma \ref{localextinction}). Rather one expects
the process to live on a moving front travelling through space, and behind this 
moving front the nutrient is used up and the process dies away. It is therefore 
more analogous to the weak survival results for 
the contact process on a homogeneous tree (see \cite{liggett} section I.4). 

\noindent \textbf{3.}
A key difficulty in dealing with many reaction diffusion systems 
is the lack of pathwise comparison results (which were crucial 
to the arguments for the scalar equation studied in \cite{mueller+tribe}). 
For example if two solutions 
$u,\tilde{u}$ satisfy $u_0 \leq \tilde{u}_0$ one 
should not expect that $u_t \leq \tilde{u}_t$ for $t>0$.  
To replace these arguments, the key new method in this paper is
the use of comparison results for the total occupation measure
$\int^{\infty}_0 u_t dt$ and total exit measure $u_{\infty}^{\partial D}$ 
of solutions on a domain $D$. 
Some intuition behind this is given in section \ref{s1.2}.

\noindent  \textbf{4.} 
Using scaling (as in Lemma \ref{scaling})
one can investigate the dependence on other parameters. For example the equation
\[ 
\left\{ 
\begin{array}{l}
\partial_t u = \Delta u + u v - \gamma_0 u +  \sqrt{\sigma u} \, \dot{W},  \\ 
\partial_t v = -uv, \quad v_0 = 1.
\end{array}
\right. 
\]
can be reduced to the standard form (\ref{spde}), by a linear change of variables, with
the parameters becoming $\gamma = \gamma_0 \sigma^{-2/(4-d)}$  
and $\beta= \sigma^{-2/(4-d)}$. The main theorem 
shows that possible life occurs, for example when $d=3$, if
$0 \leq \gamma_0 < \Phi(\sigma)$ for some
$\Phi:(0,\infty) \to (0,1)$ and
satisfying $\lim_{\sigma \to 0} \Phi(\sigma) = 1$ and
$\lim_{\sigma \to \infty} \Phi(\sigma) = 0$ at certain rates.
One reason for choosing $\beta,\gamma$ as our parameters is that we have 
natural monotonicity in these parameters which is unclear for other choices - 
for example we do not know if $\Phi$ is decreasing. 

\noindent \textbf{5.} 
One can investigate, via phase plane methods, the possibility of travelling 
waves for the corresponding deterministic system
$ \partial_t u =  \Delta u + uv -\gamma u, \; \partial_t v = -uv$ with $ v_0 = 1$.
(With no noise term and $\beta >0$ we can remove one more parameter by linear scaling, 
leaving only $\gamma \geq 0$.) We require $\gamma \leq 1$, else the death term exceeds 
the reaction term resulting in death. Fix $k \in \R^d$ with $|k|=1$. Looking for a 
travelling wave solution of the form $u_t(x)=U((x\cdot k)-ct), \, v_t(x)=V((x \cdot k)-ct)$ 
with $c>0$, we set $W=U'$ to obtain the first order system
$$
U' =  W, \quad
V'= \frac{1}{c} UV, \quad
W'= \gamma U-cW-UV.
$$
The steady states are $(U,V,W) = (0,a,0)$ for any $a$. By scaling we assume the untouched
nutrient level is $1$ and seek a path in phase space connecting $(0,0,0)$ to
$(0,1,0)$, so that $u$ is a travelling `hump', behind the wave the nutrient is used up
and ahead it is untouched. Linearizing around $(0,0,0)$ we find the eigenvalues satisfy
$\lambda=0 $ or $\lambda^2 + c \lambda -\gamma=0$. So for any $\gamma>0$ we have
real eigenvalues $\lambda_1 < 0=\lambda_2 < \lambda_3$. Linearizing around $(0,1,0)$ we
obtain $\lambda=0$ or $\lambda^2 + c \lambda + (1 - \gamma) =0$. For
$\gamma < 1$ we get $\lambda_1 < \lambda_2 < 0 = \lambda_3$ when $c^2 \geq 4(1-\gamma)$. 
Complex roots, which imply the impossibility of a (non-negative) travelling wave, occur when 
$c^2 < 4(1-\gamma)$. 
The equations are three dimensional and degenerate and we have not done a 
rigorous analysis. Based on the above local picture, 
the following seems the reasonable best guess: there is a unique 
path connecting $(0,0,0)$ to $(0,1,0)$ which stays in the  good region $U,V>0$
precisely if $c^2 \geq 4(1-\gamma)$ and hence a travelling wave for all
$\gamma <1$ and with speeds greater than or equal to $2 \sqrt{1-\gamma}$.
The main result supports the conjecture that the deterministic travelling waves are 
stable to small enough individual based noise noise in dimensions $d=2,3$. 
\subsection{A description of key techniques and the layout of the paper} \label{s1.2}
In this section we explain some of the intuition and tools used in the proofs,
emphasising the new ideas. 

When $\beta=0$, solutions $u$ to (\ref{spde})
are Dawson-Watanabe processes with underlying Brownian motion and mass annihilation
at rate $\gamma$. In dimensions $d \geq 2$, Dawson-Watanabe processes 
are singular measures $(u_t(dx):t \geq 0)$. 
However, as shown by Sugitani \cite{sugitani},
when $d \leq 3$ the  occupation measures $\int^t_s u_r dr$ 
are absolutely continuous with repsect to Lebesgue measure, 
with densities $u(s,t,x)$. This will remain true for 
for solutions to (\ref{spde}) and hence the the formal solution 
$v_t(x)= v_0(x) \exp(-u(0,t,x))$ will make sense 
as a function. Using this idea, we construct solutions 
to (\ref{spde}) in $d \leq 3$ by a change of measure starting from  
a Dawson-Watanabe process; a reverse argument
yields uniqueness in law of solutions.

The results on possible life use a construction of a supercritical oriented 
percolation process, coupled to a solution of (\ref{spde}), 
so that if the oriented percolation survives so too does 
the solution.  
The use of an embedded oriented percolation
to show survival was first developed for interacting
particle systems (see Durrett \cite{durrett} for many examples).  
Similar methods have been used before by the authors to
establish a one-parameter phase transition for a noisy version of the KPP
equation \cite{mueller+tribe}. The results on certain death in dimensions $d \geq 2$ 
are simpler, and the background idea is a comparison with a subcritical 
branching process (for which simple moment estimates show death).

When possible life occurs, one expects that 
the process lives on a front propagating with
asymptotically linear speed. Noisy oscillations around this linear speed mean 
that it is not possible to predict exactly where in space, at a given time $t$, 
to look for the living part of the solution (that is where it is in 
the process of exploiting the nutrient). This leads to the wish to apply the 
percolation comparison in space at random times.
To overcome this technical problem, we freeze the process 
at random stopping times, by means of exit measures. 
The exit measures can 
be constructed  by solving (\ref{spde}) on a domain $D$ with Dirichlet boundary 
conditions and measuring the total flux of mass out of the boundary of 
the domain.
Suitable comparison results do hold for this model, not 
for $u_t(dx)$ but for the total exit measures $u^{\partial D}_{\infty}(dx)$.
For example if 
$u_0(dx)$ and $\tilde{u}_0(dx)$ are supported inside $D$ and $u_0(dx) \leq \tilde{u}_0(dx)$,
then one can construct solutions for which $u^{\partial D}_{\infty}(dx) 
\leq \tilde{u}^{\partial D}_{\infty}(dx)$. It turns out that the total exit measures
are stochastically monotone both in the initial conditions and in the parameters $\beta,\gamma$.
A spatial Markov property, analogous to that
found for Dawson-Watanabe processes by Dynkin in \cite{dynkin}, shows that
exit measures can be constructed iteratively on larger and larger domains, and 
is used to set up the percolation comparison. It is this spatial Markov property especially, that  
means the use of the special form of noise $\sqrt{u}\dot{W}$ is needed for our methods.

Intuition as to why comparison results hold for the total exit measures 
is easiest for particle approximations to (\ref{spde}), where
population particles may react with nutrient particles leading to
death of a nutrient particle and the creation of extra population particles. 
The total exit measures do not depend on the time at which particles hit the exit 
of the domain $D$. It does not affect the final
level of the nutrient $v_{\infty}(x)$ if we change the order that particles pass
through position $x$. For example, we may think of the initial mass $\tilde{u}_0(dx)$ as being
composed of blue particles, identical to those from $u_0(dx)$ and some
extra red particles represented by $\tilde{u}_0(dx) - u_0(dx)$. We freeze the
extra red particles until all the blue particles have died or are frozen on 
$\partial D $. Then we may run the red particles through whatever is left of the
nutrient, but this can only lead to a larger exit measure. The reader may wish to 
draw a diagram illustrating the possible 
interactions for the simplest non-trivial case of two population 
particles and a single nutrient particle - it was such a diagram that
initiated this work. 

\textbf{Layout.} In section \ref{s2} we state the existence, 
uniqueness and comparison results for occupation and exit measures. 
Section \ref{s3} contains the change of measure arguments
needed for the proof of existence and uniqueness, and 
the proof of the existence of a critical curve $\Psi(\beta)$. 
Approximations to solutions, using a discretized nutrient, are given in section \ref{s4}, together with 
tightness and convergence results.
Comparison results for the approximations are established in section \ref{s5}
together with the weak convergence arguments allowing us to pass to the continuous limit.
However, the arguments for life (in section \ref{s6}) and death (in section \ref{s7})
can be read after the statement of the comparison results in section \ref{s2}. 
\section{Statements of existence, uniqueness and comparison results} \label{s2}
\subsection{Definition of solutions} \label{s2.1}

We need to consider solutions to (\ref{spde}) on domains $D$. Throughout the paper our 
domains are taken to be either
(i) connected open sets with smooth boundary; or (ii) open 
boxes, that is of the form
$D = \prod_{i} (a_i,b_i)$ for $-\infty \leq a_i<b_i \leq \infty$. 
This restriction allows us to quote estimates on the
Green's functions for the domain. Moreover, boxes will be all that we need
for the blocking arguments in sections \ref{s6} and \ref{s7}.
We will look for solutions with zero Dirichlet boundary conditions,
that is where the population is killed on reaching the boundary. This will
be encoded by the test functions for the martingale problem below.  
We let $C^k_b(\overline{D})$ be the space of $k$-times continuously differentiable 
functions $\phi:D \to \R$ for which $\phi$ and its derivatives can be extended to be 
continuous bounded functions on the closure $\overline{D}$.  
Let $C^2_0(\overline{D})$ be the subset of those $\phi \in C^2_b(\overline{D})$ 
for which $\phi =0$ on $\partial D$. 
Note that solutions on $D=\R^d$ are a special case of the definition below, 
and in this case $C^2_0(\overline{D})$ reduces to the space of bounded functions 
with two bounded continuous derivatives. 

\noindent
\textbf{Notation.} We write $\mathcal{M}(D)$ for the space of finite
measures on $D$ with the topology of weak convergence. For $\mu \in \mathcal{M}$ and 
(suitable) measurable functions $f,g:D \to R$ we will write
$\mu(f)$ for $\int_D f(x) \mu(dx)$ and $\langle f,g \rangle$ for 
$\int_D f(x) g(x) dx$.

We look for solutions $(u,v)=(u_t(dx),v_t(x):t \geq 0, x \in D)$ 
to (\ref{spde}) defined on some filtered probability space 
$(\Omega,(\mathcal{F}_t),\mathcal{F}, P)$. The solutions will be in the sense of 
martingale problems. The adapted process $u$ will have 
continuous paths in $\mathcal{M}(D)$. Let $\mathcal{P}$ be the $(\mathcal{F}_t)$ predictable 
sets in $\Omega \times [0,\infty)$ and
$\mathcal{B}(D)$ be the Borel subsets of $D$. 
The process $v$ will be predictable and non-increasing, that is the map $(\omega,t,x) \to v_t(x)(\omega)$ 
will by $\mathcal{P} \times \mathcal{B}(D)$ measurable, and such that the paths  
$t \to v_t(x) \in [0,1]$ are non-increasing for all $x$, almost surely. 
We say that such a pair $(u,v)$ is a solution to (\ref{spde}) on $D$ if  
for $\phi \in C^2_0(\overline{D})$
\begin{eqnarray} 
u_t(\phi) &=& u_0(\phi) + \int^t_0 u_s \left(\Delta \phi + \beta v_s \phi 
- \gamma \phi \right)ds + m_t(\phi), \quad \mbox{$P$ a.s.} \label{MGP1} \\
\langle v_t,\phi \rangle &=& \langle v_0,\phi \rangle - \int^t_0 u_s(v_s \phi) ds, \quad \mbox{$P$ a.s.} 
\label{MGP2}
\end{eqnarray} 
where $m_t(\phi)$ is a continuous $(\mathcal{F}_t)$ local martingale, with 
$m_0(\phi)=0$ $P$ a.s., and with quadratic variation 
\begin{equation} \label{MGP3}
\left[m(\phi)\right]_t = \int^t_0 u_s(\phi^2) ds, \quad \mbox{$P$ a.s.} 
\end{equation}
If, in addition, $P[u_0=\mu, v_0=f]=1$, for some measurable
$f:D \to [0,1]$ and $\mu \in \mathcal{M}(D)$, we say 
the solution has initial condition $(\mu,f)$.

We will wish to keep track of the flux of mass that exits to the boundary $\partial D$ 
by time $t$, which creates the so called exit measure $u_t^{\partial D}(dx)$
on $\partial D$.
For a deterministic equation, this flux would be given by a surface integral
$\int^t_0 \int_{\partial D} \nabla u \cdot \hat{N} dS \, ds$ for an outward unit normal 
$\hat{N}$. A simple way to treat this process in the noisy setting, where these
derivatives do not exist, is to extend the martingale problem formulation to
test functions $\phi \in C^2_b(\overline{D})$. The exit measure will take values in
$\mathcal{M}(\partial D)$, the space of finite measures on $\partial D$ with the topology
of weak convergence. We look for an adapted 
$\mathcal{M}(\partial D)$ valued process $t \to u^{\partial D}_t$, with 
non-decreasing continuous paths, so that $u^{\partial D}_0 =0$ $P$ a.s. and 
for $\phi \in C^2_b(\overline{D})$
\begin{equation} \label{MGPD1+}
u_t(\phi) = u_0(\phi)  
+ \int^t_0 u_s \left( \Delta \phi + \beta v_s \phi - \gamma \phi \right)ds 
- u^{\partial D}_t(\phi) + m_t(\phi), \quad \mbox{$P$ a.s.}
\end{equation}
where $m_t(\phi)$ is a continuous  $(\mathcal{F}_t)$ local martingale, with 
$m_0(\phi)=0$ $P$ a.s., and with quadratic variation as in (\ref{MGP3}).

The following existence and uniqueness result is shown in section \ref{s3.2}, using 
a change of measure argument starting with a Dawson-Watanabe process on $\overline{D}$. 

\noindent
\textbf{Notation.} We let 
$C([0,\infty),E)$ be the space of continuous functions with values in a metric space $E$, 
with the topology of uniform convergence on compacts. Both
$\mathcal{M}(D)$ and $ \mathcal{M}(\partial D)$ are metrizable as Polish spaces.
Let $\Omega_D$ be the space $C([0,\infty), \mathcal{M}(D) \times \mathcal{M}(\partial D))$.
We write $(U_t,U_t^{\partial D})$ for the canonical 
random variables on $\Omega_D$, $\mathcal{U}$ for the Borel subsets of 
$\Omega_D$ and $(\mathcal{U}_t)$ for the canonical filtration. 
Note that in the case $D=\R^d$ where $\partial D= \emptyset$ the second component
of the path space becomes trivial.

\noindent
\textbf{Notation.} Let $\mathcal{B}(D,[0,1])$ be the space of Borel measurable $f:D \to [0,1]$.
We give $\mathcal{B}(D,[0,1])$ the sigma field generated by $ \langle f,\I(A) \rangle$ for 
all bounded Borel $A \subseteq D$.
\begin{Theorem} \label{E!}
Suppose $d \leq 3$ and $D$ is a domain as described above. Fix $\beta,\gamma \geq 0$. 
\begin{enumerate}
\item[(i)] For any $ \mu \in \mathcal{M}(D)$ and $f \in \mathcal{B}(D,[0,1])$, 
there exists a solution $(u,v)$ to (\ref{spde}) on $D$ with initial conditions $(\mu,f)$.
\item[(ii)] For any solution, the law of $(u_t:t \geq 0)$ on 
$C([0,\infty), \mathcal{M}(D))$ is uniquely determined 
by the law of $(u_0,v_0)$.  
\item[(iii)] For any solution and $0 \leq s \leq t$, the occupation measure
$\int^t_s u_r dr$ is absolutely continuous with respect to Lebesgue measure, 
and there is a continuous
version of the density $(s,t,x) \to u(s,t,x)$ on $0<s \leq t, x \in D$, almost surely. 
Moreover, setting $u(0,t,x) = \lim_{s \downarrow 0} u(s,t,x)$ we have, 
\[
v_t(x) = v_0(x) e^{-u(0,t,x)} \quad \mbox{for all $t > 0$, for almost all
$x \in D$, $P$ a.s.}
\]
\item[(iv)] For any solution $(u,v)$ there exists an adapted exit measure process
$u^{\partial D} = (u^{\partial D}_t:t \geq 0)$, with non-decreasing continuous paths 
in $\mathcal{M}(\partial D)$, satisfying the extended martingale problem
(\ref{MGPD1+}).
Moreover, for all $t \geq 0$, there exists a measurable map 
$R_t:C([0,t],\mathcal{M}(D)) \to \mathcal{M}(\partial D)$ so that
$u^{\partial D}_t = R_t((u_s: s \leq t))$ almost surely. 
\item[(v)] Let $Q^{D,\beta,\gamma}_{\mu,f}$ be the law of $(u,u^{\partial D})$ on
$\Omega_D$ for a solution with initial condition $(\mu,f)$. Then
$(\mu,f) \to Q^{D,\beta,\gamma}_{\mu,f}(B)$ is measurable for any Borel 
$B \subseteq  \Omega_D$. Moreover,
the family of laws  $(Q^{D,\beta,\gamma}_{\mu,f}: \mu \in \mathcal{M}(D), 
f \in \mathcal{B}(D,[0,1]))$ satisfies the following strong Markov property: for any solution
$(u,v)$ defined on $(\Omega,(\mathcal{F}_t),\mathcal{F},P)$, and any finite $(\mathcal{F}_t)$
stopping time $\tau$,
\[
E\left[ h(u_{\tau+\cdot}) | \mathcal{F}_{\tau} \right]
= Q_{u_{\tau}, v_0 \exp(-u(0,\tau))}^{D,\beta,\gamma} 
\left[ h(U_{\cdot}) \right] \qquad
\mbox{$P$ a.s,}
\] 
for all bounded measurable $h:C([0,\infty),\mathcal{M}(D)) \to \R$.
\end{enumerate}
\end{Theorem}
We note one consequence of these results: 
using part (iii), whenever $(u,v)$ is a solution we may replace
$v$ by the process $v_0(x) \exp(-u(0,t,x))$ and it will still be a solution.

The following scaling lemma will be very useful in allowing us to transfer 
large or small parameters to terms in the equations that are convenient. 
We find our intuition stronger when working in regimes with small or large 
values of parameters, rather than small or large space and time scales. 
This lemma allows us to transfer between these two regimes. 
\begin{Lemma} \label{scaling} 
Suppose that $(u,v)$ is a solution to (\ref{spde}) on $D$. 
For $a,b,c,e>0$ define 
\[
\tilde{u}_t(A) = \frac{a}{c^d} u_{bt}(cA), \quad \tilde{v}_t(x)=
e v_{bt}(cx)
\]
where $cA= \{cx: x \in A\}$. Then $(\tilde{u}, \tilde{v})$ is a
solution on $c^{-1}D$, in the sense described above, to the equation 
\[
\left\{ 
\begin{array}{l}
\partial_t \tilde{u} = \frac{b}{c^2} \Delta \tilde{u} + \frac{b \beta}{e} 
\tilde{u} \, \tilde{v} - b \gamma \tilde{u} + \sqrt{\frac{a b}{c^d} \tilde{u} }\, \dot{W}, \\ 
\partial_t \tilde{v} = -\frac{b}{a} \tilde{u} \, \tilde{v}.
\end{array}
\right. 
\]
Furthermore, for $A \subseteq \partial (c^{-1}D)$, we have $ \tilde{u}_t^{\partial (c^{-1} D)}(A) = (a/c^d) 
u_{bt}^{\partial D}(cA)$. 
\end{Lemma}
The equation for $(\tilde{u},\tilde{v})$ follows immediately by scaling the martingale problem (\ref{MGP1},\ref{MGP2}).
The scaling of the exit measures follows by matching the finite variation parts of the semimartingale decomposition
of $\tilde{u}_t(\phi)$ and of $(a/c^d) u_{bt}(\phi(\cdot/c))$ for $\phi \in C^2_b(c^{-1} \overline{D})$ given in the 
corresponding extended martingale problems.  
\subsection{Decomposition results for the total occupation and exit measures} \label{s2.2}
\textbf{Notation.} 
For a solution $(u,v)$ to (\ref{spde}) we write $u_{[s,t]}$ for the occupation measure
$\int^t_s u_r dr$ (so that $u_{[s,t]}(dx) = u(s,t,x)dx$ $P$ a.s.) and write 
\[
u_{[0,\infty)} = \lim_{t \to \infty} u_{[0,t]}, \quad
u(0,\infty,x) = \lim_{t \to \infty} u(0,t,x),
 \quad \mbox{and} \quad 
u^{\partial D}_{\infty} = \lim_{t \to \infty}
u^{\partial D}_{t}
\]
for the total occupation measure on $D$ and its density, and for the 
total exit measure on $\partial D$. 

The lemmas below are stated under the restriction that $D$ is bounded, which implies
(see section \ref{s3.3}) that the total occupation and exit measures are almost surely finite.
The first lemma describes the monotonicity of the total occupation and exit measures 
with respect to the initial conditions $(\mu,f)$. We write $\equalD$ for equality in 
distribution.
\begin{Lemma} \label{comparison1} 
Fix a bounded domain $D$ and 
initial conditions $\mu=\mu^- + \mu^+$ and $f = f^- + f^+$. 
The law of the total occupation and exit measures of a solution $(u, v)$ to
(\ref{spde}) on $D$ with initial conditions $(\mu,f)$ can be decomposed as 
\begin{equation} \label{decomposition}
\left( u_{[0,\infty)}, u^{\partial D}_{\infty} \right)
\equalD   \left( u^-_{[0,\infty)} + u^+_{[0,\infty)}, 
u^{\partial D,-}_{\infty}+ u^{\partial D,+}_{\infty} \right) 
\end{equation} 
in either of the following two ways:

\noindent
\textbf{(i)} $(u^-,v^-)$ is a solution with initial 
conditions $(\mu^-,f)$ and, conditional on
$ \sigma\{u^-\}$, the process $(u^+,v^+)$ is a solution with initial conditions 
$ u^+_0 = \mu^+ $ and $v^+_0= v^-_{\infty}$ (where $v^-_{\infty} = \lim_{t \to \infty} v^-_t$);

\noindent
\textbf{(ii)} $(u^-,v^-)$ is a solution with initial conditions $(\mu,f^-)$ and, conditional on
$\sigma\{u^-\}$, the process $(u^+,v^+)$ is a solution 
with initial conditions $ u^+_0 = \beta f^+ (1-\exp(-u^-(0,\infty)))\,dx $ and 
$v^+_0= f \exp(-u^-(0,\infty))$.
\end{Lemma}
The next lemma shows monotonicity in the parameters $\beta$, $\gamma$.
\begin{Lemma} \label{comparison2}
Fix a bounded domain $D$, and parameter 
values $\beta= \beta^- + \beta^+$ and $\gamma^- = \gamma + \gamma^+$.
The law of the total occupation and exit measures of a solution $(u, v)$ to
(\ref{spde}) on a bounded domain $D$ with initial conditions $(\mu,f)$ 
and parameter values $(\beta,\gamma)$ can be decomposed as in
(\ref{decomposition}) in either of the following two ways:

\noindent
\textbf{(i)} $(u^-,v^-)$ is a solution with initial conditions $(\mu,f)$ and with
parameter values $\beta^-,\gamma$ and, conditional on $\sigma\{u^-\}$, the process $(u^+,v^+)$ is a solution 
with initial conditions $u^+_0 = \beta^+ (f-v^-_{\infty}) dx $ 
and $v^+_0= v^-_{\infty}$ and with parameter values $\beta,\gamma$;

\noindent
\textbf{(ii)} $(u^-,v^-)$ is a solution to (\ref{spde}) 
with initial conditions $(\mu,f)$ and parameter values $\beta,\gamma^-$ and, conditional on
$\sigma\{u^-,v^-\}$, the process $(u^+,v^+)$ is a solution
with initial conditions $u^+_0 = \gamma^+ u^-_{[0,\infty)}$
and $v^+_0= v^-_{\infty}$ and with parameter values $\beta,\gamma$.
\end{Lemma}
A special case of the Lemma \ref{comparison2} (i) will be particularly useful. The total exit 
measure can be built up in two stages: first run a process with $\beta=0$ (which is the
well understood Dawson-Watanabe process) and then run a second solution
started with the mass that the nutrient would have produced from the first process.

The following simpler comparison will often be useful.
It states, roughly, that if we convert some of the nutrient available 
into the equivalent amount of initial mass at time zero then we obtain a healthier process,
in the sense of stochastic ordering.  
\begin{Lemma} \label{comparison3}
Fix a bounded domain $D$ and initial conditions $\mu$ and $f,g$.
Let $(u,v)$ and $(\tilde{u},\tilde{v})$ be solutions to (\ref{spde}) on $D$ 
with  initial conditions $(\mu,f+g)$ and $(\mu + \beta g \, dx,f)$ respectively.
Then for any bounded measurable $F: \mathcal{M}(D) \times \mathcal{M}(\partial D) \to \R$
that is non-decreasing in both variables 
$$
E[F(u_{[0,\infty)}, u^{\partial D}_{\infty})] \leq 
E[F(\tilde{u}_{[0,\infty)}, \tilde{u}^{\partial D}_{\infty})].
$$
\end{Lemma}
The final lemma is a spatial Markov property, analogous to that known for 
Dawson-Watanabe processes. We choose two domains $D^- \subseteq D^+$. We allow the boundaries
$\partial D^-$ and $\partial D^+$ to intersect, but we need this intersection not to be 
too complicated, and for this we ask that
\begin{equation}
\mbox{$(\partial D^- \cap \partial D^+) \cap \overline{(\partial D^- \setminus \partial D^+)}$
has surface measure zero in $\partial D^-$.}
\label{nullset}
\end{equation}
This is satisfied for example if both domains are boxes.
We write $\mu |_A$ for the restriction of a measure $\mu$ to a set $A$.     
\begin{Lemma} \label{sMp}
Fix bounded domains $D^- \subseteq D^+$ satisfying (\ref{nullset}).
Suppose $(u,v)$ is a solution to (\ref{spde}) on $D^+$ started at $(\mu,f)$. 
Then the law of the total occupation and exit measures satisfies
\[ 
\left( u_{[0,\infty)}, u^{\partial D^+}_{\infty} \right)
\equalD   \left( u^-_{[0,\infty)} + u^+_{[0,\infty)}, 
u^{\partial D^-,-}_{\infty} |_{\partial D^+} + u^{\partial D^+,+}_{\infty} \right) 
\]
where $(u^-,v^-)$ is a solution to (\ref{spde}) 
on $D^-$ with initial conditions $(\mu \I(D^-), f \I(D^-))$ and, 
conditional on $ \sigma\{u^-\}$, the process $(u^+,v^+)$ is a
solution on $D^+$ with initial conditions
\[ 
u^+_0 = u^{\partial D^-,-}_{\infty}|_{D^+} + \mu \I(D^+ \setminus D^-) \quad \mbox{and} \quad 
v^+_0 = v^-_{\infty} \I(D^-) + f \I(D^+ \setminus D^-).
\] 
\end{Lemma}
\textbf{Remark.} As indicated in the introduction, these lemmas will be proved
for discretized approximations to solutions, where the nutrient is used up in finitely many 
steps, and then via a passage to the limit. We spent considerable energy trying to find
a proof using the continuum equation (\ref{spde}) directly. Approaches using a change of measure
(since the lemmas are easy to prove for a Dawson-Watanabe process) typically failed because the
decomposition (\ref{decomposition}) holds only for the total occupation and exit measures
and fails for these measures at times $t<\infty$. Another approach, using Laplace functionals, 
is sketched in remark (b) of section \ref{s3.1}. On the other hand, we found the discretized approximations,
in particular particle pictures, were the easiest testing ground for exploring possible 
extensions of these results to other models.
\section{Change of measure arguments} \label{s3}
In section \ref{s3.1} we review some results on Dawson-Watanabe processes 
and indicate when these imply estimates for solutions of (\ref{spde}). Section \ref{s3.2}
gives the proof of Theorem \ref{E!}, and section \ref{s3.3} contains the proof of the 
existence of the critical curve $\Psi(\beta)$.
\subsection{Results from Dawson-Watanabe processes} \label{s3.1}
A Dawson-Watanabe process, with underlying Brownian spatial motion, killed at the exit of
a domain $D$, and with a (deterministic) bounded, measurable mass annihilation/creation rate
$\eta: D \to \R$, is the solution on $D$, in the same sense as in section \ref{s2.1}, 
to the equation  
\begin{equation}  \label{DWSPDE}
\partial_t u =  \Delta u  - \eta u + \sqrt{u} \, \dot{W}.
\end{equation}
We shall abbreviate Dawson-Watanabe processes as DW processes.
In the case $\eta(x) = \gamma \in \R$ is constant we will call the process
a DW($D,\gamma$) process. We describe below some results for 
DW processes that we shall need, and we mostly indicate where to find the proofs, 
or how to adapt them, from the literature. The surveys by Perkins 
\cite{perkins1} and Dawson \cite{dawson} are our main sources.
However, we give proofs for some of the results on exit measures, that are perhaps not
in the literature, at the end of this subsection. 

A number of the estimates (for example moment bounds, estimates on the speed
of the support, estimates on the extinction probability) also hold for solutions
to (\ref{spde}). One way to see this is to
set up a pathwise coupling for $(u,v)$ a solution to (\ref{spde}) 
on $D$ so that $u^{(1)}_t \leq u_t \leq u^{(2)}_t$ where 
$u^{(1)}$ is a DW($D,\gamma$) process and $u^{(2)}$ is a DW($D,\gamma-\beta$)
process. We do not go through the steps to 
construct this coupling since, as we indicate below, the proofs of the desired estimates 
are typically based on stochastic calculus and are easily adapted to follow from
the martingale problem for solutions to (\ref{spde}). 

\noindent \textbf{1. Martingale problems.} We quickly repeat here 
the martingale problem definition for solutions. A continuous, adapted
$\mathcal{M}(D)$ valued process $u$, on some
filtered space $(\Omega, (\mathcal{F}_t), \mathcal{F},P)$,  
is a solution to (\ref{DWSPDE}) if for $\phi \in C^2_0(\overline{D})$
\begin{equation} \label{DWMGPD1+}
u_t(\phi) = u_0(\phi) + \int^t_0 u_s \left(\Delta \phi  - \eta \phi \right)ds 
+ m_t(\phi), \quad \mbox{$P$ a.s.}
\end{equation}
where $m_t(\phi)$ is a continuous $(\mathcal{F}_t)$ local martingale with quadratic variation 
$[ m(\phi)]_t = \int^t_0 u_s(\phi^2) ds$ almost surely.
Solutions, starting from $\mu \in \mathcal{M}(D)$, exist and are unique in law on 
$C([0,\infty), \mathcal{M}(D))$. Furthermore, for any solution 
there will exist an adapted process $(u^{\partial D}_t)$ 
with non-decreasing continuous paths in $ \mathcal{M}(\partial D)$,  
called the exit measure process, so that for $\phi \in C^2_b(\overline{D})$
\begin{equation} \label{DWMGPD1++}
u_t(\phi) = u_0(\phi) + 
\int^t_0 u_s \left( \Delta \phi  - \eta \phi \right)ds 
- u^{\partial D}_t(\phi) + m_t(\phi), \quad \mbox{$P$ a.s.} 
\end{equation}
where $m_t(\phi)$ is a continuous local martingale with variation as before.
The almost sure linearity $m_t(\phi+\psi) = m_t(\phi)+m_t(\psi)$, implies that if
$\phi_n \in C^2_0(\overline{D})$ converge bounded pointwise to 
$\phi \in C^2_b(\overline{D})$ then $m_t(\phi_n) \to m_t(\phi)$ 
in probability (since $[m(\phi-\phi_n)]_t \to 0$).
Now (\ref{DWMGPD1+}) and (\ref{DWMGPD1++}) imply that 
$u_t^{\partial D}(\phi)$ is measurable function of the path $(u_s:s \leq t)$. 
Choose is a countable family $(\phi_k)_{k=1,\ldots}$ in $C^2_b(\overline{D})$ so that the map
$\mu \to (\mu(\phi_k))_{k=1,\ldots}$ is a continuous injection from 
$\mathcal{M}(\partial D)$ into $\R^{\infty}$ (with the topology of uniform convergence on compacts). 
Using the measurable inverse of this map, one one can construct a measurable map 
$R_t:C([0,t],\mathcal{M}(D)) \to \mathcal{M}(\partial D)$ so that
\begin{equation} \label{Rmaps}
u^{\partial D}_t = R_t(\{u_s: s \leq t\}) \quad \mbox{$P$ a.s.}
\end{equation}
These maps show that the law of $(u,u^{\partial D})$ on $\Omega_D$ is determined by the law of
$u$ on $C([0,\infty),\mathcal{M}(D))$. Also, for any DW process on $D$, 
the exit measure process can be constructed by using the maps $R_t$ and 
then regularizing the path over the rationals.

\noindent \textbf{Notation.}
We denote by $Q^{D,\eta}_{\mu}$ the law on $\Omega_D$ of a DW process and its exit measure 
$(u,u^{\partial D})$ on $D$, with mass annihilation/creation rate $\eta$ as above, 
and starting at $\mu \in \mathcal{M}(D)$.
We write $Q^{D,\gamma}_{\mu}$ in the case that  $\eta(x)=\gamma$.

\noindent The map $\mu \to Q^{D,\eta}_{\mu}(B)$ is measurable for any
Borel $B \subseteq \Omega_D$; furthermore the family of laws 
$(Q^{D,\eta}_{\mu}: \mu \in \mathcal{M}(D))$ form a strong Markov family
in that for any solution
$u$ defined on $(\Omega,(\mathcal{F}_t),\mathcal{F},P)$, and any finite $(\mathcal{F}_t)$
stopping time $\tau$,
\[
E\left[ h(u_{\tau+\cdot}) | \mathcal{F}_{\tau} \right]
= Q_{u_{\tau}}^{D,\eta} \left[ h(U_{\cdot}) \right] \qquad
\mbox{$P$ a.s,}
\] 
for all bounded measurable $h:C([0,\infty),\mathcal{M}(D)) \to \R$. The measurability
and strong Markov property follow from the Laplace functional 
(\ref{LFD1}) below by, for example, the arguments used in the proof
of Theorem II.5.1 of \cite{perkins1}. 

For a DW process $u$, one can extend the local martingales to a
local martingale measure, which we write as  $m_t(\phi) = \int^t_0 \int_D \phi_s m(dx,ds)$, and is
defined for $\phi:[0,\infty) \times \Omega \times D \to \R$ that are $\mathcal{P} \times \mathcal{B}(D)$
measurable and satisfy $\int^t_0 u_s(\phi^2_s) ds < \infty$ for all $t \geq 0$. The integral
produces a continuous $(\mathcal{F}_t)$ local martingale with quadratic 
variation $[m(\phi)]_t = \int^t_0 u_s(\phi_s^2) ds$.
This follows as in \cite{perkins1} Proposition II.5.4. The same construction applies to
give a local martingale measure for solutions to (\ref{spde}), and this will be needed
in our change of measure arguments. One simple use of these extensions is to allow
time dependent test functions in the martingale problem. Since the exit measures $u^{\partial D}_t$ are non-decreasing, they 
induce a unique measure $du^{\partial D}_t(dx)$ 
on $[0,\infty) \times \partial D$ characterized by the measure of the rectangles
$(s,t] \times A$ being $u^{\partial D}_t(A)-u_s^{\partial D}(A)$. We will write the integral
of $\psi:[0,t] \times \partial D \to \R$ with respect to this measure as
$\int^t_0 du^{\partial D}_s(\psi_s)$. Then, if $\phi:[0,T] \times \overline{D} \to \R$
is in the class $C^{1,2}_b([0,T] \times \overline{D})$, where $\phi$ and one partial 
derivative in time (denoted $\dot{\phi}$) and two partial derivatives in space 
exist and have continuous bounded extensions to $[0,T] \times \overline{D}$, 
the following decomposition holds:
\begin{equation} \label{DWMGPD1+++}
u_t(\phi_t) = u_0(\phi_0) + 
\int^t_0 u_s \left( \Delta \phi_s + \dot{\phi}_s  - \eta \phi_s \right)ds 
- \int^t_0 du_s^{\partial D}(\phi_s)
+ m_t(\phi) 
\quad \mbox{$P$ a.s.} 
\end{equation}
This follows using easier versions of the arguments in \cite{perkins1} Proposition II.5.7. Moreover, the 
analogous extension works for solutions to (\ref{spde}).  

\noindent \textbf{2. Laplace functionals.}
Suppose first that $D$ has a smooth boundary and that $\eta: \overline{D} \to \R$ is smooth and bounded. 
The Laplace functional of the solution, the occupation measure and the exit measure is 
given, for smooth bounded $h^{(1)}: \overline{D} \to \R$, $h^{(2)}:[0,t] \times \overline{D} \to \R$ and
$h^{(3)}:[0,t] \times \partial D \to \R$ by 
\begin{equation} \label{LFD1}
Q^{D,\eta}_{\mu} \left[ e^{- U_t(h^{(1)}) - \int^t_0 U_s(h^{(2)}_s) - \int^t_0 dU_s^{\partial D}(h^{(3)}_s)} \right] 
= e^{ - \mu(\phi_t)}, 
\end{equation}
where $(\phi_s(x): s \in [0,t], \, x \in D)$ is the unique solution, smooth on $[0,t] \times \overline{D}$, to
\begin{equation} \label{LFD1+}
\left\{ \begin{array}{l}
\partial_s \phi_s = \Delta \phi_s - \eta \phi_s - \frac{\phi^2_s}{2} + h^{(2)}_{t-s} \quad \mbox{on $D$}, \\
\quad \mbox{$\phi_s = h^{(3)}_{t-s}$ on $\partial D$ for $s \in [0,t], \quad$ and $\phi_0 = h^{(1)}$ on $D$.} 
\end{array} \right.
\end{equation}
This follows by checking that 
$s \to \exp(-U_s(\phi_{t-s}) - \int^t_0 U_r(h^{(2)}_s)dr - \int^t_0 dU^{\partial D}_r(h^{(3)}_r))$ 
is a martingale for $s \in [0,t]$. 
If in addition $\eta \geq 0$ and $D$ is bounded, the total occupation 
$U_{[0,\infty)} = \lim_{t \uparrow \infty} U_{[0,t]}$ and total exit measures 
$U^{\partial D}_{\infty}= \lim_{t \uparrow \infty} U_t^{\partial D}$ are finite almost surely 
and have Laplace functional
\begin{equation} \label{LFD2}
Q^{D,\eta}_{\mu} \left[ e^{-U_{[0,\infty)}(h_1)- U^{\partial D}_{\infty}(h_2)} 
\right] = e^{-\mu(\phi)}, 
\end{equation}
for smooth bounded $h_1:\overline{D} \to [0,\infty)$ and $h_2:\partial D \to [0,\infty)$, where 
$(\phi(x):x \in D)$ is the unique solution, smooth on $\overline{D}$, to 
\begin{equation} \label{alfeqn}
\Delta \phi = \frac{\phi^2}{2} + \eta \phi - h_1 \quad \mbox{on $D$}, \quad
\mbox{$\phi = h_2$ on $\partial D$}.
\end{equation}
This follows by examining the martingale 
$t \to \exp(-U_t(\phi) - \int^t_0 U_s(h_1) - U^{\partial D}_t(h_2))$.
On bounded domains, and when $\eta \geq 0$, there exist finite solutions to (\ref{alfeqn}) when $h_1,h_2$ are negative and sufficiently small. 
In particular, the exponential moment $Q^{D,\eta}_{\mu}[ \exp(+\theta U^{\partial D}_{\infty}(1) + \theta U_{[0,\infty})(1))]$ is finite
for sufficiently small $\theta>0$. 

When the domain $D$ is a box the associated p.d.e.'s have unique solutions that are smooth inside 
$D$ and continuous on $\overline{D}$, and the formulae for the Laplace functionals still hold, and 
can be established by an approximation argument. Indeed, one concrete way to do this is to suppose,
without loss, that $D$ contains the origin and to replace the test functions above by
$\phi_t(x/(1+\epsilon))$ and $\phi(x/(1+\epsilon)$, which are smooth in $\overline{D}$. 

\noindent
Remark (a). From the pde viewpoint, it is natural to use  
(\ref{DWMGPD1++}) as a means of characterizing the exit measures 
$t \to u^{\partial D}_t(dx)$, that is purely as exit fluxes (as
is suggested in \cite{perkins2}).
Although we do not need such a characterization, since we have other constructions of the exit measure, 
we briefly sketch such an approach, and for convenience we consider deterministic initial 
condition $\mu$ and smooth bounded domain $D$.
Construct the martingale measure $m(dx,ds)$ from a DW($D,\gamma$) process. Then for 
smooth $h: \partial D \to \R$ we may uniquely find $\phi^h \in C^2_b(\overline{D})$ solving 
the Dirichlet problem $\Delta \phi^h =0$ on $D$ and $\phi^h=h$ on $\partial D$.
Use the formula (\ref{DWMGPD1++}) with the test function $\phi^h$ 
to define a continuous path $t \to u^{\partial D}_t(h)$ (although
we have yet to establish that the values $u^{\partial D}_t(h)$ arise from a measure $u^{\partial D}_t$ on 
$\partial D$).
By adding the decompositions (\ref{DWMGPD1++}) for $\phi^g$ and $\phi^h$ we find that
$u^{\partial D}_t(g+h) = u^{\partial D}_t(g) + u^{\partial D}_t(h)$ almost surely. 
Since $\|\phi^h\|_{\infty} \leq \|h\|_{\infty}$, we can deduce a first moment bound
$ \sup_{s \leq t} E[|u^{\partial D}_s(h)|] \leq C(\mu,t,\gamma) \|h\|_{\infty}$ from the corresponding bounds
for $u$. For $\phi \in C^2_b(\overline{D})$ with $\phi=h$ on $\partial D$, the decomposition
(\ref{DWMGPD1++}) holds by adding the decompositions for $\phi^h$ and $\phi-\phi^h$.
Then the first moment bounds allow the decompositions to be extended to 
time dependent test functions $\phi(t,x)$ as above. In particular we find,
for $h \geq 0$, that $ E [ \exp(-\lambda u^{\partial D}_t(h)) ] \leq \exp(- \mu(\phi_{t}))$
where $\partial_t \phi = \Delta \phi - \gamma \phi - \frac{\phi^2}{2}$ on $D$, and 
$\phi_t = \lambda h$ on $\partial D$, and $ \phi_0 = 0$. The inequality here arises 
since we do not yet know that $u^{\partial D}_t(h) \geq 0$ and so we have to localize 
via the stopping times $\tau_n= \inf \{t: u_t(1) \geq n\}$, and we obtain an upper bound by passing
to the limit using Fatou's Lemma. However letting $\lambda \to \infty$ now implies that
$u^{\partial D}_t(h) \geq 0$ almost surely. 
It is then not hard to establish the existence of a measure $u^{\partial D}_t$ so that
$u^{\partial D}_t(h) = \int_{\partial D} h \, du^{\partial D}_t$. 
Finally the Markov property implies that
$u^{\partial D}_t - u^{\partial D}_s$ is non-negative and hence the paths of 
$u^{\partial D}_t$ are non-decreasing.

\noindent 
Remark (b). Mimicking the calculus that leads to the Laplace functional above, one 
can show that for solutions to (\ref{spde}) that
\[
\exp \left(-u_t(\phi) - u_{[0,t]}(h_1) - u^{\partial D}_t(h_2)+ \beta \int^t_0 u_s(v_s \phi) ds \right)
\]
is a martingale when $\phi$ solves (\ref{alfeqn}) with $\eta = \gamma$.
One deduces, by letting $t \to \infty$, that
\[
E \left[ e^{-u_{[0,\infty)}(h_2) -u^{\partial D}_{\infty}(h_3)
+ \beta \langle \phi f, 1-\exp(-u(0,\infty)) \rangle} \right] = e^{-\mu(\phi)}. 
\]
We do not know if this formula, as one
varies $h_1,h_2$, characterizes the law of $(u_{[0,\infty)},u^{\partial D}_{\infty})$.
If it did characterize the law, a simple proof of all the 
decomposition lemmas from section \ref{s2.2} would result, since it is straightforward to
use calculus for both parts of the decompositions to yield the above formula
for the sum. 

\noindent \textbf{3. Some moments.}
Estimates on moments for the total mass process $u_t(1)$, starting from a deterministic initial condition, follow 
by choosing $\phi=1$ in (\ref{DWMGPD1++}), using $u^{\partial D}_t(1) \geq 0$, by the usual arguments 
(localizing by the stopping times $\tau_n = \inf \{t: u_t(1) \geq n\}$ and applying Gronwall's inequality). 
By rewriting (\ref{DWMGPD1++}), again with $\phi=1$, as an equation for $u^{\partial D}_t(1)$, one can then 
deduce moments for the exit measures. These arguments show, for instance, that for for $p \geq 1$ 
\begin{equation} \label{totalmass1}
Q^{D,\gamma}_{\mu} \left[ \sup_{t \leq T} (U_t(1))^p  + |U^{\partial D}_T(1)|^p \right] 
\leq C(p,T,\gamma) \left( (\mu(1) + (\mu(1))^p \right). 
\end{equation}
In particular, starting from a deterministic initial condition the processes 
$m_t(\phi)$ are true martingales. The same arguments apply to solutions to
(\ref{spde}), and lead to 
\begin{equation}
Q^{D,\beta,\gamma}_{\mu,f} \left[ \sup_{t \leq T} (U_t(1))^p 
+ |U^{\partial D}_T(1)|^p \right]
\leq Q^{D,\gamma - \beta}_{\mu} \left[ \sup_{t \leq T} (U_t(1))^p 
+ |U^{\partial D}_T(1)|^p \right].
\label{totalmass2}
\end{equation}
(We abuse logical order here, and below, by writing $Q^{D,\beta,\gamma}_{\mu,f}$ before
checking that these laws are uniquely defined. In such cases we mean that the arguments apply to 
all solutions to (\ref{spde}).)

We will need the first and second moment formulae, for measurable $g,h \geq 0$, 
\begin{eqnarray} 
Q^{D,\gamma}_{\mu} \left[ U_t(g) \right] & = & \mu(G^{D,\gamma}_t g), \nonumber  \\
Q^{D,\gamma}_{\mu} \left[ U_t(g) U_s(h) \right] & = &  
 \mu(G^{D,\gamma}_t g)  \mu(G^{D,\gamma}_s h)  + 
\int^{t \wedge s}_0 \mu(G^{D,\gamma}_r(G^{D,\gamma}_{t-r}g \, G^{D,\gamma}_{s-r}h))dr, 
\label{2ndmoments}
\end{eqnarray}
where $G^{D,\gamma}_t g(x) = \int_D G^{D,\gamma,x}_t(y) g(y) dy$ and $G^{D,\gamma,x}_t(y)$ is the killed
Green's function $G^{D,\gamma,x}_t(y)= e^{-\gamma t} G^{D,x}_t(y)$ 
for the domain $D$ (that is $G^{D,x}_t(y)$ is the transition density of a 
Brownian motion killed on its exit from $D$). 
These can be proved by using the time dependent test functions $\phi^{(g)}$ and $\phi^{(h)}$, where
$\phi^{(g)}_s = G^{D,\gamma}_{t-s} g$ for $s \in [0,t]$, in (\ref{DWMGPD1+++}). 
For smooth $g,h$ compactly supported in $D$, this test function has the required regularity,
and the moment formulae can then be extended to general $g,h$ 
by monotone class methods. In particular the same
methods apply to the martingale problem for solutions to (\ref{spde}) 
and give the bounds, for $g \geq 0$ and $k=1,2$, 
\begin{equation} \label{12momentbounds}
Q^{D,\gamma}_{\mu} \left[(U_t(g))^k \right]
\leq Q^{D,\beta,\gamma}_{\mu,f} \left[(U_t(g))^k \right] \leq 
Q^{D,\gamma-\beta}_{\mu} \left[(U_t(g))^k \right],
\end{equation} 
illustrating the natural intuition that less killing leads to larger solutions. 

\noindent 
\textbf{4. Extinction probabilities and speed of propagation.} 
The probability of extinction for DW($\R^d,\gamma$) processes is given by 
$Q^{\R^d,\gamma}_{\mu} [U_t= 0] = \exp(- \lambda^{(\gamma)}_t \mu(1))$
where $\lambda^{(\gamma)}_s$ is determined for $s \in (0,t]$ by
\[
\dot{\lambda}^{(\gamma)} = - \frac12 (\lambda^{(\gamma)})^2 - \gamma \lambda^{(\gamma)},
\qquad \mbox{and $\quad \lambda^{(\gamma)}_s \to \infty$ as $s \downarrow 0$}. 
\]
This follows by examining the martingale 
$s \to \exp(-\lambda^{(\gamma)}_{t-s} U_s(1))$ for $s \in (0,t]$ (see the derivation of 
(\ref{fixedtsupport}) below). This same argument applies to solutions of (\ref{spde}) to give the bounds
\begin{equation} \label{extinctionprob}
Q^{\R^d,\gamma-\beta}_{\mu}[U_t=0] \leq 
Q^{D,\beta,\gamma}_{\mu,f}[U_t=0] \quad \mbox{and} \quad
Q^{\R^d,\beta,\gamma}_{\mu,f}[U_t=0] \leq Q^{\R^d,\gamma}_{\mu}[U_t=0].
\end{equation}

DW processes, started from finite measures, are compactly
supported at all times $t>0$. The same will
apply for solutions $u$ to (\ref{spde}), by the 
absolute continuity of the law of $u$ with respect to a DW process (see section
\ref{s3.2}). However we need a quantitative bound on the size of the support and 
we use two estimates for this behaviour, 
based on the results from Dawson, Iscoe and Perkins \cite{dip}.
For $R>0$ let $D_R = (-R,R)^d$. Choose $0 \leq \xi_R \leq 1$ smooth and satisfying $\{\xi_R>0\} = 
\overline{D}^c_{R}$.  Then for $M \geq 0$ and smooth bounded $\eta:\R \to \R$, 
let $\phi_s = \phi^{(M,R,\eta)}_s$ be the unique non-negative solution to
\[
\dot{\phi} =  \Delta \phi - \frac12 \phi^2 - \eta \phi,
\quad \psi_0 = M \xi_R.
\]
Then, for $s>0$, $\phi^{(M,R,\eta)}_s \uparrow \phi^{(R,\eta)}_s < \infty$
as $M \uparrow \infty$ and moreover $\phi^{(R,\eta)}_s(x) \to 0$ as 
$R \to \infty$ (use a comparsion argument and \cite{dip} Lemma 3.6). Taking expectations of the 
martingale $\exp(-U_s(\phi^{(M,R,\eta)}_{t-s}))$ at times $s=0$ and $s=t$, and then letting $M \to \infty$, 
one finds that $ Q^{\R^d,\eta}_{\mu} [ U_t(\overline{D}^c_{R}) = 0 ] 
= \exp( - \mu( \phi^{(R,\eta)}_t))$. The same proof shows, when 
$\gamma - \beta f \geq \eta$, the bound 
\begin{equation} \label{fixedtsupport}
Q^{\R^d,\beta,\gamma}_{\mu,f} \left[U_t(\overline{D}^c_{R})=0 \right] 
\geq Q^{\R^d,\eta}_{\mu} \left[U_t(\overline{D}^c_{R})=0 \right].
\end{equation}
Letting $R \to \infty$ shows that solutions to (\ref{spde})
started from compactly supported $\mu \in \mathcal{M}(\R^d)$ must 
have compact support at any time $t>0$. 

Similarly, let $\psi_s = \psi^{(M,R,\eta)}_s$ be the unique 
non-negative solution to
\[
\dot{\psi} = \Delta \psi - \frac12 \psi^2 - \eta \psi
+ M \xi_R, \quad \psi_0 =0.
\]
Then, for $|x| < R$, $\psi^{(M,R,\eta)}_s(x) \uparrow \psi^{(R,\eta)}_s(x) < \infty$
as $ M \uparrow \infty$ and moreover $\psi^{(2,R,\gamma)}_s(x) \to 0$ as 
$R \to \infty$ (argue as in \cite{dip} Theorem 3.3 (ii)). 
Examining the process $\exp(-U_s(\psi^{(M,R,\eta)}_{t-s}))$ as above one finds that 
$Q^{\R^d,\eta}_{\mu} [ U_{[0,t]}(\overline{D}^c_{R}) = 0 ]
= \exp( - \mu( \psi^{(R,\eta)}_t))$
and the same proof shows, when $\gamma - \beta f \geq \eta$, the bound
\begin{equation} \label{globaltsupport}
Q^{\R^d,\beta,\gamma}_{\mu,f} \left[U_{[0,t]}(\overline{D}^c_{R}) = 0 \right] 
\geq Q^{\R^d,\eta}_{\mu} \left[ U_{[0,t]}(\overline{D}^c_{R})=0 \right]. 
\end{equation}
Letting $R \to \infty$ shows that solutions to (\ref{spde})
started from compactly supported $\mu$ remain 
compactly supported at all times. Comparisons can be used to estimate $\psi^{(R,\eta)}$. For example
suppose that $\Delta w \leq (1/2) w^2 - \eta w$ on $D_R$ and $\inf \{ w(x): x \in \partial D_{R-\epsilon} \}
\uparrow \infty$ as $\epsilon \downarrow 0$. Then $\psi^{(R,\eta)}_t \leq w$ on $D_R$ for all $t$ so that,
when $\eta \geq \beta f - \gamma $, for $ \mu$ supported inside $D_R$
\begin{equation} \label{globaltsupport2}
Q^{\R^d,\beta,\gamma}_{\mu,f} \left[U_{[0,\infty)}(\overline{D}^c_{R}) = 0 \right] 
\geq Q^{\R^d,\eta}_{\mu} \left[ U_{[0,\infty)}(\overline{D}^c_{R})=0 \right] \geq 
\exp(-\mu(w)).
\end{equation}

\noindent
\textbf{5. Occupation density.}
Define a putative density $U(s,t,x)$ for $U_{[s,t]}$ on path space $\Omega_D$ 
by setting, for $0 < s \leq t, \, x \in D$,
\begin{equation} \label{putdensity}
U(s,t,x) = \liminf_{\epsilon \downarrow 0} \frac{1}{\epsilon^d} \int^t_s 
U_r \left( [x, x+\epsilon)^d \cap D \right) dr. 
\end{equation}
We also set $U(0,t,x)= \lim_{s \downarrow 0} U(s,t,x)$. 

Sugitani \cite{sugitani} studied the occupation measure  
under the law $Q^{\R^d,0}_{\mu}$ of a DW process on $\R^d$. 
His main result extends to our processes on more general domains, and 
with bounded annihilation/creation rates, as follows:
$Q^{D,\eta}_{\mu}$ almost surely, the process $(s,t,x) \to U(s,t,x)$ is continuous
over $0 <s \leq t, \, x \in D$, and $x \to U(s,t,x)$ acts as a density for $U_{[s,t]}$.
One way of representing the density is via a Green's function formula, as follows. Let
$G^{D,x}_{[0,t]} = \int^t_0 G^{D,x}_s ds$ and $G^{D,\gamma,x}_{[0,t]} = \int^t_0 G^{D,\gamma,x}_s ds$. 
Then, for a fixed $t > \delta >0$ and $x \in D$, 
\begin{equation} \label{densityrep}
U(\delta,t,x) = U_{\delta}(G^{D,x}_{[0,t-\delta]}) + \int^t_{\delta} \int_D G^{D,x}_{[0,t-s]}(z)
(-\eta(z) U_s(dz) ds + M(dz,ds)) \quad \mbox{$Q^{D,\eta}_{\mu}$ a.s.}
\end{equation}
where $M(dz,ds)$ is the martingale measure associated to $U$. 
Moreover, for $k \in \N$,
\begin{equation} \label{densitymoment1}
Q^{D,\eta}_{\mu} [ U(\delta,t,x)^k] \leq C(k,\delta,\mu, \|\eta\|_{\infty}, T) \quad \mbox{for all $x \in D, t \in [\delta,T]$.} 
\end{equation}
When the initial condition is more one may take $\delta = 0$ in (\ref{densityrep}), and indeed this is the case if $\mu$ has a bounded density. 
In particular, when $\mu(dx) = \beta dx$ we have, for $k \in \N$,
\begin{equation} \label{densitymoment2}
Q^{D,\eta}_{\beta dx} [ U(0,t,x)^k] \leq C(k,\beta, \|\eta\|_{\infty}, T) \quad \mbox{for all $x \in D, t \in [0,T]$.} 
\end{equation}
It is not hard to extend Sugitani's proof of the representation formula 
to hold for domains $D$ and with annihilation/creation terms, and also to bound the moments
(inductively in $k$). 
Sugitani's proof of continuity is via Kolmogorov's continuity criterion and moments of 
increments. 
We omit the details since we will follow Sugitani's arguments to establish a similar  
increment estimates for approximations to solutions in section \ref{s4.2}, and 
some of the key steps are illustrated there. 

The second moment formula (\ref{2ndmoments}) lead to the moment
\[
Q^{D,\gamma}_{\mu} \left[ U^2(0,t,x) \right] = 
\left( \mu(G^{D,\gamma,x}_{[0,t]}) \right)^2 
+ \int^t_0 \mu(G^{D,\gamma}_s((G^{D,\gamma,x}_{[0,t-s]})^2))ds.
\]
We use the bound $G^{D,\gamma,x}_t \leq e^{-\gamma t} G^{\R^d,x}_t$ and the estimate
\[
q(t) : = \int (G^{\R^d,x}_{[0,t]}(y))^2 dx \leq C(d) t^{(4-d)/2} \quad \mbox{for all $x,y$ when $d \leq 3$.}
\] 
One can now deduce that
\begin{equation} \label{densityatzero}
Q^{D,\gamma}_{\mu} \left[ \int_D (U(0,t,x))^2 dx \right] \leq C(d,\gamma) \left( \mu^2(1) q(t) +
\mu(1) \int^t_0 q(s) ds \right).
\end{equation}
Again, one has that  $Q^{D,\beta,\gamma}_{\mu,f} \left[ \langle U^2(0,t),1 \rangle \right] \leq
Q^{D,\gamma-\beta}_{\mu} \left[ \langle U^2(0,t),1 \rangle \right]$.

\noindent
\textbf{6. A construction of the exit measure process and the spatial Markov property.} 
A convenient way to construct the exit measure processes was developed in
Perkins \cite{perkins2}, which will be useful in our proof of a 
coupling version of the spatial Markov property lemma below.
This exploits a historical DW process, which is a process of random measures 
on the space of paths, where intuitively the paths trace the positions 
of the ancestors of the particles in a DW process.
To explain the construction, we recall some notation from  \cite{perkins2}
for the historical Brownian DW process (which is 
only required in this one subsection):  
on the path space $C([0,\infty),\R^d)$ let $\mathcal{C}$
be the Borel $\sigma$-field and 
$(\mathcal{C}_t)$ the canonical filtration; 
$\mathcal{M}(C)$ is the space of finite measures on $C([0,\infty),\R^d)$ 
with the topology of weak convergence; on $\Omega_H = C([0,\infty),\mathcal{M}(C)) $ let 
$\mathcal{H}$ be the Borel sets and $(H_t)$ the canonical
coordinate variables. Write $Q_{H,\mu}^{\eta}$ for the law on $(\Omega_H, \mathcal{H})$ of a
historical Brownian DW process, with branching rate one, with bounded measurable annihilation/creation rate
$\eta:\R^d \to \R$, and with initial condition $\mu \in M_F(\R^d)$.
Letting $\mathcal{N}$ be the $Q^{\eta}_{H,\mu}$ null sets we define a
filtration by $\mathcal{H}_t = \cap_{s >t} \sigma \{H_r: r \leq s\} \vee 
\mathcal{N}$. The canonical process $H=(H_t)$ is a historical DW process 
with creation/annihilation $\eta$ on 
$(\Omega_H, \mathcal{H}, \mathcal{H}_t, Q^{\eta}_{H,\mu})$. In the notation from
\cite{perkins2}, it solves the martingale problem $(MP)^{0,\mu}_{1,2I,0,\hat{\eta}}$,
where $\hat{\eta}(s,y) = \eta(y_s)$.

For $(\mathcal{C}_t)$ predictable
$C:[0,\infty) \times C([0,\infty),\R^d) \to \R$, having left continuous paths of bounded variation,  
Theorem 2.23 of Perkins \cite{perkins2}
constructs the integral of $C$ along the paths of a historical DW process as
the limit (in probability, or almost surely along a subsequence)
\begin{equation} \label{pathwiseem}
\int^t_0 H_s(dC_s) = \lim_{n \to \infty} \sum_{i=1}^{\infty} 
\I(t^n_i < t) H_{t^n_i} (C_{t^n_i}-C_{t^n_{i-1}} ), \quad \mbox{where $t^n_i=i2^{n}$.}
\end{equation}
Note that the process $\int^t_0 H_s(dC_s)$ is non-decreasing whenever $C_s$ is also non-decreasing.

Fix a bounded domain $D$ and define $\tau= \tau(y)= \inf \{t \geq 0: y_t \not \in D \}$
for $y \in C([0,\infty),\R^d) $, and where $\inf\{\emptyset\} = 0$.  
Define a measure $u_t(dx)$ on $D$ by 
\[
u_t(dx) = \int_{C([0,\infty),\R^d)} \I(y_t \in dx, \, t \leq \tau, \, \tau >0) H_t(dy).
\]
If $\mu$ is supported inside $D$ we may omit the indicator $\I( \tau>0)$. 
For bounded measurable $\phi: \partial D \to \R$ define
\[
u^{\partial D}_t(\phi) =  \int^t_0 H_s(dC^{\phi}_s), \quad \mbox{where $C^{\phi}_s(y) = \phi (y_{\tau}) 
\I(s > \tau > 0)$.}
\]
By definition $u^{\partial D}(\phi) \geq 0$ if $\phi \geq 0$ and 
$u^{\partial D}(\phi + \psi) = u^{\partial D}(\phi)+u^{\partial D}(\psi)$ almost surely. 
Using these properties, it is not hard to show that
there is an ($\mathcal{H}_t$ measurable) random measure $u_t^{\partial D}$ on $\partial D$ so that
$u^{\partial D}_t(\phi) = \int_{\partial D} \phi(x) du_t^{\partial D}(dx)$ almost surely. 

Let $Z(dy,dt)$ denote the martingale measure associated to $H$ (constructed for example in 
section 2 of  \cite{perkins2}) and use the notation $Z^0(dy,ds) = Z(dy,ds) - \eta(y_s) H_s(dy) ds$, 
that is an integral against $Z^0(dy,ds)$ means the difference of two integrals. 
This notation is borrowed from \cite{perkins2} and will
save space since these two integrals frequently appear with the same integrand.
We may apply the historical Ito formula (\cite{perkins2} Theorem 2.14) to the process
$Y_t(y) = y_{t \wedge \tau}$, since the stopped path $y_{t \wedge \tau}$.
For $\phi \in C^{1,2}_b([0,\infty) \times \R^d)$ 
this yields 
\begin{eqnarray}
\int \phi_{t \wedge \tau} (y_{t \wedge \tau}) H_t(dy) 
&=& \mu(\phi_0 ) 
+ \int^t_0 \int (\dot{\phi} + \Delta \phi)_{s \wedge \tau}(y_{s \wedge \tau}) \I(s \leq \tau) H_s(dy) ds 
\nonumber \\
&& + \int^t_0 \int \phi_{s \wedge \tau} (y_{s \wedge \tau}) Z^0(dy,ds)
\quad \mbox{for $t \geq 0$, a.s.}
\label{hc1}
\end{eqnarray}
Using Proposition 2.7 from \cite{perkins2} we have also, $Q^{\eta}_{H,\mu}$ a.s. for all $t \geq 0$,
\begin{equation}
\int \phi_0 (y_0) \I(\tau=0) H_t(dy)  =  \mu( \I(D^c) \phi_0)  
+ \int^t_0 \int \phi_0 (y_0) \I(\tau=0) Z^0(dy,ds). \label{hc2}
\end{equation}
Finally, Theorem 2.23 of \cite{perkins2} shows, when $C^{\phi}_t(y) = \phi_{\tau}(y_{\tau}) I(t > \tau >0)$,
\begin{equation} 
H_t(C^{\phi}_t)  =  \int^t_0 H_s(dC^{\phi}_s) + \int^t_0 \int \phi_{\tau} (y_{\tau}) \I(s > \tau>0) Z^0(dy,ds)
\quad \mbox{for $t \geq 0$ a.s.} \label{hc3}
\end{equation}
Note that 
\[
u_t(\phi_t) = \int \phi_{t \wedge \tau} (y_{t \wedge \tau}) H_t(dy) 
- \int \phi_0(y_0) \I(\tau = 0) H_t(dy) - H_t(C^{\phi}_t).
\]
Combining this with (\ref{hc1},\ref{hc2},\ref{hc3}) shows that
\begin{eqnarray}
u_t(\phi_t) & = & \mu(\phi_0 \I(D)) - \int^t_0 H_s(dC^{\phi}_s) 
+ \int^t_0 \int \phi_s(y_s) \I(s \leq \tau, \tau >0) Z(dy,ds) \nonumber \\
&& + \int^t_0 \int (\dot{\phi} + \Delta \phi - \eta \phi)_s(y_s) \I(s \leq \tau) H_s(dy) ds
\quad \mbox{for $t \geq 0$ a.s.}
\label{hif} 
\end{eqnarray}
An approximation argument shows that (\ref{hif}) continues to hold for
$\phi \in C^{1,2}([0,t] \times \overline{D})$. 
Taking time independent $\phi$ shows that
the extended martingale problem (\ref{DWMGPD1++}) holds for $(u,u^{\partial D})$.
To check the regularity of the paths, note that 
the path $t \to u_t(\phi)$ is continuous
for $\phi \in C^2_0(\overline{D})$ by (\ref{DWMGPD1+}). 
To establish continuity when $\phi \in C^2_b(\overline{D})$,
we use an estimate of the amount of mass near the boundary of $D$.
For $\epsilon>0$ let $\phi_{D,\epsilon}(x) = (1 - \epsilon^{-1} d(x, \partial D)) \vee 0$.
If $\mu(\partial D) =0$ then $U_t(\partial D) =0$ for all 
$t \geq 0$, $Q^{\R^d,\eta}_{\mu}$ almost surely.
This follows since $t \to U_t(\partial D)$ is continuous (see \cite{perkins1} Theorem III.5.1
and use absolute continuity between $Q^{\R^d,\eta}_{\mu}$ and $Q^{\R^d,0}_{\mu}$) 
and has zero first moment.
Since $\phi_{D,\epsilon} \downarrow \I(\partial D)$, for such $\mu$ we have
\begin{equation} \label{massnrbdy}
\mbox{$\sup_{t \leq T} U_t(\phi_{D,\epsilon}) \downarrow 0$ as $\epsilon \downarrow 0$, 
$\quad Q^{\R^d,\eta}_{\mu}$ almost surely.}
\end{equation}
Note the simple coupling 
\begin{equation} \label{simplecoupling}
u_t \leq \tilde{u}_t \quad \mbox{for all $t \geq 0$}
\end{equation}
where $\tilde{u}$ is a DW($\R^d$) process with creation/annihilation $\eta$, both with the same initial
condition. This is natural since $\tilde{u}$ does not have the killing on the boundary of $D$, and 
follows by setting $\tilde{u}_t(dx) = \int \I(y_t \in dx) H_t(dy)$. 
This coupling ensures that (\ref{massnrbdy}) is also true for $\sup_{t \leq T} u_t(\phi_{D,\epsilon})$.
This implies the continuity of $t \to u_t(\phi)$ for $\phi \in C^2_b(\overline{D})$, and hence of
$t \to u_t \in \mathcal{M}(D)$.
The continuity of $t \to u^{\partial D}_t(h)$ for $h \in C^2_b(\partial D)$ follows
from the continuity of all other terms in (\ref{DWMGPD1++}). 
This implies the paths $t \to u_t^{\partial D}$ are continuous, and they are 
non-decreasing by definition. So $u$ is a DW process on $D$, with annihilation/creation rate $\eta$, 
with exit measure $u^{\partial D}$, and with initial condition $\mu \I(D)$. 

Dynkin \cite{dynkin} established a spatial 
Markov property for DW processes. We will need the following
slight variant of the standard statement. The lemma is reasonably clear at the 
approximation level of discrete branching trees, but we choose to give a 
proof in the continuum setting using historical calculus. 

\begin{Lemma} \label{smlemma}
Let $D^- \subseteq D^+$ be bounded domains.
Fix $\mu \in \mathcal{M}(D^+)$. There exists a coupling of  
three processes: $u$ a DW($D^+,\gamma$)
process with initial condition $\mu$; $u^-$ a DW($D^-,\gamma$)
process with initial condition $\mu \I(D^-)$; and $u^+$ which, conditional
on $\sigma\{u^-\}$, is a DW($D^+,\gamma$) process with initial condition 
$\mu \I(D^+ \setminus D^-) + u^{\partial D^-,-}_{[0,\infty)}|_{D^+}$; and moreover these 
processes satisfy the splitting
\begin{equation} \label{smsplitting}
\left(u_{[0,\infty)}, u^{\partial D^+}_{\infty} \right) = 
\left(u^{-}_{[0,\infty)} + u^+_{[0,\infty)}, 
u^{\partial D^-,-}_{\infty} |_{\partial D^+} + u^{\partial D^+,+}_{\infty} \right) \qquad
\mbox{almost surely.}
\end{equation}
\end{Lemma}

\noindent
\textbf{Proof.}
Define $\tau^{\pm}= \inf \{t: y_t \not \in D^{\pm}\}$.
As above, we suppose the canonical process $(H_t)$ is a historical DW process
with constant annihilation rate $\gamma$ on 
$(\Omega_H, \mathcal{H}, \mathcal{H}_t, Q^{\gamma}_{H,\mu})$.
Define
\begin{eqnarray*}
&& u_t(dx) = \int  \I(y_t \in dx, \,t \leq \tau^+) \, H_t(dy),  \quad 
u^-_t(dx) = \int \I(y_t \in dx, \,t \leq \tau^-, \, \tau^->0) \, H_t(dy) \\
&& u^{\partial D^+}_{t}(\phi)   =  \int^t_0 H_s(dC^{\phi,+}_s), \qquad 
u^{\partial D^-,-}_{t}(\phi)   =  \int^t_0 H_s(dC^{\phi,-}_s),
\end{eqnarray*}
where $C^{\phi,+}_s(y) = \phi(y_{\tau^+}) \I(s > \tau^+)$ for bounded measurable
$\phi:\partial D^+ \to \R$ and $C^{\phi,-}_s(y) = \phi(y_{\tau^-}) 
\I(s > \tau^- > 0)$ for bounded measurable $\phi: \partial D^- \to \R$. 
The earlier arguments leading to (\ref{hif}) 
show that this defines processes $(u,u^{\partial D^+})$ and $(u^-,u^{\partial D^-,-})$ 
having laws $Q^{D^+,\gamma}_{\mu}$ and $Q^{D^-,\gamma}_{\mu \I(D^-)}$ as desired.

There exist random measures $ \mathbf{I} $ on $[0,\infty) \times D^+$
and $\mathbf{I}^{\partial D^+}$ on $[0,\infty) \times \partial D^+$ satisfying, 
for bounded measurable $\phi:[0,\infty) \times D^+ \to \R$ and 
$\psi:[0,\infty) \times  \partial D^+ \to \R$, 
\begin{eqnarray*}
\mathbf{I}(\phi) & = & \int^{\infty}_0 \int \phi_{t-\tau^-}(y_t) \I(t \in (\tau^-,\tau^+]) H_t(dy) dt \\
\mathbf{I}^{\partial D^+}(\psi) & = & \int^{\infty}_{0} H_s(d\hat{C}^{\psi}_s) \quad \mbox{where
$\hat{C}^{\psi}_s(y) = \psi_{\tau^+-\tau^-}(y_{\tau^+}) \I(s > \tau^+ > \tau^-)$.}
\end{eqnarray*}
We will show that conditional on $\sigma\{u^-\}$, the law
of $(\mathbf{I},\mathbf{I}^{\partial D^+})$ is that of $(U_t(dx)dt, dU^{\partial D}_{t}(dx))$
under $Q^{D^+,\gamma}_{\mu \I(D^c) + u^{\partial D^-,-}_{\infty}|_{D^+}}$.
This allows us to define $(u^+,u^{\partial D^+,+})$ with the required law 
so that 
\[
\mathbf{I}(\phi) = \int^{\infty}_0 u^+_t(\phi_t) dt \quad \mbox{and} \quad 
\mathbf{I}^{\partial D^+}(\psi) = \int^{\infty}_0 du^{\partial D^+,+}_s(\psi_s). 
\]
The splitting (\ref{smsplitting})
for the total occupation measures holds since  
\[
u^-_{[0,\infty)}(A) = \int^{\infty}_0 \int \I(y_t \in A, t \leq \tau^-, \tau^->0) H_t(dy) dt
= \int^{\infty}_0 \int \I(y_t \in A, t \leq \tau^-) H_t(dy) dt
\]
while
\[
u_{[0,\infty)}(A) = \int^{\infty}_0 \int \I(y_t \in A, t \leq \tau^+) H_t(dy) dt, \quad
\mathbf{I}([0,\infty) \times A) = \int^{\infty}_0 \int \I(y_t \in A, t \in (\tau^-,\tau^+]) H_t(dy) dt.
\]
Note that for $A \subseteq \partial D^+$
\[ 
u^{\partial D^+}_{\infty}(A) = \int^{\infty}_0 H_s(dC^{A,+}_s), \quad
u^{\partial D^-,-}_{\infty}(A) = \int^{\infty}_0 H_s(dC^{A,-}_s), \quad
\mathbf{I}^{\partial D^+}([0,\infty) \times A) = \int^{\infty}_0 H_s(d\hat{C}^{A}_s)
\]
where
\[
C^{A,+}_s(y) = \I(y_{\tau^+} \in A, s > \tau^+), \quad 
C^{A,-}_s(y) = \I(y_{\tau^-} \in A, s > \tau^- >0), \quad
\hat{C}^{A}_s(y) = \I(y_{\tau^+} \in A, s > \tau^+ > \tau^-). 
\]
Note also that $C^{A,-}_s(y) = \I(y_{\tau^+} \in A, s > \tau^+=\tau^- >0)$ since $A \subseteq \partial D^+$.
Then the splitting for the total exit measures holds since 
\[
C^{A,+}_s = \hat{C}^{A}_s + C^{A,-}_s + \I(y_{\tau^+} \in A, s>\tau^+ = 0)
\]
and $H_t(\tau^+ =0) = 0$ for all $t \geq 0$, $Q^{\gamma}_{H,\mu}$ a.s.

To identify the conditional law of $(\mathbf{I},\mathbf{I}^{\partial D^+})$ we will use the Laplace functionals.
Choose non-negative, continuous bounded  $h^{(2),\pm}:[0,\infty) \times D^{\pm}$ and
$h^{(3),\pm}: [0,\infty) \times \partial D^{\pm}$ which are zero for $t \geq T$. Let $\hat{\phi}^{\pm} 
\in C^{1,2}([0,T] \times D^{\pm}) \cap C([0,T] \times \overline{D^{\pm}})$ be the unique non-negative mild solution 
to
\[
\mbox{$- \dot{\hat{\phi}}^{\pm} = \Delta \hat{\phi}^{\pm} - \gamma \hat{\phi}^{\pm} - \frac12 (\hat{\phi}^{\pm})^2 + h^{(2),\pm}$ on 
$D^{\pm}$, $\quad \hat{\phi}^{\pm}_T = 0$ and  $\hat{\phi}^{\pm} = h^{(3),\pm}$ on $\partial D^{\pm}$.} 
\]
We extend these by letting $\hat{\phi}^{\pm}_t =0$ for $t \geq T$.  It is enough to show 
\begin{equation} 
E\left[ e^{- \mathbf{I}(h^{(2),+}) - \mathbf{I}^{\partial D^+}(h^{(3),+})} \left| \sigma\{u^-\} 
\right. \right] = e^{-u^+_0(\hat{\phi}^+_0)}
\label{smpeq10}
\end{equation}
where $u^+_0$ is defined as $\mu \I(D^+ \setminus D^-) + u^{\partial D^-,-}_{\infty} |_{D^+}$.
This in turn will be implied by 
\begin{equation} 
E\left[ e^{+ u^+_0(\hat{\phi}^+_0)} \; e^{- \mathbf{I}(h^{(2),+}) 
- \mathbf{I}^{\partial D^+}(h^{(3),+})} \; 
e^{-\int^{\infty}_0 u^-_s(h^{(2),-}_s) ds - 
\int^{\infty}_0 du_s^{\partial D^-,-}(h^{(3),-}_s)}  \right] = 
e^{-\mu(\hat{\phi}^-_0)}
\label{smpeq20}
\end{equation}
for all $h^{(2),-},h^{(3),-}$. 

A comparison of the martingale problem for time dependent test functions (\ref{DWMGPD1+++}) with (\ref{hif}) shows 
that if $\phi:[0,\infty) \times \partial D^+ \to \R$ and $C^{\phi,+}_s(y) = \phi_{\tau^+}(y_{\tau^+}) \I(s > \tau^+)$
then $\int^t_0 H_s(dC^{\phi,+}_s) = \int^t_0 du^{\partial D^+}_s(\phi_s)$ (and a similar
representation for a $du^{\partial D^-,-}_t$ integral). 
Suppose first that $\hat{\phi}^-$ is in $C^{1,2}_b([0,T] \times \overline{D^-})$. 
Then using (\ref{hif}) on the domain $D^-$, one finds that $\Xi^-_t$ defined by 
\[
\Xi^-_t = e^{-u^-_t(\hat{\phi}^-_t) - \int^t_0 u^-_s(h^{(2),-}_s) ds - \int^t_0 du^{\partial D^-,-}_s(h^{(3),-}_s)}
\]
is a martingale which satisfies 
$\Xi^-_{\infty} = e^{-\int^{T}_0 u^-_t(h^{(2),-}_t) dt - \int^{T}_0 du^{\partial D^-,-}_t(h^{(3),-}_t)}$ and 
\begin{equation}
\Xi^-_t =  e^{-\mu(\hat{\phi}^-_0)} - \int^t_0 \int \Xi_s^- \, \hat{\phi}^-_s(y_s)
\I(s \leq \tau^-, \tau^->0) Z(dy,ds) \quad \mbox{for $t \geq 0$ a.s.}
\label{smpeq30}
\end{equation}
An approximation argument shows (\ref{smpeq30}) continues to hold when  
$\hat{\phi}^{-} \in C^{1,2}([0,T] \times D^{-}) \cap C([0,T] \times \overline{D^{-}})$.
We now develop a similar representation for the Laplace functional of $(\mathbf{I},\mathbf{I}^{\partial D^+})$.
Set 
\[
Y_t = ((t-\tau^-)_+ - (t-\tau^+)_+), y_{t \wedge \tau^+}) 
= \left(\int^t_0 \I(s \in (\tau^-,\tau^+]) ds, \int^t_0 \I(s \leq \tau^+) dy(s)\right).
\]
For $\phi \in C^{1,2}_b([0,\infty) \times \R^d)$ we have, using the historical Ito formula again,
\begin{eqnarray*}
\int \phi(Y_t) H_t(dy) &=& \mu(\phi_0) + \int^t_0 \int \phi(Y_s) Z^0(dy,ds)  \\
&& + \int^t_0 \int \dot{\phi}(Y_s) \I(s \in (\tau^-,\tau^+]) + \Delta \phi(Y_s) \I(s \leq \tau^+) H_s(dy) ds.
\end{eqnarray*}
On $\{t > \tau^+\}$ we have $Y_t = (\tau^+-\tau^-,y_{\tau^+})$ so that
\begin{eqnarray}
\int \phi(Y_t) \I(t > \tau^+) H_t(dy) &=& \int \phi_{\tau^+-\tau^-}(y_{\tau+}) \I(t > \tau^+) H_t(dy)  \nonumber \\
&=& H_t(\hat{C}^\phi_t) + \int \phi_0(y_{\tau^-}) \I(t>\tau^+=\tau^-) H_t(dy).
\label{smpeq40b}
\end{eqnarray}
 Using Theorem 2.23 of \cite{perkins2},
\[
H_t(\hat{C}^{\phi}_t)
= \int^t_0 H_s(d\hat{C}^{\phi}_s) + \int^t_0 \int \phi_{\tau^+-\tau^-}(y_{\tau^+}) \I(s > \tau^+ >\tau^-) Z^0(dy,ds).
\]
On $\{t \leq \tau^-\}$ we have $Y_t = (0,y_t)$ so that 
\begin{eqnarray}
\int \phi(Y_t) \I(t \leq \tau^-) H_t(dy)  &=& \int \phi_0(y_{t \wedge \tau^-}) \I(t \leq \tau^-) H_t(dy) \nonumber  \\
&=& \int \phi_0(y_{t \wedge \tau^-}) H_t(dy) + \int \phi_0(y_{\tau^-}) \I(t>\tau^-) H_t(dy).
\label{smpeq40d}
\end{eqnarray}
By the historical Ito formula again
\[
\int \phi_0(y_{t \wedge \tau^-}) H_t(dy) =
\mu(\phi_0) + \int^t_0 \int \phi_0(y_{s \wedge \tau^-}) Z^0(dy,ds) + \int^t_0 \int \Delta \phi_0(y_s) \I(s \leq \tau^-) H_s(dy) ds.
\]
We may combine some terms from (\ref{smpeq40b}) and (\ref{smpeq40d}) as follows, again using
Theorem 2.23 of \cite{perkins2},
\begin{eqnarray*}
&& \hspace{-.4in} \int \phi_0(y_{\tau^-}) \I(t>\tau^-) H_t(dy) - \int \phi_0(y_{\tau^-}) \I(t>\tau^+=\tau^-) H_t(dy) \\
& = & \int \phi_0 \I(D^+)(y_{\tau^-}) \I(t > \tau^-) H_t(dy) \\
&=& H_t(C^{\phi_0 \I(D^+),-}_t) + \int \phi_0 \I(D^+)(y_0) \I(t > \tau^- =0) H_t(dy) \\
& = & \int^t_0 H_s(dC^{\phi_0\I(D^+),-}_s) + \mu(\phi_0 \I(D^+ \setminus D^-)) 
+ \int^t_0 \int \phi_0 \I(D^+)(y_{\tau^-}) \I(t > \tau^-) Z^0(dy,ds).
\end{eqnarray*}
The last six displayed equations hold for all $t \geq 0$ almost surely, and a little book-keeping combines them to 
yield
\begin{eqnarray*}
&& \hspace{-.4in} \int \phi_{t-\tau^-}(y_t) \I(t \in (\tau^-,\tau^+]) H_t(dy) + \int^t_0 H_s(d\hat{C}^{\phi}_s)
- \left(\mu(\phi_0\I(D^+ \setminus D^-)) + \int^t_0 H_s(dC^{\phi_0 \I(D^+),-}_s)\right) \nonumber \\
& = & \int^t_0 \! \int \phi_{s-\tau^-}(y_s) \I(s \in (\tau^-,\tau^+]) Z(dy,ds) +
\int^t_0 \! \int (\dot{\phi}+\Delta \phi - \gamma \phi)_{s-\tau^-}(y_s) \I(s \in (\tau^-,\tau^+]) H_s(dy) ds.
\label{smpeq50}
\end{eqnarray*}
If we now suppose that $\dot{\phi} + \Delta \phi - \gamma \phi - (1/2) \phi^2 + h =0$ for some $h \in C_b([0,\infty) \times \R^d)$
then $\Xi^+_t$ defined by 
\begin{eqnarray*}
\Xi^+_t & = & \exp( + \mu(\phi_0 \I(D^+ \setminus D^-)) + \int^t_0 H_s(dC^{\phi_0 \I(D^+),-}_s)
- \int \phi_{t-\tau^-}(y_t) \I(t \in (\tau^-,\tau^+]) H_t(dy)) \\
&& \hspace{.4in} \exp(- \int^t_0 H_s(d\hat{C}^{\phi}_s) - \int^t_0 \int h_{s-\tau^-}(y_s) 
\I(s \in (\tau^-,\tau^+]) H_s(dy) ds)
\end{eqnarray*}
is a non-negative local martingale and satisfies
\begin{equation}
\Xi^+_t = 1 - \int^t_0 \int \Xi^+_s \phi_{s-\tau^-}(y_s) \I(s \in (\tau^-,\tau^+]) Z(dy,ds)
\quad \mbox{for $t \geq 0$ a.s.}
\label{smpeq60}
\end{equation}
An approximation argument shows that (\ref{smpeq60}) continues to hold with $\phi = \hat{\phi}^+$
and $h=h^{(2),+}$.
With this choice, and extending $\hat{\phi}^+_t =0$ for $t \geq T$, we let $t \to \infty$ 
to find that
\[
\Xi^+_t \to \exp\left( + u^+_0(\hat{\phi}^+_0) - \mathbf{I}(h^{(2),+})
 - \mathbf{I}^{\partial D^+}(h^{(3),+}) \right).
\] 
Combining (\ref{smpeq30}) and (\ref{smpeq60}) we have 
\begin{eqnarray}
&& \hspace{-.6in}
 e^{+ u^+_0(\hat{\phi}^+_0) - \mathbf{I}(h^{(2),+}) - \mathbf{I}^{\partial D^+}(h^{(3),+})} \; 
e^{-\int^{\infty}_0 u^-_t(h^{(2),-}_t) dt - \int^{\infty}_0 du_t^{\partial D^-,-}(h^{(3),-}_t)} \nonumber \\
& = & e^{-\mu(\hat{\phi}^-_0)} 
- \int^{\infty}_0 \int \Xi^-_t \Xi^+_t \hat{\phi}^+_{t-\tau^-}(y_t) \I(t \in (\tau^-,\tau^+]) Z(dy,dt) \nonumber \\
&&  - \int^{\infty}_0 \int \Xi_t^+ \Xi_t^-  \hat{\phi}^-_t (y_t)
\I(t \leq \tau^-, \tau^->0) Z(dy,dt)
\label{smpeq70}
\end{eqnarray}
(noting that the cross variation of the two stochastic integrals is zero). 
Note that $\Xi^+_t \leq \exp(+\|\hat{\phi}^+\|_{\infty} u^+_0(1))$. The quadratic 
variation of the first stochastic integral on the right hand side of (\ref{smpeq70}) is
therefore bounded 
\[
 \|\hat{\phi}^+\|_{\infty}^2 \exp(+ 2 \|\hat{\phi}^+\|_{\infty}(\mu(1) + u^{\partial D^-,-}_{\infty}(1)) u_{[0,\infty)}(1)
\]
The finiteness of small positive exponential moments for the total exit measure  
$u^{\partial D^-,-}_{\infty}(1)$ implies that the stochastic integral has mean zero 
if the norm $\|\hat{\phi}^+\|_{\infty}$ is small, which in turn is satisfied
if the norms $\|h^{(2),+}\|_{\infty}$ and $\|h^{(3),+}\|_{\infty}$ are small enough. 
Similarly the other stochastic integral is mean zero. 
Taking expectations in (\ref{smpeq70}) establishes the conditional Laplace functional (\ref{smpeq10}) when
 $\|h^{(2),+}\|_{\infty}$ and $\|h^{(3),+}\|_{\infty}$ are sufficiently small. 
But this is sufficient to characterize the conditional law. \qed
\subsection{Existence and uniqueness via change of measure} \label{s3.2}
\textbf{Proof of Theorem \ref{E!}.} 
Many of the change of measure arguments we will use can be found in 
section IV.1 of Perkins \cite{perkins1}.
Fix $f \in \mathcal{B}(D,[0,1])$, and on the path space $(\Omega_D, \mathcal{U}, \mathcal{U}_t, Q^{D,0}_{\mu})$
use the path space density from (\ref{putdensity}) to set $ V_t(x) = f(x) \exp(- U(0,t,x))$, where
we use the convention $\exp(-\infty)=0$.
Note that $V$ is $\mathcal{P} \times \mathcal{B}(D)$ measurable and has non-increasing paths as
required. Define the stochastic integral 
\begin{equation} \label{cofmM}
M^{\beta,\gamma,f}_t =  \int^t_0  \int_D (\beta V_s(x)- \gamma) M(dx,ds)
\end{equation}
where $M(ds,dx)$ is the martingale measure on $[0,\infty) \times D$ constructed from $U$ 
under $Q^{D,0}_{\mu}$ (as explained in section 3.1.1). 
A simple version of the arguments in \cite{perkins1} Theorem IV.1.6.(b)
implies that the stochastic exponential $\mathcal{E}_t(M^{\beta,\gamma,f})$ is a true 
martingale, so that there exists a measure
$Q_{\mu,f}^{D,\beta,\gamma}$ on the path space $(\Omega_D, \mathcal{U})$ satisfying 
$dQ_{\mu,f}^{D,\beta,\gamma}/dQ^{D,0}_{\mu}= \mathcal{E}_t(M^{\beta,\gamma,f})$ on $\mathcal{U}_t$.
It is straightforward to check via Girsanov's theorem that
the triple $(U_t,U^{\partial D}_t, V_t)$ solve, under $Q_{\mu,f}^{D,\beta,\gamma}$, 
the extended martingale problem (\ref{MGP3}) and (\ref{MGPD1+}) for solutions to (\ref{spde}).  
The regularity of $U(s,t,x)$ will imply that (\ref{MGP2}) holds, namely for all bounded
measurable $\phi$
\begin{equation} \label{part2}
(V_t,\phi) = (f,\phi) - \int^t_0 U_s(V_s \phi) ds,  
\quad \mbox{$Q_{\mu,f}^{D,\beta,\gamma}$ a.s.} 
\end{equation}
This would be trivial, by the fundamental theorem of calculus, were $U_s(dx)$ to have  
continuous density in $(s,x)$. We show at the end of this section that it follows in our 
more singular setting. This completes the proof of part (i) of Theorem \ref{E!}. 

To show uniqueness of solutions, we take a solution $(u,v)$ to (\ref{spde}) on $D$, 
on a filtered space $(\Omega, (\mathcal{F}_t), \mathcal{F}, P)$,
with an arbitrary $\mathcal{F}_0$ measurable starting condition.
Let $m(dx,dt)$ be the local martingale measure extending the 
local martingales $m_t(\phi)$. We now 
reverse the change of measure argument. The local martingale defined by 
\[
m^{\beta,\gamma,f}_t = \mathcal{E}_t \left( \int \! \int (\gamma -\beta v_s(x)) m(dx,ds) \right)
\] 
is in fact a true martingale (argue as in \cite{perkins1} 
Theorem IV.1.6.(a), noting that the deterministic initial conditions there are not needed for
this argument). Defining $Q$ by $dQ/dP= m^{\beta,\gamma,f}_t$ on $\mathcal{F}_t$, one finds that
$(u_t)$ solves the martingale problem for a DW($D,0$) process under $Q$, and so $Q$ is determined 
on $\sigma\{u_s:s \leq t\}$ by the law of $u_0$. 
By Sugitani's results listed in section \ref{s3.1}, the 
occupation measures $u_{[s,t]}$ have densities (under $P$ or $Q$), and moreover 
that there is a density $(t,x) \to u(s,t,x)$ that is 
continuous on $\{(s,t,x): 0< s \leq t, x \in D\}$. 
Knowing this one shows that $v$ is given by (see the argument at the end of this section)
\begin{equation} \label{vformula}
v_t(x) = v_0(x) e^{- u(0,s,x)} \quad \mbox{for all $t \geq 0$, for almost all $x$, $P$ a.s.}
\end{equation}
where $u(0,s,x) = \lim_{r \downarrow 0} u(r,s,x)$. Now we may replace the change of measure by
\[
dQ/dP = \mathcal{E}_t \left( \int \! \int \left(\gamma-\beta v_0(x) e^{-u(0,s,x)} \right) m(dx,ds) \right)
 \quad \mbox{on $\mathcal{F}_t$}.
\]
since the two exponential martingales 
are equal upto a null set. Inverting one finds 
\[
dP = \mathcal{E}_t \left( \int \! \int \left(\beta v_0(x) e^{-u(0,s,x)} - \gamma \right) \tilde{m}(dx,ds) \right)
dQ \quad \mbox{on $\mathcal{F}_t$}.
\]
where $\tilde{m}(dx,ds)$ is the local martingale measure determined by the local
martingales found in the martingale problem of $u$ under $Q$.
The stochastic exponential is a measurable function of
$(v_0, \{u_s: s \leq t\})$. Applying the Markov property of $u$ under $Q$ at time $t=0$, we find that 
$P$ is determined on $\sigma\{u_s:s \leq t\}$ by the law of $(u_0,v_0)$. Indeed, for each $t \geq 0$ there 
is a measurable kernel $p_t(\mu,f,B)$ for 
$f \in \mathcal{B}(D,[0,1]),\, \mu \in \mathcal{M}(D)$ and Borel $B \subseteq \mathcal{M}(D)$
so that
\begin{equation} \label{kernel}
P[ u_t \in dB] = \int_{\mathcal{B}} \int_{\mathcal{M}(D)} p_t(\mu,f,B) P[u_0 \in d\mu, \, v_0 \in df]. 
\end{equation}
Under $Q$, and therefore under $P$, there exists an exit measure process
$u^{\partial D}$ satisfying $u^{\partial D}_t=R_t(\{u_s: s \leq t\})$ almost surely.
Another Girsanov calculation starting with
the extended martingale problem under $Q$, shows that the exit measures satisfy the
extended martingale problem for (\ref{spde}) under $P$. 
This completes the proof of parts (ii),(iii) and (iv) of Theorem \ref{E!}.
For part (v), note that 
the measurability of the maps $(\mu,f) \to Q_{\mu,f}^{D,\beta,\gamma}(B)$ follows from 
the analogous measurability of $Q^{D,0}_{\mu}$ and the fact that, for $t \geq 0$, the derivative
$dP/dQ^{D,0}_{\mu}$ on $\mathcal{F}_t$ is a measurable map of $(f, \{u_s: s \leq t\})$.
The strong Markov property for the laws $Q^{D,\beta,\gamma}_{\mu,f}$ can be deduced either
from the strong Markov property for the DW laws $Q^{D,\gamma}_{\mu}$, or more easily 
from the uniqueness of the one dimensional distributions for the martingale problem 
(use (\ref{kernel}) and follow the argument in \cite{perkins1} Theorem II.5.6). \qed

\noindent \textbf{Remark.} It is straightforward to show that uniform integrability of the 
martingale $t \to \mathcal{E}_t(M^{\beta,\gamma,1})$ under $Q^{\R^d,0}_{\mu}$ for some $\mu \neq 0$
(or that of $t \to \mathcal{E}_t(M^{\beta,0,1})$ under $Q^{\R^d,\gamma}_{\mu}$) 
is equivalent to certain death for the parameters $(\beta,\gamma)$. See the remark in
section \ref{s7.1}. 

\textbf{Proof of (\ref{part2}).}
Suppose $(\nu_s:0 \leq s \leq t)$ is a continuous path in $\mathcal{M}(\R^d)$. Suppose
also that $(n(s,x):s \leq t, x \in \R^d)$ is bounded and continuous and that $\int^s_0 \nu_r dr
= n(s,x)dx$ as measures for $s \leq t$. Let $T_{\delta}$ denote the heat semigroup with 
generator $\Delta$, and write $T^*_{\delta}$ for the dual semigroup
acting on measures. For continuous compactly supported $\psi$, we shall argue that
\begin{eqnarray*}
\int (1 - e^{- n(t,x)}) \psi(x) dx & = & 
\lim_{\delta \downarrow 0} \int (1 - e^{- T_{\delta} n(t,x)}) \psi(x) dx \\
& = & \lim_{\delta \downarrow 0} \int^t_0  \int 
e^{- T_{\delta} n(s,x)}  \psi(x) T^*_{\delta} \nu_s(dx) ds \\
&=& \int^t_0 \int e^{-n(s,x)} \psi(x) \nu_s(dx) ds.
\end{eqnarray*} 
Since the measures $T^*_{\delta} \nu_s$ have densities, the paths 
$s \to T_{\delta}n(s,x)$ are absolutely continuous for 
almost all $x$. So the second equality above follows by 
expanding $s \to \exp( - T_{\delta} n(s,x))$
as an integral of its derivative over $[0,t]$ and then interchanging
$s$ and $x$ integrals. The third equality holds since $T_{\delta}n(s,x) \to n(s,x)$ uniformly
over $s \leq t$ and $x$ in the support of $\psi$, and since 
$T^*_{\delta} \nu_s \to \nu_s$ weakly.

The assumptions above hold, $Q^{D,\beta,\gamma}_{\mu,f}$ almost surely,  
if we take $\nu_s(dx) = \overline{\psi}(x) U_{\epsilon+s}(dx)$,
for some $\epsilon>0$ and continuous $\overline{\psi}$, with 
compact support in $D$, and take the 
corresponding density $n(s,x)=\overline{\psi}(x) U(\epsilon,\epsilon+s,x)$.
Choosing continuous $\psi$ compactly supported in $D$ and $\overline{\psi}\equiv 1$ on the support of $\psi$ 
we obtain, almost surely,  
\[
\int e^{-U(\epsilon,t,x)} \psi(x) dx = \int \psi(x) dx - \int^t_{\epsilon} 
\int e^{-U(\epsilon,s,x)} \psi(x) U_s(dx) ds.
\]
To pass to the limit $\epsilon \downarrow 0$ in this equation note that for $t \geq \epsilon$
\[
\int \left|e^{-U(\epsilon,t,x)} - e^{-U(0,t,x)} \right| dx  
 \leq  \int U(0,\epsilon,x) dx =  \int^{\epsilon}_0 U_s(1) ds \to 0 \quad \mbox{as $\epsilon \downarrow 0$,}
\]
and similarly 
\begin{eqnarray*}
&& \hspace{-.3in} \left| \int^t_{\epsilon} \int e^{-U(\epsilon,s,x)} \psi(x) U_s(dx) ds - 
\int^t_{\epsilon} \int e^{-U(0,s,x)} \psi(x) U_s(dx) ds \right| \\
& \leq & \|\psi\|_{\infty} \int^t_{\epsilon} U_s(U(0,\epsilon)) ds \\
& \leq & \|\psi\|_{\infty} \int U(0,t,x) U(0,\epsilon,x) dx 
\end{eqnarray*}
which converges in $\mathcal{L}^2$ to zero, as $\epsilon \downarrow 0$, 
by (\ref{densityatzero}). This leads to 
\[
\int e^{-U(0,t,x)} \psi(x) dx = \int \psi(x) dx - \int^t_0 
\int e^{-U(0,s,x)} \psi(x) U_s(dx) ds
\]
for continuous compactly supported $\psi$.
We may extend to general bounded measurable $\psi:D \to \R$ by monotone class arguments,
and choosing $\psi=f \phi$ for $\phi \in C^0_b(\overline{D})$ shows that (\ref{part2}) holds. \qed

\textbf{Proof of (\ref{vformula}).} First, we may use a density argument to
ensure $ \langle v_t,\phi \rangle = \langle v_0,\phi \rangle - \int^t_0 u_s(v_s \phi) ds $ for all $t \geq 0$ and all 
continuous $\phi$ with compact support in $D$, $P$ a.s.
Now we can argue pathwise. Fix a sample point $\omega$ so that the above
happens and so that also $(r,s,x) \to u(r,s,x)$ is continuous over $x \in D, 0<r \leq s \leq t$. 
Choose a partition $0<s=t_0 <t_1 < \ldots <t_N =t$ with $\max_i (t_i-t_{i-1})=\delta$.
Choosing $\phi(x) = \exp(-u(t_{n-1},t,x)) \psi(x)$, for some continuous $\psi \geq 0$ with compact support
$A \subseteq D$, we find
\[
\left \langle v_{t_n}, e^{-u(t_{n-1},t)} \psi \right \rangle
= \left \langle v_{t_{n-1}}, e^{-u(t_{n-1},t)}\psi \right \rangle
- \int^{t_n}_{t_{n-1}} u_s \left( v_s e^{-u(t_{n-1},t)} \psi \right) ds. 
\]
A little rearrangement leads to 
\[
\left \langle v_{t_n}, e^{-u(t_{n},t)} \psi \right \rangle
= \left \langle  v_{t_{n-1}}, e^{-u(t_{n-1},t)} \psi \right \rangle
- E^1_n + E^2_n 
\]
where $E^i_n \geq 0$ are given by  
\[
 E^1_n = \int^{t_n}_{t_{n-1}} u_s \left( e^{-u(t_{n-1},t)} (v_s-v_{t_n}) \psi \right) ds 
\leq \|\psi\|_{\infty} \left \langle u(t_{n-1},t_n) \I(A), v_{t_{n-1}}-v_{t_n} \right \rangle
\]
and, using $0 \leq e^{z}-1-z \leq 2z^2$ for $z \in [0,1]$,
\[
E^2_n = 
\left \langle v_{t_{n}} \psi e^{-u(t_{n-1},t)}, e^{u(t_{n-1},t_n)} -1-u(t_{n-1},t_n) \right  \rangle
\leq 2 \|\psi\|_{\infty}  
\left \langle \I(A), u^2(t_{n-1},t_n) \wedge 1 \right \rangle.
\]
Now sum over the partition to find
\[
\left \langle v_{t}, \psi \right \rangle
= \left \langle v_s, e^{-u(s,t)} \psi \right \rangle
- \sum_n E^1_n + \sum_n E^2_n. 
\]
The error bounds above, and the fact that 
$\max \{ u(s',t',x): s \leq s' \leq t' \leq t, |s'-t'| \leq \delta, x \in A\} \to 0$ as
$ \delta \to 0$, imply that $\sum_n E^i_n \to 0$  as $\delta \to 0$.
Thus $v_t(x) = v_s(x) \exp(-u(s,t,x))$ for almost all $x$. 
The non-increasing paths $t \to v_t(x)$ ensures that 
$v_t(x) = v_s(x) \exp(-u(s,t,x))$ for all $0<s<t$, for almost all $x$. Since
$\lim_{s \downarrow 0} v_s(x) = v_0(x)$ for almost all $x$ (use the non-increasing paths of
$v$ and (\ref{MGP2})) we reach the conclusion (\ref{vformula}).  \qed

We conclude this section with the proof of local extinction of solutions. 
\begin{Lemma} \label{localextinction} 
Suppose that $\gamma >0$. Then, for compact $A \subseteq \R^d$, 
\[
Q^{\R^d,\beta,\gamma}_{\mu,f} \left[ U_{[t,\infty)}(A)=0 \right] \to 1 \quad \quad 
\mbox{as $ \; t \to \infty$.}
\]
\end{Lemma}

\textbf{ Proof.} Although not directly used, we first give a weaker fixed $t$ extinction 
result whose proof is similar to the support results in section \ref{s3.3}.
Fix $R>0$ and let $0 \leq \psi_R \leq 1$ be smooth and satisfy $\{\psi_R >0\} = B_R =
\{x: |x| <R\}$.
Let $\phi^{(M)}_t(x)$ be the solution, for $x \in \R^d, \, t \geq 0$, of
\[
\partial_t \phi = \Delta \phi - \gamma \phi - \frac{\phi^2}{2} + \beta \phi 
\I(t \leq 1), \quad  \phi_0 = M \psi_R.
\]
Then stochastic calculus shows that, under $Q^{\R^d,\beta,\gamma}_{\mu,f}$, 
\[
d \left( e^{-U_s(\phi^{(M)}_{t-s})} \right)  \geq - \beta e^{-U_s(\phi^{(M)}_{t-s})}
U_s(V_s \phi^{(M)}_{t-s}) I(s \leq t-1)ds \quad + \; \mbox{martingale increments}
\]
so that, 
\begin{eqnarray} 
Q^{\R^d, \beta,\gamma}_{\mu,f} \left[ U_t(B_R)=0 \right]
& = & \lim_{M \to \infty} Q^{\R^d, \beta,\gamma}_{\mu,f} \left[ e^{-M U_t(\psi_R)} \right] 
\nonumber \\
& \geq & e^{-\mu(\phi_t)} - \beta Q^{\R^d,\beta,\gamma}_{\mu,f} \left[
\int^{\infty}_0 U_s(V_s \phi_{t-s}) I(s \leq t-1) ds \right], \label{loex1}
\end{eqnarray}
where $\phi = \lim_{M \to \infty} \phi^{(M)}$. As in the extinction
results in section 3.1.4, a comparison argument shows that $\phi_t$ is bounded 
for all $t>0$. Moreover, then $\|\phi_t\|_{\infty} \to 0$ as $t \to \infty$, so that the term
$\exp(-\mu(\phi_t)) \to 1$. Similarly, for fixed $(s,\omega)$ we have
$U_s(V_s \phi_{t-s}) \to 0$ as $t \to \infty$. 
To justify the conclusion that
$Q^{\R^d,\beta,\gamma}_{\mu,f} [ U_t(B_R)=0 ] \to 0$ as $t \to \infty$ we will
dominate the second term on the right hand side of (\ref{loex1}) by
\[
Q^{\R^d,\beta,\gamma}_{\mu,f} \left[
\int^{\infty}_0 U_s(V_s \phi^*)  ds \right] = 
Q^{\R^d,\beta,\gamma}_{\mu,f} \left[
\int \phi^*(x) f(x) (1- e^{-U(0,\infty,x)}) dx \right] \leq \langle \phi^*,1 \rangle
\]
where $\phi^*(x) = \sup_{t \geq 1} \phi_t(x)$, and where
we are argue as in (\ref{part2}) for the first equality. But, for $\gamma >0$ and where
$T_t$ is convolution with the heat kernel $p_t(x)$, 
\begin{eqnarray*}
\phi^*(x) & \leq & \sup_{t \geq 0} e^{- \gamma t} T_{t} \phi_1(x) \\
& \leq & \sup_t \left\{ \int_{|y| \leq 1} \phi_1(x+y) e^{-\gamma t} p_t(y)  dy \right\}
+ \int_{|y| \geq 1} \phi_1(x+y) \sup_t \{e^{-\gamma t} p_t(y)\} dy \\
& \leq & \sup \{ \phi_1(z): |z-x| \leq 1 \} + C(\gamma) 
\int \phi_1(x+y) e^{-\gamma^{1/2} |y|} dy. 
\end{eqnarray*}
The fact that $ \langle \phi^*,1 \rangle < \infty$ can now be deduced from the boundedness and exponential decay of
$\phi_1$, which follow from comparison arguments as in \cite{dip} Lemma 3.1.

The stronger conclusion of the lemma follows by a similar argument using 
a slightly more complicated test function. Suppose that $\psi_R$ satisfies in addition
$\psi_R=1$ on $B_{R-1}$. Fix $t_0,t_1 \geq 0$ and redefine $\phi$ as the solution to 
\[
\partial_t \phi = \Delta \phi - \gamma \phi - \frac{\phi^2}{2} + M \psi_R \I(t \leq t_1) +
\beta \phi \psi_{R+2} \I(t \leq t_1+1), \qquad \phi_0 = 0.
\]
Some calculus shows that, under $Q^{\R^d,\beta,\gamma}_{\mu,f}$,
\begin{eqnarray*}
&& \hspace{-.3in} 
d( \exp( -U_s(\phi_{t_0+t_1-s}) - M U_{[t_0, t_0 \vee s]}(\psi_R)) ) \\
& = & \exp( -U_s(\phi_{t_0+t_1-s}) - M U_{[t_0, t_0 \vee s]}(\psi_R))
\left( - \beta U_s(V_s \phi_{t_0 +t_1-s}) + \beta U_s(\phi_{t_0 +t_1-s} \psi_{R+2})
\I(s \geq t_0-1) \right) dt 
\end{eqnarray*}
up to martingale increments. Taking expectations we have
\begin{eqnarray} 
&& \hspace{-.3in} Q^{\R^d,\beta,\gamma}_{\mu,f} \left[ e^{- M U_{[t_0, t_0+t_1]}(\psi_R)} \right]
\nonumber \\
& \geq & e^{-\mu(\phi_{t_0+t_1})} - \beta Q^{\R^d,\beta,\gamma}_{\mu,f} \left[
\int^{t_0+t_1}_0 U_s(\phi_{t_0+t_1-s} V_s)  - U_s(\phi_{t_0+t_1-s} \psi_{R+2}) \I(s \geq t_0-1) ds \right].
\label{loex2}
\end{eqnarray}
There is a unique non-negative solution $\overline{\phi}$ (see the arguments of \cite{iscoe}) to 
\[
\Delta \overline{\phi} = \gamma \overline{\phi}
 + \frac{\overline{\phi}^2}{2} - \beta \overline{\phi} \psi_{R+2} \quad \mbox{on $B^c_R$ and
$\overline{\phi}(x) \uparrow \infty$ as $|x| \downarrow R$.}
\]
A comparison argument shows that $\phi_t(x) \leq \overline{\phi}(x)$ for $ t \leq t_1+1$ and $|x|>R$. 
Also $\phi_{t_1+t} \leq \lambda_t$ for $t >0$ where 
\[
\dot{\lambda}_s = -\gamma \lambda_s - \frac{\lambda^2_s}{2} + \beta \lambda_s \I(s \leq 1) \quad 
\mbox{and $\lambda_s \uparrow \infty$ as $s \downarrow 0$.}
\] 
Moroever, as in the fixed $t$ argument, 
\[
\phi_{t_1+1+s} \leq e^{-\gamma s} T_s \phi_{t_1+1} \leq e^{-\gamma s} T_s (\overline{\phi}
\wedge \lambda_1) := \overline{\phi}_s.
\]
Using these bounds in (\ref{loex2}) we obtain for $t_0 >1$
\begin{eqnarray*}
&& \hspace{-.3in}
Q^{\R^d,\beta,\gamma}_{\mu,f} \left[ e^{- M U_{[t_0, t_0+t_1]}(\psi_R)} \right] \\
& \geq & e^{-\lambda_{t_0} \mu(1)} - \beta Q^{\R^d,\beta,\gamma}_{\mu,f} \left[
\int^{t_0-1}_0 U_s( V_s \overline{\phi}_{t_0-s-1}) ds
+ \int^{\infty}_{t_0-1} U_s (V_s \overline{\phi} \I(B^c_{R+1})) ds \right].
\end{eqnarray*}
Now we let $M,t_1 \to \infty$ to obtain the same upper bound for
$  Q^{\R^d,\beta,\gamma}_{\mu,f}[ U_{[t_0,\infty)}(B^c_R) =0]$. As $t_0 \to \infty$ we have
$\lambda_{t_0} \to 0$. The conclusion of the lemma follows as in the fixed $t$ result
using the domination $ \int^{\infty}_0 U_s( V_s \overline{\phi}^*) ds 
\leq \langle \overline{\phi}^*,1 \rangle $
where $\overline{\phi}^* = \max\{ \overline{\phi} \I(B^c_{R+1}), \sup_s \overline{\phi}_s \}$.
The integrability $ \langle \overline{\phi}^*,1 \rangle < \infty$ follows as in the fixed $t$ result using the 
exponential decay of the function $\overline{\phi}$. \qed
\subsection{Existence of the critical curve $\Psi(\beta)$} \label{s3.3}
For a DW($D,\gamma$) process the death time 
time $\tau=\inf \{t:U_t=0\}$ is finite, and 
then $U_{[0,\infty)} = U_{[0,\tau]}$ and $U^{\partial D}_{\infty} = U^{\partial D}_{\tau}$
are finite measures, $Q^{D,\gamma}_{\mu}$ almost surely.
Using the change of measure from section \ref{s2.2} the same holds for
the law $Q^{D,\beta,\gamma}_{\mu,f}$ when $D$ is a bounded domain. Indeed the 
change of measure martingale $M^{\beta,0,f}$ from (\ref{cofmM}) satisfies
\begin{eqnarray*}
\left[ M^{\beta,0,f} \right]_t & = & \beta^2 \int^t_0 \int_D f^2(x) e^{-2 U(0,s,x)} U_s(dx) ds \\
& = &  \frac{\beta^2}{2} \int_D f^2(x) \left(1- e^{-2U(0,t,x)} \right) dx
\end{eqnarray*}
by arguing as in the proof of (\ref{part2}). 
This is bounded independently of $t$, for a bounded domain $D$,
 so that $\mathcal{E}_t(M^{\beta,0,f})$ is a
uniformly integrable martingale, ensuring the above almost sure properties carry over
to the law $Q^{D,\beta,\gamma}_{\mu,f}$.

A similar argument leads to the following lemma, which  
shows that either certain death or possible life must occur, for each
pair of parameter values $\beta,\gamma$.
\begin{Lemma} \label{anymu}
If $Q^{\R^d,\beta,\gamma}_{\mu,1}[ \mbox{$U_t \neq 0$ for all $t$}]$ is strictly positive 
for some $\mu \in \mathcal{M}(\R^d)$ then it is strictly positive 
for all non-zero $\mu \in \mathcal{M}(\R^d)$.
\end{Lemma}

\noindent
\textbf{ Proof.} Step 1: We claim that  the laws of solutions
$ Q^{\R^d,\beta,\gamma}_{\mu,f}$ and $ Q^{\R^d,\beta,\gamma}_{\mu,g}$ are mutually absolutely continuous
on the entire sigma field $\mathcal{U}_{\infty} = \sigma\{U_t: t \geq 0\}$, provided that
$f$ and $g$ differ only on a compact set. To see this, consider the Radon-Nikodym derivative
between the two laws, which is the stochastic exponential 
$\mathcal{E}(M)$ arising from the martingale
$ M_t = \beta \int^t_0 \int_{\R^d} (f-g)(x) e^{- U_{[0,s]}(x)} M(dx,ds)$.  
Then, if $0\leq f,g\leq 1$ are supported inside the compact set $A$,
\begin{eqnarray}
 \left[ M \right]_t & = & \beta^2
\int^t_0 \int_{\R^d} (f-g)^2 (x) e^{-2 U_{[0,s]}(x)} U_s(dx) ds \nonumber \\ 
& \leq & \beta^2 \int^t_0 \int_{A} e^{-2 U_{[0,s]}(x)} U_s(dx) ds \nonumber \\ 
& = & \frac{\beta^2}{2} \int_A \left( 1- e^{- 2 U_{[0,t]}(x)} \right) dx 
\leq \frac{\beta^2}{2} |A|. \label{criticalbetabound}
\end{eqnarray}
The exponential martingale is therefore uniformly integrable  
which establishes the claim. 

Step 2: A result of Evans and Perkins \cite{evans+perkins} shows that,
for any fixed $t>0$, 
the laws $Q^{D,0}_{\mu}[U_t \in \cdot]$ and 
$Q^{D,0}_{\tilde{\mu}}[U_t \in \cdot]$ of two DW($\R^d,0$) processes started
at non-zero $\mu, \, \tilde{\mu}$ are mutually absolutely continuous on
$\mathcal{M}(\R^d)$. The change of measure ideas used to show
existence and uniqueness imply that the same absolute continuity holds for
$Q^{\R^d,\beta,\gamma}_{\mu,f}[U_t \in \cdot]$ and 
$Q^{\R^d,\beta,\gamma}_{\tilde{\mu},f}[U_t \in \cdot]$. 

Now we combine the two steps. Suppose $\mu, \, \tilde{\mu}$
are non-zero and compactly supported, 
and that started at $(\mu,1)$ the solutions may survive. Then
\begin{eqnarray*}
0 & < &  Q^{\R^d, \beta, \gamma}_{\mu,1} [ U_t \neq 0 \; \mbox{for all $t$}] \\
& = &
\int_{\mathcal{B}(\R^d,[0,1]) \times \mathcal{M}(\R^d)} 
Q^{\R^d, \beta, \gamma}_{\nu, g} [U_t \neq 0 \; \mbox{for all $t$}] \;
Q^{\R^d, \beta, \gamma}_{\mu,1} [U_{t_0} \in d \nu, \, f e^{-U(0,t_0)} \in dg]
\end{eqnarray*}
by the Markov property at time $t_0>0$. By the compact support
property of solutions, $f \exp(-U(0,t_0))$ is 
identically one outside a compact set, almost surely. So, using step one,  
\begin{eqnarray*}
0 & < & \int_{\mathcal{B}(\R^d,[0,1]) \times \mathcal{M}(\R^d)}  
Q^{\R^d, \beta, \gamma}_{\nu, 1} [U_t \neq 0 \; \mbox{for all $t$}]  \,
Q^{\R^d, \beta, \gamma}_{\mu,1} [U_{t_0} \in d \nu, \, f e^{-U(0,t_0)} \in dg] \\
& = & \int_{\mathcal{M}(\R^d)} Q^{\R^d, \beta, \gamma}_{\nu, 1} 
[U_t \neq 0 \; \mbox{for all $t$}] \,
Q^{\R^d, \beta, \gamma}_{\mu,1} [U_{t_0} \in d \nu].
\end{eqnarray*}
Now using step two above we obtain
$ 0 <  \int_{\mathcal{M}} Q^{\R^d, \beta, \gamma}_{\nu, 1} 
[U_t \neq 0 \; \mbox{for all $t$}] \,
Q^{\R^d, \beta, \gamma}_{\tilde{\mu},1} [U_{t_0} \in d \nu]$
and undoing the steps we find that the process with initial conditions $(\tilde{\mu},1)$ 
may survive.

To remove the restrictions that $\mu,\tilde{\mu}$ be compactly supported, recall 
that the solutions are compactly supported at any time $t>0$ almost surely. 
Hence, if solutions may live from some general $\mu \neq 0$, they may, by 
applying the Markov property at time $t>0$, live from some, and hence all, non-zero
compactly supported initial conditions. Any solution started at non-zero $\tilde{\mu}$ will have
positive probability of being alive and compactly supported at small times $t>0$ 
(use the continuity of the total mass process) and hence, by the Markov property again, may live for all
time. \qed

The next two lemmas give a characterization of certain death, which allows us to 
reduce the question of possible life/certain death to the study of the total
exit measures solutions on bounded domains. The existence of a 
non-decreasing critical curve follows
from this, since these exit measures are stochastically 
monotone in $\beta,\gamma$.
\begin{Lemma} \label{deathcharacterization1}
Suppose $\mu \in \mathcal{M}(\R^d)$ is compactly supported. Then 
$Q^{\R^d,\beta,\gamma}_{\mu,f}$ almost surely  
\[ 
\{ U_t = 0 \; \mbox{for large $t$} \} = 
\{ U_{[0,\infty)} \; \mbox{is compactly supported} \}. 
\]
\end{Lemma}

\noindent
\textbf{Proof.} The inclusion 
$\{ U_t = 0 \; \mbox{for large $t$} \} \subseteq 
\{ U_{[0,\infty)} \; \mbox{is compactly supported} \}$
follows from the compact support property 
(\ref{globaltsupport}). Let $D_L = (-L,L)^d$. 
We claim for any $L$, $P$ a.s., there exists a $t>0$ so that 
either $U_t(dx)=0$ or $U_t(D_L^c)>0$
must occur. This implies the opposite inclusion. 

Applying the first and second moment bounds (\ref{2ndmoments},\ref{12momentbounds}),
and the simple estimate $ P[Z>0] \geq (E[Z])^2/E[Z^2] $ for a non-negative 
variable $Z$, we may find $\epsilon_1>0$, depending only $L,\beta,\gamma$, so that
$ Q^{D,\beta,\gamma}_{\mu,f}[ U_1(D^c_L)>0 ] \geq \epsilon_1 $ for all
$f \in \mathcal{B}(D,[0,1])$ and  $\mu(1) \geq 1$. 
Using the extinction probability (\ref{extinctionprob}) we may find $\epsilon_2>0$,
depending only on $\beta,\gamma$, so that 
$Q^{D,\beta,\gamma}_{\mu,f}[U_1=0] \geq \epsilon_2$ for all 
$f \in \mathcal{B}(D,[0,1])$ and  $\mu(1) \leq 1$. 
Applying the Markov property at time $t$, we obtain
\[ 
Q^{\R^d,\beta,\gamma}_{\mu,f} 
\left[\mbox{$U_{t+1}(dx)=0$ or $U_{t+1}(D_L^c)>0$} \, | \, \mathcal{U}_t \right]
 \geq \epsilon_1 \wedge \epsilon_2.
\]
Iterating this estimate over the integers $t=1,2,\ldots$ 
and applying Borel-Cantelli completes the argument.
\qed.
\begin{Lemma} \label{deathcharacterization2}
Suppose $\mu \in \mathcal{M}$ is compactly supported in $D$. Then
\[
Q^{\R^d,\beta,\gamma}_{\mu,1}\left[ U_{[0,\infty)} (D^c) =0 \right]
=  Q^{D, \beta,\gamma}_{\mu,1} [U^{\partial D}_{\infty} = 0].
\]
\end{Lemma}

\noindent
\textbf{Proof.} We may 
construct a coupling of solutions $(u^-,v^-)$ and $(u^+,v^+)$ as follows:
$(u^-,v^-)$ is a solution to (\ref{spde}) on $D$ with initial conditions $(\mu,1)$;
$u^{\partial D,-}$ is the associated exit measure process;
$\tau = \inf\{t: u^{\partial D,-}_t >0\}$; and conditional on $\sigma\{u^-,v^-\}$
the process $(u^+,v^+)$ has the law of a solution on $\R^d$ started at the random 
initial condition $u^+_0 = u^-_{\tau}$ and $v^+_0 = v^-_{\tau}$.
One way to do this is to use 
the measurability of $(\mu,f) \to Q^{D,\beta,\gamma}_{\mu,f}$ and construct a
skew-product measure on the product space $\Omega_D \times \Omega_{\R^d}$.
We may also suppose that $v^-_t = \exp(-u^-(0,t))$ and 
$v^+_t= \exp(-u^-(0,\tau) - u^+(0,t))$ for all $t \geq 0$.  
Now define $(u_t,v_t) = (u^-_t,v^-_t)$ for $ t \leq \tau$ and
$(u_t,v_t) = (u^+_{t-\tau},v^+_{t-\tau})$ for $ t \geq \tau$. Then $v_t$ is a measurable function
of $(u_s:s \leq t)$ and this can be used to check that, with respect to the 
filtration $\sigma \{u_s:s \leq t\}$, the process $v$ is predictable and non-increasing. 
For $\phi \in C^2_b(\R^d)$, we may apply the extended martingale problem for
$(u^-,v^-)$ to the restriction of $\phi$ on $\overline{D}$. 
Combining with the martingale 
problem for $(u^+,v^+)$ we find that $(u,v)$ is a solution on $\R^d$
with initial conditions  $(\mu,1)$. Moreover on the set $\{\tau = \infty\}$ we have that
$u_{[0,\infty)}(D^c) = u^-_{[0,\infty)}(D^c)=0$. Thus
\[
Q^{\R^d,\beta,\gamma}_{\mu,1} \left[ U_{[0,\infty)}(D^c)=0 \right] 
 \geq  P \left[ \tau = \infty \right] 
 =  Q^{D,\beta,\gamma}_{\mu,1} \left[ U^{\partial D}_{\infty} =0 \right].
\]
For the converse inclusion we construct 
$(u^-,v^-)$ a solution to (\ref{spde}) on $\R^d$ with initial conditions $(\mu,1)$;
$\tau = \inf\{t: u^-_t(D^c) > 0\}$; and conditional on $\sigma\{u^-,v^-\}$
a process $(u^+,v^+)$ with the law of a solution on $D$ started at the random 
initial condition $u^+_0 = u^-_{\tau}$ and $v^+_0 = v^-_{\tau}$. 
A result of Perkins (see \cite{perkins1} Theorem III.5.1 and use absolute continuity) implies
that the paths of $t \to u_t^-(D^c)$ are almost surely continuous. So we may 
define $(u_t,v_t)$ as above and the process $u$ is almost surely continuous. 
A test function $\phi \in C^2_0(\overline{D})$ may be extended to $\tilde{\phi} \in C^2_b(\R^d)$.
Applying the martingale problem for
$(u^-,v^-)$ to $\tilde{\phi}$ and the martingale 
problem for $(u^+,v^+)$ to $\phi$, we find that $(u,v)$ is a solution on $D$
with initial conditions  $(\mu,1)$. Moreover, choosing $u^{\partial D}_t = u^{\partial D,+}_{(t-\tau)_+}$
on $t > \tau$ and $u^{\partial D}_t=0$ for $t \leq \tau$ one finds  
the extended martingale problem is solved. Then on the set $\{\tau = \infty\}$ we have that
$u^{\partial D}_{\infty} = 0$ and thus
\[
Q^{D,\beta,\gamma}_{\mu,1} \left[ U^{\partial D}_{\infty} =0 \right]
\geq P \left[ \tau = \infty \right] 
= Q^{\R^d,\beta,\gamma}_{\mu,1} \left[ U_{[0,\infty)}(D^c)=0 \right]
\]
completing the proof. \qed

Using the above results and the comparison results stated in section \ref{s2.2} we will now deduce 
the existence of an non-decreasing critical curve $\Psi$ as in Theorem \ref{mainresults}. 
\begin{Corollary}
There exists an non-decreasing function $\Psi(\beta):[0,\infty) \to [0,\infty]$ 
so that for $0 \leq \gamma < \Psi(\beta)$ possible life occurs and for 
$\gamma  > \Psi(\beta)$ certain death occurs.
\end{Corollary}

\noindent
\textbf{ Proof.} 
Lemmas \ref{anymu}, \ref{deathcharacterization1} and \ref{deathcharacterization2} 
show, for $\mu \neq 0$, that $Q^{\R^d,\beta,\gamma}_{\mu,1}[U_t =0 \; \mbox{for large $t \geq 0$}] =1 $
if and only if $Q^{\R^d,\beta,\gamma}_{\delta_0,1}[U_t =0 \; \mbox{for large $t \geq 0$}] =1 $
if and only if $Q^{\R^d,\beta,\gamma}_{\delta_0,1}[U_{[0,\infty)} \; \mbox{is compactly supported}] =1$
if and only if 
\[
\sup_n Q^{D_n,\beta,\gamma}_{\delta_0,1}[ U^{\partial D_n}_{\infty} =0] = 1.
\]
But Lemma \ref{comparison2} shows that these probabilities are non-increasing in $\beta$
and non-decreasing in $\gamma$. The result follows by setting
$ \Psi(\beta) = \sup \{ \gamma \geq 0: 
\sup_n Q^{D_n,\beta,\gamma}_{\delta_0,1}[ U^{\partial D_n}_{\infty} =0] < 1 \}$.  \qed

\noindent
\textbf{Remark.} 
The change of measure (\ref{cofmM}) can be used to show 
that $(\beta,\gamma) \to Q^{D,\beta,\gamma}_{\mu,1}[U^{\partial D}_{\infty} = 0]$ 
is continuous when $D$ is bounded, and hence that  
$Q^{\R^d, \beta,\gamma}_{\mu,1}[\mbox{$U_t =0$ for large $t$}]$ is lower semicontinuous
in $\beta,\gamma$. This 
does not seem to have any immediate implications for the continuity of $\Psi$.
\section{Approximations} \label{s4}
The proof of the decomposition results in section \ref{s2.2} 
is rather clear for particle systems that approximate
our reaction diffusion system. Our first proofs used a full particle approximation with
population and nutrient particles living on discrete lattices. We later realized that the key
to establishing the comparisons was to discretise the effect of the nutrient, and we here
present an approximation where the reaction with the nutrient occurs in a finite number of
discrete `packages' but the population still evolves as a continuum SPDE. This intermediate
approximation makes passage to the limit easier to establish. This passage to the limit 
is broadly similar to many in the literature, 
and the points of most interest are perhaps (i) that the nutrient interaction is singular
as a function of the population measure $u$, and requires tightness of the occupation densities;
(ii) that convergence of the exit measures follows quite simply from
convergence of the population measures, via the extended martingale problem (\ref{MGPD1+}).
\subsection{Construction of the approximation} \label{s4.1}
The approximations will depend on a parameter $N$, which we will suppress in the notation
in this section.
Fix a domain $D \subseteq \R^d$ and partition it as a disjoint union of sets
$D= \cup D_j$, where each set $D_j$ has diameter at most $N^{-1}$. 
Choose a finite number of functions $(\psi_k: 1 \leq k \leq K)$, each of the form 
$N^{-1} I(x \in D_j)$ for some $D_j$ in the partition. Each $\psi_k$ will represent a small 
package of nutrient which may be triggered to produce new population. 
The function $f = \sum_k \psi_k \leq 1$ will be our approximation 
to the initial nutrient level. 

Choose a probability space equipped with the following independent families:
an i.i.d. family of rate one exponential variables $e_k$ for $k \geq 1$;
independent DW($D,\gamma$) processes $(u_{k,t}:t \geq 0)$, for $k \geq 0$, with initial conditions
\[
u_{0,0} = \mu, \quad \mbox{and} \quad u_{k,0} = \beta \psi_k(x) dx \; \; 
\mbox{for $k \geq 1$.}
\]
Let $(u^{\partial D}_{k,t}:t \geq 0)$ be the associated exit measures processes.
Given realizations of these variables, the approximation can be constructed pathwise.
We will list as $\mathcal{S}_t \subseteq \{1,\ldots,K\}$ the 
labels of those nutrient packages that have been triggered by time $t$, 
and denote by $\tau_k \in [0,\infty]$ the time at which the $k$th nutrient package is triggered. 
Thus $\mathcal{S}_0 = \emptyset$ and $\mathcal{S}_t= \{k:\tau_k \leq t\}$. 
The approximation will satisfy
\begin{equation} \label{approximationdefn}
u_t = u_{0,t} + \sum_{k \in \mathcal{S}_t} u_{k,t - \tau_k}, \quad
u^{\partial D}_{t} = u^{\partial D}_{0,t} + \sum_{k \in \mathcal{S}_t} 
u^{\partial D}_{k,t - \tau_k}, \quad
v_t = \sum_{k \not \in \mathcal{S}_t} \psi_k.
\end{equation} 
Moreover we want the approximation process to be triggered in such a way that 
\[
\tau_k = \inf \left\{t: u_{[0,t]}(\hat{\psi}_k) > e_k \right\} \quad \mbox{for $k \geq 1$}
\]
where $\hat{\psi}_k = \psi_k/ \langle \psi_k,1 \rangle$.
This uniquely specifies a process, which
can be constructed pathwise by defining the triggering times $\tau_k$ 
and the processes $u_t,v_t$ inductively over the intervals between triggers
(and where, on a null set, packages may be simultaneously triggered). 

We now derive a martingale problem for the approximation.
On the intervals $[\tau_k,\tau_{k+1})$ the approximation process evolves as a finite sum of
DW processes. So it satisfies the following martingale problem: for $\phi \in C^2_b(\overline{D})$
\begin{equation} \label{AMGPD1-}
u_t(\phi) = \mu(\phi)  
+ \int^t_0 u_s \left( \Delta \phi - \gamma \phi \right)ds 
- u^{\partial D}_{t}(\phi) + \sum_{s \leq t} D_s u_s(\phi) + m_t(\phi), 
\end{equation}
where $m_t(\phi)$ is a continuous martingale with quadratic variation 
$\int^t_0 u_s(\phi^2) ds$ and where $D_s u_s (\phi)$ is the jump $u_s(\phi)- u_{s-}(\phi)$ at time $s$. 
We can compensate the jump term as
\begin{equation} \label{compensator}
\sum_{s \leq t} D_s u_s(\phi) = \beta \int^t_0 \sum_{k \not \in \mathcal{S}_s} \langle \phi,\psi_k \rangle 
u_s(\hat{\psi}_k) ds + E^{(1)}_t(\phi) 
\end{equation}
where $E^{(1)}_t(\phi)$ is a martingale with 
jumps bounded by $\beta N^{-d-1} \|\phi\|_{\infty}$ and with predictable
brackets process given by
\begin{eqnarray} 
d \langle E^{(1)}(\phi) \rangle_t &=&  \beta^2 \sum_{k \not \in \mathcal{S}_t}
\langle \phi,\psi_k \rangle^2 u_t(\hat{\psi}_k) dt \nonumber \\
& \leq & \beta^2 \| \phi\|_{\infty}^2 N^{-d-1} \sum_{k \not \in \mathcal{S}_t}
u_t(\psi_k) dt \leq 
\beta^2 \| \phi\|_{\infty}^2 N^{-d-1} u_t(1) dt. \label{error1}
\end{eqnarray}
The term $N^{-d-1}$ arises from bounding $\langle 1,\psi_k \rangle \leq N^{-d-1}$ (since $0 \leq \psi_k \leq N^{-1}$ and
the diameter of the support of $\psi_k$ is at most $N^{-1}$). Examining the variation of the function $\phi$ 
over the support of $\psi_k$, we can approximate
\[
\left|\langle \phi,\psi_k \rangle  u_s(\hat{\psi}_k) - u_s(\psi_k \, \phi ) \right| \leq 
 N^{-1} \| \nabla \phi \|_{\infty} u_s(\psi_k)
\]
and rewrite the compensator in (\ref{compensator}) as 
$\beta \int^t_0 u_s(v_s \phi) ds$ up to an error of size $O(N^{-1})$.
Combining these estimates we have
 \begin{equation} \label{AMGPD1}
u_t(\phi) = \mu(\phi)  
+ \int^t_0 u_s \left( \Delta \phi + \beta v_s \phi - \gamma \phi \right)ds 
- u^{\partial D}_{t}(\phi)  + m_t(\phi) +  E^{(1)}_t(\phi)+  E^{(2)}_t(\phi)
\end{equation}
where the second error term $E^{(2)}_t(\phi)$ is controlled by 
\begin{eqnarray} 
\left| E^{(2)}_t(\phi) - E^{(2)}_t(\phi) \right|  & = &
\left| \beta \int^t_s \sum_{k \not \in \mathcal{S}_r} \left( \langle \phi,\psi_k \rangle  
u_r(\hat{\psi}_k) -  u_r(\phi \psi_k) \right) dr \right| \nonumber \\
& \leq & \beta \left( (N^{-1}  \| \nabla \phi \|_{\infty}) \wedge \|\phi\|_{\infty} \right)  \int^t_s u_r(1) dr. \label{error2}
\end{eqnarray}
Taking $\phi =1 $ in (\ref{AMGPD1}) we find that
\[
u_t(1) \leq \mu(1) + (\beta-\gamma) \int^t_0 u_s(1) ds + m_t(1) + E^{(1)}_t(1).
\]
From this, via the standard argument (localizing, applying Burkholder's inequality to control the martingales and 
Gronwall's inequality), one can derive the moment bound, for any $p>0$ adnd $T < \infty$,
\begin{equation} \label{approxmoment} 
E \left[ \sup_{t \leq T} |u_t(1)|^p \right] \leq C(p,\beta,\gamma,T) \left( 1 + (\mu(1))^p \right). 
\end{equation}

We now consider the martingale decomposition of $t \to v_t(x)$. 
For $x \in D$ let $D_x$ be the partition element 
that contains $x$, let $\psi_x(y) = N^{-1} I(y \in D_x)$ and $\hat{\psi}_x = \psi_x / \langle \psi_x,1 \rangle$.
Then a nutrient package satisfies $\psi_k(x)>0$ only if $\psi_k=\psi_x$.
The process $t \to v_t(x)$ is a pure jump process with compensator given as in
\begin{eqnarray}
v_t(x) & = & f(x) - \int^t_0 \sum_{k \not \in \mathcal{S}_s} 
\psi_k(x) u_s(\hat{\psi}_k) ds + m^{v}_t(x) \nonumber \\
 & = & f(x) - \int^t_0 v_s(x)
u_s(\hat{\psi}_x)ds + m^{v}_t(x) \label{aveqn}
\end{eqnarray}
for a martingale $m^{v}_t(x)$ with jumps bounded by $N^{-1}$
and with predictable brackets process
\begin{equation} \label{error3}
d \langle m^{v}(x) \rangle_t = \sum_{k \not \in \mathcal{S}_t} 
\psi_k^2(x) u_t(\hat{\psi}_k) dt \leq N^{-1} v_t(x) u_t(\hat{\psi}_x) dt.
\end{equation}
We can solve the equation (\ref{aveqn}) for $v_t(x)$ as
\begin{equation} \label{approxvrep}
v_t(x) = e^{ - u_{[0,t]}(\hat{\psi}_x)}
\left( f(x) + \int^t_0 e^{- u_{[0,s]}(\hat{\psi}_x)} dm^{v}_s(x) \right).
\end{equation}
\subsection{Tightness of the approximations} \label{s4.2}
To pass to the limit in the martingale problems for the approximation we 
require tightness estimates. We tried deriving these from
known tightness estimates for each of the DW processes used as
building blocks, however the random start times caused some difficulties in
combining the estimates, so we proceed by repeating Sugitani's estimates for our 
discretized approximation. 

We now include the dependence on $N$ and write the 
approximation constructed in section \ref{s4.1} as $(u^{(N)},v^{(N)})$.
\begin{Proposition} \label{tightnessprop}
Let $D$ be a bounded domain. 
Suppose the approximation $(u^{(N)},v^{(N)})$ has initial conditions 
$u^{(N)}_0=\mu^{(N)}$ and $v^{(N)}_0 = f^{(N)} \leq 1$ satisfying 
that $\mu^{(N)} \to \mu$ weakly and 
$f^{(N)} \to f$ in $L^1(D)$ for some $f \in \mathcal{B}(D,[0,1])$.
Then
\begin{enumerate}
\item[(i)] the laws of $(u^{(N)})$ are tight in $D([0,\infty), \mathcal{M}(D))$ and limit points are continuous,
\item[(ii)] for any $\delta>0$, the laws of
the occupation densities $(u^{(N)}(\delta,t+\delta,x): t \geq 0, x \in D)$ are tight in
$C([0,\infty),C(D))$
\end{enumerate}
where $C(D)$ is the space of continuous
functions with the topology of uniform convergence on compacts. 
\end{Proposition}

\textbf{Proof}. 
The moment bound (\ref{approxmoment}) allows us the uniform control on total mass moments:
\begin{equation} \label{approxmoments}
E \left[ \sup_{t \leq T} |u^{(N)}_t(1)|^p \right] \leq  C(p,T,\beta,\gamma) \left(1 + (\mu^{(N)}(1))^p \right).
\end{equation}
To control the moments of the approximate densities, the 
following simple upper bound for $u^{(N)}$ is useful. Triggering all 
the nutrient packages $(u_k:1 \leq k \leq K_N)$ used in the construction of $u^{(N)}$
at time zero gives the process $\tilde{u}^{(N)}_t = \sum_{k=0}^{K_N} u_{k,t}$. 
This satisfies
\begin{equation} \label{tildecomp}
u^{(N)}(0,t,x) \leq \tilde{u}^{(N)}(0,t,x) \quad \mbox{for all $x \in D$ and $t \geq 0$.}
\end{equation}
Moroever $\tilde{u}^{(N)}$ has the law $Q^{D,\gamma}_{\mu^{(N)} + \beta f^{(N)} dx}$ of a 
DW($D,\gamma$) process.

Now fix $\phi \in C^2_0(\overline{D})$. Using the decomposition (\ref{AMGPD1}), we 
consider the increment $|u^{(N)}_t(\phi)-u^{(N)}_s(\phi)|$ for 
$0 \leq s <t$. This leads, via the total mass moment bounds and Burkholder's inequality, 
to the estimate
\begin{equation} \label{tight20}
E \left[ \sup_{r \in [s,t]} \left|u^{(N)}_r(\phi)-u^{(N)}_s(\phi)\right|^p \right] 
\leq  C \left( |t-s|^{p/2} + N^{-(d+1)p} \right) \quad \mbox{for all $N$ and $s,t \in [0,T]$,}
\end{equation}
where $C< \infty$ depends only on $\beta,\gamma,\phi,p,T$ and $\sup_N \mu^{(N)}(1)$ 
(we suppress this dependence below). 
The term $N^{-(d+1)p}$ comes from the size of the jumps in the Burkholder inequality, for a cadlag martingale $M$, 
in the form $E[\sup_{s \leq t} |M_s|^p] \leq C_p E[[M]_t^{p/2}] + C_p E[\sup_{s \leq t} |D_sM|^p]$. 
The estimate (\ref{tight20}) gives, via a Chebychev inequality,
\[
P \left[ \sup_{s \in [t,t+N^{-(d+1)}]} 
\left| u^{(N)}_s(\phi) - u^{(N)}_t(\phi) \right| \geq N^{-1/2} \right]
\leq C N^{-dp/2}.
\] 
By summing over the grid $t_j = j N^{-(d+1)}$ and applying Borel Cantelli, we find that
\begin{equation} \label{final5}
\left| u^{(N)}_s(\phi) - u^{(N)}_t(\phi) \right| \leq 2N^{-1/2} \quad \mbox{
whenever $0 \leq s,t \leq T$ and $|s-t| \leq N^{-(d+1)}$,
for all large $N$, $P$ a.s.} 
\end{equation}

We claim that (\ref{tight20}) is sufficient to
imply that $u^{(N)}(\phi)$ is tight in $D([0,\infty),\R)$ and has
continuous limit points. To see this we smooth the sample paths by setting
$ \overline{u}^{(N)}_t = N^{d+1} u^{(N)}_{[t,t+N^{-(d+1)}]}$.
Then $t \to \overline{u}^{(N)}(\phi)$ is continuous and
\[
\left| \overline{u}^{(N)}_t(\phi) - \overline{u}^{(N)}_t(\phi) \right|
\leq 2 N^{(d+1)} |t-s| \|\phi\|_{\infty} \sup_{r \leq T+1}u^{(N)}_r(1)
\quad \mbox{for all $N$ and $s,t \in [0,T]$.}
\]
This together with (\ref{tight20}) implies that
$ E [ | \overline{u}^{(N)}_t(\phi)- \overline{u}^{(N)}_s(\phi)|^p] 
\leq C |t-s|^{p/2}$ 
for $s,t \in [0,T]$. Hence the laws of $ \overline{u}^{(N)}(\phi)$ 
are tight in $C([0,T],\R)$. But (\ref{final5}) shows that the error
$ \sup_{t \leq T} |u^{(N)}_t(\phi) - \overline{u}^{(N)}_t(\phi)| \to 0 $ almost surely
and $u^{(N)}_t(\phi)$ has the the same continuous limit points as
$\overline{u}^{(N)}_t(\phi)$. 

We now state a compact containment condition. Let
$D_{\epsilon} = \{x \in D: d(x,D^c) > \epsilon\}$. Choose smooth
$\phi_{D,\epsilon}:D \to [0,1]$ satisfying $\phi_{D,\epsilon}=0$ on $D_{2 \epsilon}$
and $\phi_{D,\epsilon} = 1$ on $D \setminus D_{\epsilon}$. 
We also suppose that $\phi_{D,\epsilon}$ is decreasing in $\epsilon$.
We claim that for all $\delta >0$ there exists $\epsilon(\delta)>0$ satisfying
\begin{equation} \label{compactcontainment}
\lim_{\delta \to 0} \limsup_{N \to \infty} 
P \left[ \sup_{t \leq T} u^{(N)}_t( \phi_{D,\epsilon(\delta)}) \geq \delta \right] = 0. 
\end{equation}
We postpone the proof of this to the end of the section. 
This condition controls the amount of mass near the boundary and  
allows us to extend the tightness of $u^{(N)}_t(\phi)$ above to all $\phi \in C^0_b(\overline{D})$ by 
an approximation argument. It also is the compact containment hypothesis in 
\cite{perkins1} Theorem II.4.1, from which we obtain  
the tightness of $u^{(N)}$ in $D([0,\infty),\mathcal{M}(D))$ and the 
continuity of the limit points.

To establish the tightness of the occupation densities we use the following 
increment estimate.
Under the hypotheses of Proposition \ref{tightnessprop}, we claim there exist 
$0<\alpha_i< \infty$ and finite
$C$, depending on $\beta,\gamma,p,\delta,\delta',T,D$ and $\sup_N \mu^{(N)}(1) + (f^{(N)},1)$, so that
\begin{equation} \label{densitytightness}
E \left[ |u^{(N)}(\delta,s,x)-u^{(N)}(\delta,t,y)|^p \right] 
\leq C \left(|t-s|^{\alpha_1 p} + 
|x-y|^{\alpha_2 p} \right) 
\end{equation}
for all $\delta \leq s,t \leq T$ and $x,y \in D \setminus D_{\delta'}$.
This implies tightness of $(t,x) \to u^{(N)}(\delta,t,x)$ in 
$C([\delta,\infty), C(D))$, where $C(D)$ is the space of continuous functions with
the topology of uniform convergence on compacts. 

To establish (\ref{densitytightness}) we use a Green's function representation for the density
(note the density exists since it exists for each component part $u_{k,t}$).
It is convenient to break the approximation into two parts
\begin{equation} \label{densitysplit}
u^{(N)}_t = u_{0,t} + \hat{u}^{(N)}_t \quad \mbox{where} \quad \hat{u}^{(N)}_t = \sum_{k \in \mathcal{S}_t} u_{k,t-\tau_k}. 
\end{equation}
For the process
$\hat{u}^{(N)}$ we may consider $\hat{u}^{(N)}(0,t,x)$ since the nutrient packages do not start 
at a singular initial condition. Indeed
combining the Green's function representation (\ref{densityrep}) for each $u_{k,t-\tau_k}$ and summing 
over $k \in \mathcal{S}_t$ we obtain
\begin{equation} \label{densityapproxrep}
\hat{u}^{(N)}(0,t,x) = \int^t_0 \int_D G^{D,x}_{[0,t-s]}(z) \left(-\gamma \hat{u}^{(N)}_s(dz) + \hat{m}^{(N)}(dz,ds)\right) 
+ \sum_{s \leq t} D_s u_s(G^{D,x}_{[0,t-s]})
\end{equation}
where $\hat{m}^{(N)}$ is the martingale measure associated to $\hat{u}^{(N)}$. 
We now consider increments of this representation. There are many terms, so we choose to illustrate the 
idea on one key noise term and on the jump term.
We shall also state all the underlying estimates on the kernel $G^{D,x}_{[0,t]}$. First
\[
G^{D,x}_{[0,t]}(z) \leq G^{\R^3,x}_{[0,t]}(z) \leq C|x-z|^{-1} \quad \mbox{for $t >0$ and $x,z \in \R^3$}
\]
and the analogous bounds $G^{\R^2,x}_{[0,t]}(z) \leq C(T) (1 \vee \ln(1/|x-z|)$ and 
$G^{\R,x}_{[0,t]}(z) \leq C(T)$ when $t \leq T$. Also
\[
|G^{D,x}_{[0,t]}(z) - G^{D,x}_{[0,t']}(z)|  \leq  C(d,\delta) |t-t'| |x-z|^{-d} 
\]
for all $t,t'$ and $x,z \in D$ and 
\[
|G^{D,x}_{[0,t]}(z) - G^{D,y}_{[0,t]}(z)| \leq  C(d,\delta,T) |x-y| (1+|x-z|^{1-d})
\]
for $t,t' \leq T$ and $x,y \in D_{\delta}$, $z \in D$ satisfying $|x-z| \geq 2|x-y|$.
These estimates are explicit calculations when $D = \R^d$. One way to obtain them for 
general $D$ is via a suitable coupling argument. 

Consider the term $K(t,x) = \int^t_0 \int G^{D,x}_{[0,t-s]}(z) \hat{m}^{(N)}(dz,ds)$. The time 
increment $K(t+s,x) - K(t,x)$ can be split into two terms, defined by 
\[
K_1 = \int^t_0 \int G^{D,x}_{[0,t+s-r]}(z)- G^{D,x}_{[0,t-r]}(z) \hat{m}^{(N)}(dz,dr), \quad
K_2 = \int^{t+s}_t \int G^{D,x}_{[0,t-r]}(z) \hat{m}^{(N)}(dz,dr).
\]
Using a Burkholder inequality, and the kernel estimates above we find, when $d=3$,
\[
E[|K_1|^p] \leq C(p) E \left[ \left| \int_D \hat{u}^{(N)}(0,t,z) \left(|z-x|^{-2} \wedge s^2 |z-x|^{-6} \right) 
dz \right|^{p/2} \right].
\]
A comparison as in (\ref{tildecomp}) combined with the moments (\ref{densitymoment2}) shows that
$\sup_z E[\hat{u}^{(N)}(0,t,z)^{p/2}] \leq C(p,\beta,T)$ when $t \leq T$ and $z \in D$. This leads to 
$E[|K_1|^p] \leq C(p,\beta,T,D) s^{p/4}$. Similarly, again when $d=3$ and $t+s \leq T$, apply Holder's inequality
to find 
\begin{eqnarray*}
E[|K_2|^p] & \leq & C(p) E \left[ \left|\int^{T}_0 \I(r \in [t,t+s]) \hat{u}^{(N)}_r(|x-\cdot|^{-2}) \right|^{p/2} \right] \\
& \leq &  C(p) E \left[ \left|\int^{T}_0 \I(r \in [t,t+s]) \hat{u}^{(N)}_r(dz) \right|^{p/8}  
\; \left|\int^{T}_0 \hat{u}^{(N)}_r(|x-\cdot|^{-8/3}) \right|^{3p/8} \right] \\
& \leq & C(p) |s|^{p/8} E\left[ | \sup_{t \leq T} \hat{u}^{(N)}_t(1)|^{p/8}  \left| 
\int \hat{u}^{(N)}(0,T,z) |x-z|^{-8/3} dz \right|^{3p/8} \right] \\
& \leq & C(p,\beta,T,D) |s|^{p/8}.
\end{eqnarray*}
In this second estimate we do not have the optimal power of $s$, but have shown how to rely only on the 
simple moments estimates for $\sup_{t \leq T} \hat{u}^{(N)}_t(1)$ and $u^{(N)}(0,T,z)$.
 
The jump term in the representation (\ref{densityapproxrep}) can also be easily handled, since the sum of all jumps
in the measure $t \to u_t$ is at most $\beta dx$. For example, an increment in $x$
is bounded using the kernel estimates above, when $d=3$ and $|x-y| \leq 1/2$,
\begin{eqnarray*}
&& \hspace{-.4in} \left| \sum_{s \leq t} D_s u_s(G^{D,x}_{[0,t-s]}) - D_s u_s(G^{D,y}_{[0,t-s]}) \right| \\
& \leq &  \beta \int_D \sup_{s \leq t} |G^{D,x}_{[0,t-s]}(z) - G^{D,y}_{[0,t-s]}(z) | dz \\
& \leq & C \beta \int_{B(x, |x-y|^{1/2})} |x-z|^{-1} + |y-z|^{-1} dz + 
C \beta |x-y| \int_{D \setminus B(x, |x-y|^{1/2})}
1+ |x-z|^{-2} dz \\
& \leq & C(\beta,D) |x-y|.
\end{eqnarray*}
Other terms in the Green's function representation of the increment $\hat{u}(0,t,x)$, and in all dimensions $d=1,2,3$, 
can be controlled by similar estimates. 

From (\ref{densitysplit}) we have  $u^{(N)}(\delta,t,x) = u_0(\delta,t,x) + \hat{u}^{(N)}(\delta,t,x)$ and we have 
controlled the increments of $\hat{u}^{(N)}(0,t,x)$, and hence of 
$\hat{u}^{(N)}(\delta,t,x) = \hat{u}^{(N)}(0,t,x) - \hat{u}^{(N)}(0,\delta,x)$. The increment estimates for 
$u_0(\delta,t,x)$ are similar, except that one restricts to $\delta >0$ and uses
(\ref{densitymoment1}). One also needs estimates on the fixed time measure $u_{0,\delta}$ to control
the increments of $u_{0,\delta}(G^{D,x}_{[0,t-\delta]})$. For this, use the method
of Lemma III.3.6 in Perkins \cite{perkins1} to show, for $p \in [0,2)$ and
suitable $c_0=c_0(p)>0$, that $E[\exp(c_0 u_{0,\delta}(|x-\cdot|^{-p})]
\leq C(d,\sup_n \mu_n(1), \delta,p)< \infty$, which is sufficient. \qed

\noindent
\textbf{Proof of the compact containment (\ref{compactcontainment}).}
We start by showing that
\begin{equation} \label{final1}
\sup_N E \left[  u^{(N)}_{[0,T]}( \phi_{D,\epsilon}) \wedge 1 \right] \to 0 \quad \mbox{as $\epsilon \downarrow 0$.} 
\end{equation}
Indeed, arguing as in (\ref{massnrbdy}) we have
\begin{equation} \label{partial}
Q^{\R^d,\gamma}_{\mu} \left[ U_{[0,T]}(\phi_{D,\epsilon}) \wedge 1 \right]
\downarrow 0 \quad \mbox{as $\epsilon \downarrow 0$.}
\end{equation}
Moreover the continuity of the laws $\mu \to Q^{\R^d,\gamma}_{\mu}$ on 
$C([0,\infty), \mathcal{M}(\R^d))$ ensures that 
the expectation in (\ref{partial}) is a continuous function of $\mu$. 
Therefore the limit in (\ref{partial}) is uniform over the compact set $K_{D,L}$ of measures 
supported in $\overline{D}$ and with total mass $\mu(1) \leq L$.  
By the coupling (\ref{tildecomp}) and the simple comparison (\ref{simplecoupling})
we find 
\[
E \left[ u^{(N)}_{[0,T]}( \phi_{D,\epsilon}) \wedge 1 \right] 
\leq 
Q^{D,\gamma}_{\mu^{(N)}+ \beta f^{(N)}dx} \left[ U_{[0,T]}( \phi_{D,\epsilon}) \wedge 1 \right] 
\leq 
Q^{\R^d,\gamma}_{\mu^{(N)}+ \beta f^{(N)}dx} \left[ U_{[0,T]}( \phi_{D,\epsilon}) \wedge 1 \right] 
\]
which leads to (\ref{final1}). 

To control the supremum $\sup_{t \leq T} u^{(N)} (\phi_{D,\epsilon})$, we will 
establish a modulus of continuity that is uniform in $N$ and 
$\epsilon$. Let $G^{D}_t\phi$ denote the action of the Green's function as defined by 
$G^{D}_t\phi(x) = \int_D G^{D,x}_t(y) \phi(y) dy$.
Then, arguing as above, we have
\begin{equation} \label{final2}
\sup_N \sup_{t \leq T} \mu^{(N)}(G^D_t \phi_{D,\epsilon}) \to 0 \quad \mbox{as $\epsilon \to 0$.} 
\end{equation}
For $\phi \in \mathcal{B}(D,[0,1])$, the measurable functions from
$\overline{D}$ to $[0,1]$, define
\[
X_t(N,\phi) = u^{(N)}_t(\phi) - \mu^{(N)}(G^{D}_t \phi) \quad \mbox{and} \quad
\overline{X}_t(N,\phi) = N^{d+1} \int_t^{t + N^{-(d+1)}} X_s(N,\phi)ds.
\]
We may extend the decomposition (\ref{AMGPD1}) to time dependent test functions 
as in section 3.1.1.  For smooth $\phi$ compactly supported inside $D$, using the test function
$(s,x) \to G^D_{t-s}\phi(x)$ in (\ref{AMGPD1}), over the interval $s \in [0,t]$,
leads to the Green's function representation for $X_t(N,\phi)$. This representation in turn 
can be used to reach an increment estimate
\[
E[ | X_t(N,\phi) - X_s(N,\phi) |^p ] \leq C \left( |t-s|^{p/2} + N^{-(d+1)p} \right) 
\quad \mbox{for all $N$ and $s,t \in [0,T]$,}
\]
where $C$ depends on $p,T$ but not on $\phi$ or $N$. This requires 
(see Perkins \cite{perkins3} Corollary 5 for this argument) 
only a smoothing property of the Green's kernel, namely we use:
\[
\sup_{x \in D} \int_D | G^{D,x}_u(y) - G^{D,x}_t(y) | dy \leq C(D) t^{-1}(u-t) \quad 
\mbox{for all $0 < t < u$.}
\]
This estimate holds by direct calculation for the case where $D=\R^d$, and thereby when $D$ is a box, since 
then the Green's kernel is a finite combination of reflected copies of the free space Green's kernel. 
For domains with a smooth boundary, it can be derived from Corollary 5 of Davies \cite{davies}.

The smoothed version $\overline{X}_t(N,\phi)$ then satisfies 
$ E[ | \overline{X}_t(N,\phi) - \overline{X}_s(N,\phi) |^p ] \leq C |t-s|^{p/2}$
for all $N$ and $s,t \in [0,T]$. The argument from Perkins \cite{perkins3} 
Theorem 1 then shows that for all 
$\phi \in \mathcal{B}(D,[0,1])$, there exists a random $c(\phi,N,T)>0$ so that 
\begin{equation} \label{final3}
\left| \overline{X}_t(N,\phi) - \overline{X}_s(N,\phi) \right| \leq |t-s|^{1/3} 
\quad \mbox{whenever $s,t \in [0,T]$ and $|t-s| \leq c(\phi,N,T)$}
\end{equation}
and moreover that $P[c(\phi,N,T) \leq \delta] \to 0$ as $\delta \downarrow 0$
uniformly over $N$ and $\phi \in \mathcal{B}(D,[0,1])$. The point here is that
the estimates on $c(\phi,N,T)$ depend only on total moment bounds and on 
$\|\phi\|_{\infty}$.  

It is straightforward to combine (\ref{final1}), (\ref{final2}) and (\ref{final3})
to see that there exists $\epsilon(\delta)>0$ so that
\begin{equation} \label{final4}
\lim_{\delta \to 0}
\sup_N P \left[\sup_{t \leq T} \overline{u}^{(N)}_t(\phi_{D,\epsilon(\delta)}) \geq 2 \delta \right] =0.
\end{equation}
Indeed, consider the set 
\[
\left\{ \sup_{t \leq T} \overline{u}^{(N)}_t(\phi_{D,\epsilon}) \geq 2 \delta, \;
 c(\phi_{D,\epsilon},N,T+1) \geq \delta^3 \right\}.
\]
If $\delta >1$ and $ \sup_N \sup_{t \leq T+1} \mu^{(N)}(G^D_t(\phi_{D,\epsilon})) \leq \delta^4/8$, 
then on this set
we have $\sup_{t \leq T} \overline{X}_t(N,\phi_{D,\epsilon}) \geq \delta$
and the modulus of continuity ensures that 
$\int^T_0 \overline{X}_t(N,\phi_{D,\epsilon}) dt \geq \delta^4/4$, and thence that 
$\int^T_0 \overline{u}^{(N)}_t(\phi_{D,\epsilon}) dt \geq \delta^4/8$. But this 
has small probability by (\ref{final1}) and Markov's inequality. So (\ref{final4}) follows
by first choosing $\delta$ so that that  
$P [c(\phi_{D,\epsilon},N,T+1) \leq \delta^3]$ is small, uniformly in 
$N$ and $\epsilon$, and then choosing $\epsilon$ so that 
$ \sup_N \sup_{t \leq T+1} \mu^{(N)}(G^D_t(\phi_{D,\epsilon}) \leq \delta^4/8$
and $\sup_N E [  \overline{u}^{(N)}_{[0,T]}( \phi_{D,\epsilon}) \wedge 1 ]$ is small.

To remove the time smoothing in \ref{final4}), we would like to apply (\ref{final5}). However, since 
$\phi_{D,\epsilon}$ does not lie in $C^2_0(\overline{D})$, there is an extra exit measure term
when deriving the increment estimate (\ref{tight20}). However the exit measure is non-negative and
we can still deduce the one sided estimate $ u^{(N)}_t(\phi) \leq u^{(N)}_s (\phi) +  2N^{-1/2} $
whenever $0 \leq s < t \leq T$ and $|t-s| \leq N^{-(d+1)}$, for all large $N$, $P$ a.s,
and this is sufficient to deduce from (\ref{final4}) the desired compact containment 
estimate. \qed
\subsection{Passage to the limit} \label{s4.3}
The aim now is to prove Theorem \ref{convergence} below, showing that the approximations 
constructed in section \ref{s4.1} converge in law to solutions of
(\ref{spde}). The method is to pass to the limit in the approximate 
martingale problem (\ref{AMGPD1}) and confirm that any limit point of the 
sequence of approximations must be a solution to (\ref{spde}). Many steps in the argument
are quite standard (see, for example,
the arguments in Proposition II.4.2 in \cite{perkins1}), so we
concentrate on those that concern the interaction term and the exit measures.

Since we wish to let $N \to \infty$, we again include this dependence in the 
notation. Assume the hypotheses of the tightness Proposition \ref{tightnessprop}. 
For any subsequence, we may
choose a sub-subsequence where the approximations and their occupation densities converge
in distribution. 
By a Skorokhod embedding, we may choose versions of these
approximations that converge almost surely. Without changing 
notation for the version and labelling the convergent sub-subsequence 
still as $N$, we may suppose there exists a limit process
$u$, with continuous paths in $\mathcal{M}(D)$ and a continuous field 
$(u(s,t,x): 0 < s \leq t, x \in D)$ so that, almost surely,  
\begin{eqnarray}
\sup_{s \leq t} \left| u^{(N)}_s(\phi) - u_s(\phi) \right| \to 0, 
&& \mbox{for all $\phi \in C^0_b(\overline{D})$ and $t< \infty$,} \label{conv1} \\
\sup_{x \in A, s \in [\delta,t]} \left|u^{(N)}(\delta,s,x) - u(\delta,s,x) \right| \to 0, &&
\mbox{for all $0<\delta \leq t< \infty$ and compact $A \subseteq D$.} \label{conv2}
\end{eqnarray}
Note that $u(s,t,x)$ must act as a density for $u_{[s,t]}(dx)$ for $0< s \leq t$. Since
also $s \to u(s,t,x)$ is non-increasing we may set
$ v_t(x) = f(x) \exp(-u(0,t,x))= \lim_{s \downarrow 0} f(x) \exp(-u(s,t,x))$. 

Fix $\phi \in C^2_0(\overline{D})$. Consider the 
martingale problem (\ref{AMGPD1}) for $u^{(N)}$ at a fixed $t \geq 0$, noting that 
$u^{\partial D,(N)}_{t}(\phi)=0$. 
The terms $u^{(N)}_t(\phi)$, $\mu^{(N)}(\phi)$ and 
$\int^t_0 u^{(N)}_s(\Delta \phi-\gamma\phi)ds$ converge pathwise
as $N \to \infty$ by (\ref{conv1}).
The error terms $ E^{(1,N)}_t(\phi)$ and $E^{(2,N)}_t(\phi)$ converge to zero 
(in $L^2$ and pathwise respectively) using the estimates (\ref{error1}) and (\ref{error2})
and the moment bounds (\ref{approxmoments}). 
To handle the key term $\int^t_0 u^{(N)}_s(v^{(N)}_s \phi) ds$ we approximate, using (\ref{approxvrep}),
\[
\int^t_0 u^{(N)}_s(v^{(N)}_s \phi) ds = \int^t_0 u^{(N)}_s \left(f^{(N)} 
e^{-u^{(N)}_{[0,s]}(\hat{\psi}_{\cdot})} \phi \right) 
ds + E^{(3,N)}_t(\phi)
\]  
where the error term is bounded by 
\begin{equation}
\left| E^{(3,N)}_t(\phi) \right| = 
\left| \int^t_0 u^{(N)}_s \left( \phi e^{-u^{(N)}_{[0,s]}(\hat{\psi}_{\cdot})}
 \tilde{m}^{v,N}_s \right) ds \right|
\leq \|\phi\|_{\infty}  \int_D u^{(N)}(0,t,x) \, \sup_{s \leq t} 
\tilde{m}^{v,N}_s(x) \, dx \label{temp212}
\end{equation}
and $\tilde{m}^{v,N}_t(x) = \int^t_0 \exp ( - u^{(N)}_{[0,s]}(\hat{\psi}_x))
dm^{v,N}_s(x)$. Doob's inequality, and the estimate (\ref{error3}), shows that
\[ 
E \left[ \sup_{s \leq t} |\tilde{m}^{v,N}_t(x)|^2 \right] \leq 2 N^{-1} \int^t_0 
E \left[ e^{- 2 u^{(N)}_{[0,s]}(\hat{\psi}_x)} u^{(N)}_s (\hat{\psi}_x)  \right] ds
\leq N^{-1}. 
\]
The comparison (\ref{tildecomp}) shows that 
$ E \left[(u^{(N)}(0,t,x))^2\right] \leq 
Q^{D,\gamma}_{\mu^{(N)}+\beta \I_D dx} 
\left[ (U(0,t,x))^2 \right].$
Applying Cauchy-Schwarz to (\ref{temp212}), and using the the estimate 
(\ref{densityatzero}),
we find that $E(|E^{(3,N)}_t(\phi) |) \to 0$ as $N \to \infty$. 

For $\delta \in (0,t)$ we write 
\begin{eqnarray}
&& \hspace{-.3in} \int^t_0 u^{(N)}_s \left(f^{(N)} e^{-u^{(N)}_{[0,s]} 
(\hat{\psi}_{\cdot}) } \, \phi \right) ds  \nonumber \\
& = & \int^t_{\delta} u^{(N)}_s \left(f e^{-u^{(N)}_{[\delta,s]} 
(\hat{\psi}_{\cdot})} \, \phi \right) ds
+ \int^t_{\delta} u^{(N)}_s \left((f-f^{(N)}) e^{-u^{(N)}_{[\delta,s]} 
(\hat{\psi}_{\cdot})} \, \phi \right) ds  \nonumber\\
& & + \int^{\delta}_0 u^{(N)}_s \left(f^{(N)} e^{-u^{(N)}_{[0,s]} 
(\hat{\psi}_{\cdot})} \, \phi \right) ds 
 + \int^t_{\delta} u^{(N)}_s \left( f^{(N)} 
( e^{-u^{(N)}_{[0,s]}(\hat{\psi}_{\cdot})} - e^{-u^{(N)}_{[\delta,s]}
(\hat{\psi}_{\cdot})} ) \phi \right) ds \nonumber \\
& = & \int^t_{\delta} u^{(N)}_s \left(f e^{-u^{(N)}_{[\delta,s]} 
(\hat{\psi}_{\cdot})} \, \phi \right) ds + E^{(4,N)}_t(\phi)
+E^{(5,N)}_t(\phi) + E^{(6,N)}_t(\phi).
\label{deltaerror}
\end{eqnarray}
The error term $E^{(5,N)}_t(\phi)$ is bounded by
$\| \phi \|_{\infty} \int^{\delta}_0 u^{(N)}_s(1)ds$ and converges to zero,
uniformly in $N$, as $\delta \to 0$. 
The error term $E^{(6,N)}_t(\phi)$ is bounded by 
\[
\|\phi\|_{\infty} \int^t_0 u^{(N)}_s (u^{(N)}_{[0,\delta]} (\hat{\psi}_{\cdot}) ) ds
=  \|\phi\|_{\infty} \int_D u^{(N)}(0,t,x) u^{(N)}_{[0,\delta]}(\hat{\psi}_x)  dx.
\]
The approximate density $u^{(N)}_{[0,\delta]}(\hat{\psi}_x)$ satisfies
 $\int_D (u^{(N)}_{[0,\delta]}(\hat{\psi}_x))^2 dx
\leq \int_D (u^{(N)}(0,\delta,x))^2 dx$. This, together with
Cauchy-Schwarz as above, shows that 
$E^{(6,N)}_t(\phi)$ converges to zero in $L^1$, uniformly in $N$, 
as $\delta \to 0$. 
The error term $E^{(4,N)}_t(\phi)$ is bounded by
\[
\| \phi\|_{\infty} \int^t_0 u^{(N)}_s \left(|f^{(N)}- f| \right) 
= \|\phi\|_{\infty} \int_D u^{(N)}(0,t,x) |f^{(N)}(x)-f(x)| dx.
\]
Since $f^{(N)}$ are bounded we have $f^{(N)} \to f$ in $L^2(D)$ and combined with
(\ref{densityatzero}) we see that $E^{(4,N)}_t(\phi) \to 0$ in $L^2$. 

Finally, the first term on the right hand side of (\ref{deltaerror}) 
converges pathwise, for fixed $\delta$, 
to $\int^t_{\delta} u_s ( f \exp(-u(\delta,s)) \phi )$ as 
$N \to \infty$. 
This follows  from (\ref{conv1},\ref{conv2}) when $f \in C^0_b(D)$, 
and for general measurable $f$ we can approximate in $L^2$ by continuous $\tilde{f}$ and 
control the error as for the error term $E^{(4,N)}_t(\phi)$.
The limit $\int^t_{\delta} u_s( f \exp(-u(\delta,s)) \phi) ds$ is, by 
repeating the above approximations, close, for small $\delta$, to $\int^t_0 u_s(v_s \phi) ds$. 
The conclusion is that the key term $\int^t_0 u^{(N)}_s(v^{(N)}_s \phi) ds$ converges,
in probability, to $\int^t_0 u_s(v_s \phi) ds$. 

We may now define a continuous 
process $m_t(\phi)$, for $\phi \in C^2_0(\overline{D})$, by the 
formula (\ref{MGP1}). The convergence of all other terms ensures that
$m^{(N)}_t(\phi) \to m_t(\phi)$ in probability. Moreover, standard arguments
yield that $m_t(\phi)$ is a martingale with respect to $\sigma \{u_s:s \leq t\}$,
and with the correct quadratic variation (note that we have uniform moment control by
(\ref{approxmoments})).
Thus the limiting process $(u,v)$ is a solution to the martingale problem for 
(\ref{spde}) started from $(\mu,f)$ with respect to its natural filtration
$\sigma \{u_s:s \leq t\}$. We have almost proved the following convergence result.
\begin{Theorem} \label{convergence}
Let $D$ be a bounded domain. 
Suppose the approximation $(u^{(N)},v^{(N)})$ has initial 
$u^{(N)}_0=\mu^{(N)}$ and $v^{(N)}_0 = f^{(N)} \leq 1$ satisfying 
that $\mu^{(N)} \to \mu$ weakly and 
$f^{(N)} \to f $ in $L^1(D)$ for some $f \in \mathcal{B}(D,[0,1])$.
Then,  
\begin{enumerate}
\item[(i)] the laws $(u^{(N)},u^{\partial D,(N)})$ converges in distribution on 
$D([0,\infty),\mathcal{M}(D)) \times 
C([0,\infty),\mathcal{M}(\partial D))$ to
the limit $Q^{D,\beta,\gamma}_{\mu,f}$; 
\item[(ii)] for any $T<\infty$, the law of the triple
$(u^{(N)}_{[0,T]}, u^{\partial D,(N)}_T, v^{(N)}_T)$ on
$\mathcal{M}(D) \times \mathcal{M}(\partial D) \times 
L^1(D)$ converges to the law of $(U_{[0,T]},U^{\partial D}_T,f \exp(-U(0,T)))$ 
under $Q^{D,\beta,\gamma}_{\mu,f}$.
\end{enumerate}
\end{Theorem}

\textbf{Completion of the proof.}
The convergence of the exit measures can be deduced from the extended
martingale problem for test functions $\phi \in C^2_b(\overline{D})$. Indeed, return to 
the sub-subsequence studied before the statement of the theorem. 
Since $(u,v)$ solves (\ref{spde}) there is a continuous exit measure process
$u^{\partial D}$ solving the extended martingale problem. 
Using an approximation argument by $\phi_n \in C^2_0(\overline{D})$,
it follows that $m^{(N)}_t(\phi) \to m_t(\phi)$ in probability for all 
$\phi \in C^2_b(\overline{D})$.
Choose $\phi \in C^2_b(\overline{D})$
so that $h = \phi |_{\partial D}$ is non-negative on $\partial D$. 
Then passing to the limit in the extended 
approximate martingale 
problem (\ref{AMGPD1}) for $u^{(N)}$, we have shown convergence of all but one of the terms,
so that this last term $u^{\partial D, (N)}_t(h)$ must also converge in probability
to $u_t^{\partial D}(h)$.
The fact that $t \to u^{(N), \partial D}_{t}(h)$ and 
$t \to u^{\partial D}_{t}(h)$ are non-decreasing and continuous imply, 
at least along a further subsequence and for all $t$, that
$\sup_{s \leq t} |u^{(N), \partial D}_{s}(h) - u^{\partial D}_{s}(h)| \to 0$
in probability. A further subsequence argument allows us 
to conclude the same for a countable dense set 
in $C^0_b(\partial D)$, and this implies that 
$u^{\partial D,(N)} \to u^{\partial D}$ in $C([0,\infty), \mathcal{M}(\partial D))$.
The uniqueness in law of $(u, u^{\partial D})$ for solutions to (\ref{spde}) implies 
the convergence of $(u^{(N)},u^{\partial D,(N)})$ in part (i). 

Based on (\ref{approxvrep}), the proof of the 
$L^1$ convergence of part (ii) uses similar, but slightly simpler, tricks 
as above when dealing with the convergence of the key term 
$\int^t_0 u^{(N)}_s(v^{(N)}_s \phi) ds$, and the details are omitted. \qed
\section{Proof of the decomposition results} \label{s5}
This section contains the proofs of the decomposition results stated in section \ref{s2.2}.
In section \ref{s5.1} we show that the exit and occupation measures 
of the approximate discrete nutrient process, as constructed in section \ref{s4.1},
can be done in two stages. In section \ref{s5.2}, five different examples of 
this two-stage construction lead to five decomposition results for the approximation processes.
Passage to the limit is done in section \ref{s5.3}.
\subsection{A two-stage construction of the approximations} \label{s5.1}
The construction of the approximation process given in section \ref{s4.1}
is a pathwise construction. The approximation is constructed as a 
deterministic procedure, which we call the {\it basic 
construction}, applied to fixed realizations of the non-interacting DW($D,\gamma$) processes
$(u_k:k=0,1,\ldots)$ and exponential variables $(e_k:k=1,2,\ldots)$. 
In this section we show that the occupation and exit measures from this construction 
can be built up in two stages, where each stage applies the basic construction 
to a certain set of variables. This is described in a somewhat abstract
manner, and the reader might 
want to look ahead at one of the five examples in section \ref{s5.2} to be convinced 
it is a natural idea. For the rest of this subsection 
we act pathwise, supposing a single realization has been fixed of the underlying variables.

Integrating over $t$ in (\ref{approximationdefn}), or letting $t \to \infty$, we obtain
\begin{equation} \label{Sinfty}
u_{[0,\infty)} = u_{0,[0,\infty)} + \sum_{k \in \mathcal{S}_{\infty}} u_{k,[0,\infty)}
\quad \mbox{and} \quad 
u^{\partial D}_{\infty} = u^{\partial D}_{0,\infty} 
+ \sum_{k \in \mathcal{S}_{\infty}} u^{\partial D}_{k,\infty},
\end{equation}
where $\mathcal{S}_{\infty} = \lim_{t \to \infty} \mathcal{S}_t$ is the set of
labels of nutrient packages that are ever triggered. So the total occupation and exit measures 
are determined by the set $\mathcal{S}_{\infty}$
and the total occupation and exit measures of the non-interacting DW processes.
Moreover for $k \in \mathcal{S}_{\infty}$ we have
\[
u_{0,[0,\infty)}(\hat{\psi}_k) + \sum_{j \in \mathcal{S}_\infty} 
u_{j,[0,\infty)}(\hat{\psi}_k) > e_k
\]
while for the $k \not \in \mathcal{S}_{\infty}$ the converse inequality holds.
This is exactly the condition that $\mathcal{S}_{\infty}$ is a fixed point of the 
mapping $T$ defined in the abstract lemma below, once the choices $A=\{1,\ldots,K\}$ and
\[
e(k) = e_k, \quad f(k) = u_{0,[0,\infty)}(\hat{\psi}_k), \quad
M(j,k) = u_{j,[0,\infty)}(\hat{\psi}_k) \quad \mbox{for $j,k \in A$}
\]
have been made. 
\begin{Lemma} \label{abstractlemma1}
Suppose $A$ is a finite set, and fix
$e,f:A \to [0,\infty)$ and $M:A \times A \to [0,\infty)$.
Define, for $ B \subseteq A$,
\begin{equation} \label{Tmapping}
T(B) = \left\{a \in A: f(a) + \sum_{b \in B} M(b,a) > e(a) \right\}.
\end{equation}
Then there is a unique smallest fixed point $S \subseteq A$ of $T$, that is 
$T(S)=S$ and $S$ is contained in any other fixed point. Moreover if 
$B \subseteq S$ then $T^n(B)$ equals $S$ for large $n$. 
\end{Lemma}

\noindent
\textbf{Proof.}
Since $M(a,b) \geq 0$ we see that if 
$B \subseteq B'$ then $T(B) \subseteq T(B')$. So
$T^n(\emptyset)$ increases to a limit $S$ which must be a fixed point of 
$T$. Any other fixed point $S'$ contains $T^n(\emptyset)$ for all $n$ and 
hence contains $S$ and
the uniqueness of the smallest fixed point is clear. If $B \subseteq S$
then $T^n(\emptyset) \subseteq T^n(B) \subseteq T^n(S) = S$ and the result follows
since $T^n(\emptyset)=S$ for large $n$. \qed

It is straightforward to check that $\mathcal{S}_{\infty}$ is the 
smallest fixed point of $T$. Indeed, let $k_m$ be the label of the 
$m$th nutrient package to be triggered, and let $\tau_m$ be the time it is triggered;
then the definition of $\tau_{m}$ shows, when $\tau_{m}<\infty$,
that $k_{m} \in T(\{k_1, \ldots,k_{m-1}\})$. Inductively,   
$\{k_1,\ldots,k_m\} \subseteq T^m(\emptyset)$ must hold and letting $m \to \infty$ 
shows that $\mathcal{S}_{\infty}$ is contained in the smallest fixed point of $T$.

\textbf{Example 0.} 
To illustrate the use of this lemma, we give here the analogue of the comparison
Lemma \ref{comparison3} for the approximation processes. Suppose that $\sum_{k=1}^L \psi_k = f$ and
$\sum_{k=L+1}^K \psi_k = g$. Then running the basic construction on the processes
$(u_k:k=0,1,\ldots,K)$ and exponential variables $(e_k:k=1,\ldots,K)$
we obtain the approximation process $(u,u^{\partial D},v)$ started at
$\mu$ and with initial nutrient level $f+g$, and the triggered set $\mathcal{S}_{\infty}$ that is the 
smallest fixed point of a mapping $T$. We now set
$ \tilde{u}_0 = u_0 + \sum_{k=L+1}^K u_k$.
Now running the basic construction on the realizations of $(\tilde{u}_0, u_1,u_2,\ldots,u_L)$ 
and the exponential variables $(e_k:k=1,\ldots,L)$ we obtain
an approximation process $(\tilde{u},\tilde{u}^{\partial D}, \tilde{v})$ started at
$\mu + \beta g \, dx$ and with initial nutrient level $f$, and the triggered set 
$\tilde{\mathcal{S}}_{\infty}$ that is the 
smallest fixed point of a corresponding mapping $\tilde{T}$. Moreover, for $A \subseteq \{1,\ldots,L\}$,
\begin{eqnarray*}
\tilde{T}(A) & = & \left\{k \in \{1,\ldots,L\}:
\tilde{u}_{0,[0,\infty)}(\hat{\psi}_k) + \sum_{j=1}^L \I(j \in A) 
u_{j,[0,\infty)}(\hat{\psi}_k) > e_k \right\} \\
&=& \left\{k \in \{1,\ldots,L\}:
u_{0,[0,\infty)}(\hat{\psi}_k) + \sum_{j=1}^K \I(j \in A \cup \{L+1,\ldots,K\}) 
u_{j,[0,\infty)}(\hat{\psi}_k) > e_k \right\} \\
&=& T(A \cup \{L+1,\ldots,K\}) \setminus \{L+1,\ldots,K\}.
\end{eqnarray*}
Choosing $A=\emptyset$ and then iterating we find that
$T^n(\emptyset) \setminus \{L+1,\ldots,K\} \subseteq \tilde{T}^n(\emptyset)$
and hence, by Lemma \ref{abstractlemma1}, that $\mathcal{S}_{\infty} \setminus \{L+1,\ldots,K\} \subseteq \tilde{\mathcal{S}}_{\infty}$.
Combining this with (\ref{Sinfty}) shows that $(u_{[0,\infty)},u^{\partial D}_{\infty})$
is smaller than $(\tilde{u}_{[0,\infty)},\tilde{u}^{\partial D}_{\infty})$.

We now describe the two-stage procedure, which will build up the 
set $\mathcal{S}_{\infty}$ of triggered nutrient packages in two steps.
We suppose there is a splitting of each of the DW processes $(u_k:k=0,\ldots,K)$, and their exit measures, 
into two parts, $u^{-}_k$ and $u^{+}_k$, each with continuous 
measure valued paths, and satisfying
\begin{equation} \label{splitting}
\left(u_{k,[0,\infty)}, u^{\partial D}_{k,\infty} \right) = 
\left(u^{-}_{k,[0,\infty)} + u^{+}_{k,[0,\infty)}, 
u^{\partial D,-}_{k,\infty} +
u^{\partial D,+}_{k,\infty} \right).
\end{equation}
In stage one we apply the basic construction using the
processes $(u^{-}_k:k =0,1,\ldots,K)$ and the exponential variables 
$(e_k:k =1,\ldots,K)$ to create a process 
$(u^-,v^-)$ and a set $\mathcal{S}^-_{\infty}$ of triggered nutrient packages.
In particular
\begin{equation} \label{Sinfty-}
u^-_{[0,\infty)}  =   u^{-}_{0,[0,\infty)} + 
\sum_{k \in \mathcal{S}^-_{\infty}} u^{-}_{k,[0,\infty)}.
\end{equation}
We then define
$ \hat{u}_{0,t} = u^{+}_{0,t} + \sum_{k \in \mathcal{S}^-_{\infty}} 
u^{+}_{k,t}$ and, for $k \in \{1,\ldots,K\} \setminus \mathcal{S}^-_{\infty}$,
\begin{equation} \label{hatdefn}
\hat{u}_{k,t}  =  u_{k,t} \quad \mbox{and} \quad 
\hat{e}_k  =  e_k - u^{-}_{0,[0,\infty)}(\hat{\psi}_k) 
- \sum_{j \in \mathcal{S}^-_{\infty}} u^{-}_{j,[0,\infty)}(\hat{\psi}_k). 
\end{equation}
In stage two we run the basic construction on the processes $\hat{u}_{0}$ and
$(\hat{u}_k:k \in \{1,\ldots,K\} \setminus \mathcal{S}^-_{\infty})$,
with the nutrient packages triggered using the values
$(\hat{e}_k:k \in \{1,\ldots,K\} \setminus \mathcal{S}^-_{\infty})$.
This leads to a second process $(u^+,v^+)$ and a second set of triggered
nutrient packages $\mathcal{S}^+_{\infty} \subseteq
\{1,\ldots,K\} \setminus \mathcal{S}^-_{\infty}$ satisfying
\begin{equation} \label{Sinfty+}
u^+_{[0,\infty)}  =   u^{+}_{0,[0,\infty)} + 
\sum_{k \in \mathcal{S}^-_{\infty}} u^{+}_{k,[0,\infty)} +
\sum_{k \in \mathcal{S}^+_{\infty}} u_{k,[0,\infty)}.
\end{equation}
Adding (\ref{Sinfty-}) and (\ref{Sinfty+}) and comparing to (\ref{Sinfty}), we see that
$ u_{[0,\infty)} = u^-_{[0,\infty)} + u^+_{[0,\infty)}$
will hold provided that the equality $\mathcal{S}_{\infty} = \mathcal{S}^-_{\infty}
\cup \mathcal{S}^+_{\infty}$ holds. A similar argument shows that this equality 
is also sufficient to ensure that the total exit measures satisfy 
$u^{\partial D}_{\infty} = u^{\partial D,-}_{\infty} + u^{\partial D,+}_{\infty}$.

To verify this equality we use the following second abstract lemma.
With the choices
\[
e(k) = e_k, \quad f^{\pm}(k) = u^{\pm}_{0,[0,\infty)}(\hat{\psi}_k), \quad
M^{\pm}(j,k) = u^{\pm}_{j,[0,\infty)}(\hat{\psi}_k) \quad \mbox{for $j,k \in A$}
\]
note that the definition of $\hat{f},\hat{e},\hat{M}$ in the lemma 
corresponds to the definition of $\hat{u}_0, \hat{e}_k$ in the stage two construction.
The conclusion of the lemma then implies that 
the equality  $\mathcal{S}_{\infty} = \mathcal{S}^-_{\infty}
\cup \mathcal{S}^+_{\infty}$ does indeed hold. 

For the statement of the lemma we 
denote the unique smallest fixed point in Lemma 
\ref{abstractlemma1} by $S=S(A,e,f,M)$ to indicate its dependence.
\begin{Lemma} \label{abstractlemma2}
Suppose $A$ is a finite set, and fix $e,f^-,f^+:A \to [0,\infty)$ and $M^-,M^+:A \times A \to [0,\infty)$.
Set $M=M^- + M^+$ and $f=f^- + f^+$. Then the fixed point $S(A,e,f,M)$ can be broken into two subsets
as follows. First let $S^-=S(A,e,f^-,M^-)$. Then define, for $a \in A \setminus S^-$,
\[
\hat{f}(a)=f^+(a) + \sum_{b \in S^-} M^+(b,a), \quad 
\hat{e}(a)=e(a) - f^-(a)-\sum_{b \in S^-} M^-(b,a) 
\]
and write $\hat{M}$ for the restriction of $M$ to $(A \setminus S^-) \times (A \setminus S^-)$.
Then $\hat{e} \geq 0$ and 
if $S^+$ is the fixed point $S(A \setminus S^-,\hat{e}, \hat{f}, \hat{M})$ we have
\[
S(A,e,f,M) = S^- \cup S^+.
\] 
\end{Lemma}

\noindent
\textbf{Proof.}
Let $T$ (respectively $T^-$ and $T^+$) {be the map defined by (\ref{Tmapping}) used 
for the definition of $S$ (respectively $S^-$ and $S^+$).
The fact that $\hat{e}(a) \geq 0$ for $a \in A \setminus S^-$ follows from the fact that
$S^-$ is a fixed point of $T^-$. 

Since $f \geq f^-$ and $M \geq M^-$ we see that $T^-(B) \subseteq T(B)$.
(Note in particular $S^- \subseteq T^n(S^-)$ for all $n \geq 1$.)
Then $(T^-)^n(\emptyset) \subseteq T^n (\emptyset)$ for all $n$ and hence $S^- \subseteq S$. 

For $ B \subseteq A \setminus S^-$ we have
\begin{eqnarray*}
T^+(B) & = & 
 \left\{a \in A \setminus S^-: \hat{f}(a) + \sum_{b \in B} \hat{M}(b,a) > \hat{e}(a) \right\} \\
& = &
 \left\{a \in A \setminus S^-: f^+(a) + \sum_{b \in S^-} M^+(b,a) +
\sum_{b \in B} M(b,a) > e(a) - f^-(a) - \sum_{b \in S^-} M^-(b,a) \right\} \\
& = &
 \left\{a \in A \setminus S^-: f(a) + \sum_{b \in B \cup S^-} M(b,a)
> e(a)\right\} \\
& = & T(B \cup S^-) \setminus S^-.
\end{eqnarray*}
Applying this first to $B=\emptyset$ and iterating, we obtain 
$(T^+)^n(\emptyset) = T^n(S^-)\setminus S^-$. The equality here uses
$S^- \subseteq T^n(S^-)$. Using the final statement of 
Lemma \ref{abstractlemma1},
we may let $n \to \infty$ to obtain $S^+ = S \setminus S^-$ as desired. \qed 
\subsection{Comparison theorems for the approximations} \label{s5.2}
\noindent \textbf{Notation} We denote the law of the approximation $(u_t,u^{\partial D}_t:t \geq 0)$
constructed in section \ref{s4.1}, on the space 
$D([0,\infty),\mathcal{M}(D)) \times C([0,\infty),\mathcal{M}(\partial D))$,
by $Q^{N,D,\beta,\gamma}_{\mu,f}$. 

We now give five examples of the two stage construction of the last section, each 
leading to a decomposition theorem for the discrete nutrient approximation process. 
In each example we describe the processes $u^-_k,u^+_k$ used for the 
splitting (\ref{splitting}). We let 
\[
\mathcal{G}^- = \sigma\{u^-_{0}\} \vee \sigma \{\mathcal{S}^-_{\infty} \}
 \vee \sigma\{u^-_{k} \I(k \in \mathcal{S}^-_{\infty}) : 1 \leq k \leq K\}
\]
be the information gained by observing the stage one process. 

\vspace{.1in}

\noindent
\textbf{Example 1. Decomposition for the initial condition $\mu$.}
Suppose that $\mu=\mu^- + \mu^+$. Take independent DW($D,\gamma$) 
processes $(u^-_{0,t})$ and $(u^+_{0,t})$, with initial conditions $\mu^-$ and $\mu^+$. Then
$u_{0,t}= u^-_{0,t} + u^+_{0,t}$ is still a DW($D,\gamma$) process. Using
independent DW($D,\gamma$) processes $(u_k:1 \leq k \leq K)$ started at $u_{k,0}=\psi_k$,
we set
\[ 
u^-_{k,t} = u_{k,t} \quad \mbox{and} \quad u^+_{k,t} =0 \quad \mbox{for $1 \leq k \leq K$ and $t \geq 0$}
\]
so that the splitting (\ref{splitting}) holds trivially. 
The first stage leads to a process $(u^-,u^{\partial D,-})$ with law
$Q^{N,D,\beta,\gamma}_{\mu^-,f}$. 
The memorylessness property of exponentials implies that, conditional on $\mathcal{G}^-$,
the variables
$(\hat{e}_k: k \in \{1,\ldots,K\} \setminus \mathcal{S}^-_{\infty})$ defined in (\ref{hatdefn}) 
remain independent rate one exponential variables.
The processes $\hat{u}_k$ for $k \in \{0,1,\ldots,K\} \setminus \mathcal{S}^-_{\infty}$
are, conditionally on $\mathcal{G}^-$, independent DW($D,\gamma$) processes.
Note that $v^-_{\infty} = \sum_{k \not \in \mathcal{S}^-_{\infty}} \psi_k$
is $\mathcal{G}^-$ measurable.  
The second stage therefore produces a process $(u^+,u^{\partial D,+})$ which, conditional on 
$\mathcal{G}^-$, has law $Q^{N,D,\beta,\gamma}_{\mu^+,v^-_{\infty}}$.

\vspace{.1in}

\noindent
\textbf{Example 2. Decomposition for the initial condition $f$.}
We split the set of nutrient package labels into two parts 
$\{1,\ldots,L\}$ and $\{L+1,\ldots,K\}$ and set $f^- = \sum_{k=1}^L \psi_k, \,
f^+ = \sum_{k=L+1}^K \psi_k$.
Then we form the splitting (\ref{splitting}) out of the 
independent DW($D,\gamma$) processes $(u_k:0 \leq k \leq K)$ started
at $\mu,\psi_1,\ldots,\psi_K$ as follows: for $t \geq 0$
\[
u^-_{k,t} = \left\{ \begin{array}{ll} 
u_{k,t} & \mbox{for $0 \leq k \leq L$,} \\ 
0 & \mbox{for $L+1 \leq k \leq K$,}
\end{array} \right. \quad \mbox{and} \quad
u^+_{k,t} = \left\{ \begin{array}{ll} 
0 & \mbox{for $0 \leq k \leq L$,} \\ 
u_{k,t} & \mbox{for $L+1 \leq k \leq K$.}
\end{array} \right.
\]
The first stage of the construction produces a process
$(u^-,u^{\partial D,-})$ with law $Q^{N,D,\beta,\gamma}_{\mu,f^-}$. 
Note that $u^+_0 = \hat{u}_{0,0} = \beta \sum_{k=L+1}^K 
\psi_k I(k \in \mathcal{S}_{\infty}^-)$ and 
$v^+_0 = \sum_{k \not \in \mathcal{S}^-_{\infty}} \psi_k$ are measurable with respect to 
the sigma field $\mathcal{G}^-$. Moreover, as above, the variables
$(\hat{e}_k, \hat{u}_k: k \in \{1,\ldots,K\} \setminus \mathcal{S}^-_{\infty})$ 
are, conditionally upon $ \mathcal{G}^-$, independent exponential variables
and DW($D,\gamma$) processes. 
The second stage therefore produces a process $(u^+,u^{\partial D,+})$ which, 
conditional on $\mathcal{G}^-$, has law $Q^{N,D,\beta,\gamma}_{u^+_0,v^+_0}$.

\vspace{.1in}

\noindent
\textbf{Example 3. Spatial Markov property.} 
We fix domains $D^- \subseteq D^+$. We will choose 
the finite partition $D^+= \cup D_j$ used in the construction of the approximation so that
each $D_j$ is either a subset of $D^-$ or a subset of $D^+ \setminus D^-$.
Then a subset, which we list as $(\psi_k: 1 \leq k \leq L)$, will 
satisfy $\sum_{k=1}^L \psi_k = f I(D^-)$ and 
$\sum_{k=L+1}^K = f I(D^+ \setminus D^-)$.
We apply the spatial Markov Lemma \ref{smlemma} to find
processes $(u_k,u^-_k,u^+_k: 0 \leq k \leq K)$ so that 
$(u_k:0 \leq k \leq K)$ are independent DW($D^+,\gamma$) processes 
with initial conditions $\mu, \psi_1, \ldots, \psi_K$; 
$(u^-_k:0 \leq k \leq K)$ are independent DW($D^-,\gamma$) processes
with initial conditions $\mu \I(D^-), \psi_1 , \ldots, \psi_L,0,\ldots,0$;
and, conditional on $\sigma\{u^-_k: 0 \leq k \leq  L\}$, 
$(u^+_k:0 \leq k \leq K)$ are independent DW($D^+,\gamma$) processes
with initial conditions 
$ u^+_{0,0}=\mu \I(D^+ \setminus D^-) + u^{\partial D^-,-}_{0,[0,\infty)}|_{D^+}$,
$ u^+_{k,0}= u^{\partial D^-,-}_{k,\infty}|_{D^+}$ for $1 \leq k \leq L$ and 
$u^+_{k,0}= \psi_k$ for $L+1 \leq k \leq K$.
In addition, the Lemma \ref{smlemma} ensures the splitting
\[
\left(u_{k,[0,\infty)}, u^{\partial D^+}_{k,\infty} \right) = 
\left(u^{-}_{k,[0,\infty)} + u^{+}_{k,[0,\infty)}, u^{\partial D^-,-}_{\infty} |_{\partial D^+} +
u^{\partial D^+,+}_{k,\infty} \right)
\quad \mbox{for $0 \leq k \leq K$.}
\]
This gives a suitable splitting as in (\ref{splitting}) for the two part construction, where 
the exit measures are all evaluated as measures on $\partial D^+$. 
The first stage of the construction produces a process
$(u^-,u^{\partial D^-,-})$ with law $Q^{N,D^-,\beta,\gamma}_{\mu \I(D^-),f I(D^-)}$. 
The second stage produces a process $(u^+,u^{\partial D^+,+})$ which, 
conditional on $\mathcal{G}^-$, has law
$Q^{N,D^+,\beta,\gamma}_{u^+_0,v^+_0}$, where
$u^+_0 = \mu \I(D^+ \setminus D^-) + u^{\partial D^-,-}_{[0,\infty)}|_{D+}$ and
$v^+_0 = f \I(D^+ \setminus D^-) + v^-_{\infty} \I(D^-)$. 

\vspace{.1in}

\noindent
\textbf{Example 4. Decomposition for the birth rate $\beta$.}
Fix $\beta=\beta^- + \beta^+$. Take independent
DW($D,\gamma$) processes $(u^-_k,u^+_k: 0 \leq k \leq K)$
with initial conditions $u^-_{0,0}=\mu$, $u^+_{0,0}=0$,
$u^-_{k,0}= \beta^- \psi_k$ and $u^+_{k,0}= \beta^+ \psi_k$.
Set $u_{k}= u^-_{k}+u^+_{k}$ so that 
$(u_k: 0 \leq k \leq K)$ satisfy the splitting (\ref{splitting})
and are themselves DW($D,\gamma$) processes.
Then the first stage of the construction produces a process
$(u^-,u^{\partial D,-})$ with law $Q^{N,D,\beta^-,\gamma}_{\mu,f}$. 
The second stage produces a process $(u^+,u^{\partial D,+})$ which, 
conditional on $\mathcal{G}^-$, has law
$Q^{N,D,\beta,\gamma}_{u^+_0,v^+_0}$, where
$u^+_0 = \beta^+ (f-v^-_{\infty}) \, dx$ and
$v^+_0 = v^-_{\infty}$. 

\vspace{.1in}

\noindent
\textbf{Example 5. Decomposition for the death rate $\gamma$.}
Fix $\gamma^- = \gamma + \gamma^+$. 
Lemma \ref{gammalem} below ensures we can find 
processes $(u_k,u^-_k,u^+_k:0 \leq k \leq K)$ satisfying the splitting 
(\ref{splitting}) and so that 
$(u_k:0\leq k \leq K)$ are independent DW($D,\gamma$) processes
with initial conditions $\mu, \psi_1, \ldots, \psi_K$; 
$(u^-_k:0 \leq k \leq K)$ are independent DW($D,\gamma^-$) processes 
with initial conditions $\mu, \psi_1, \ldots, \psi_K$;
and, conditional on $\sigma\{u^-_k:0 \leq k \leq K\}$, 
$(u^+_k:0 \leq k \leq K)$ are independent DW($D,\gamma$) processes 
with initial conditions $(\gamma^+ u^-_{k,[0,\infty)}:0 \leq k \leq K)$.
Then the first stage of the construction produces a process
$(u^-,u^{\partial D,-})$ with law $Q^{N,D,\beta,\gamma^-}_{\mu,f}$. 
The second stage produces a process $(u^+,u^{\partial D,+})$ which, 
conditional on $\mathcal{G}^-$, has law
$Q^{N,D,\beta,\gamma}_{u^+_0,v^+_0}$, where
$u^+_0 = \gamma^+ u^-_{[0,\infty)}$ and
$v^+_0 = v^-_{\infty}$. 

\begin{Lemma} \label{gammalem}
Suppose  $\gamma^- = \gamma + \gamma^+$. 
There exists a coupling of three processes: 
$u$ a DW($\gamma,D$)
process with initial condition $\mu$; $u^-$ a DW($\gamma^-,D$)
process with initial condition $\mu$; and $u^+$ which, conditional
on $\sigma\{u^-\}$, is a DW($\gamma,D$) process with initial condition 
$\gamma^+ u^-_{[0,\infty)}$; and moreover these 
processes satisfy the splitting
\[
\left(u_{[0,\infty)}, u^{\partial D}_{\infty} \right) = 
\left(u^{-}_{[0,\infty)} + u^{+}_{[0,\infty)}, 
u^{\partial D,-}_{\infty} +u^{\partial D,+}_{\infty} \right).
\]
\end{Lemma}

\noindent
\textbf{Proof.}
Construct a coupling of two processes $(u^-, \tilde{u})$
as follows. Let $u^-$ be a DW($D,\gamma^-$) process with initial condition $\mu$.
Conditionally on $\sigma\{u^-\}$, let $\tilde{u}$ be a  
DW process on the space $D \times \R$, with spatial motion that is Brownian
on the first component $D$ and zero on the second component $\R$, 
with annihilation rate $\gamma$ and with initial condition $\gamma^+ u_r(dx) I(r \geq 0) dr$.   
If we define
\[
u^+_t(dx) = \tilde{u}_t(dx \times \R) \quad \mbox{and} \quad u^{\partial D,+}_t(A)
= u^{\partial(D \times \R)}_t(A \times \R) \; \mbox{for $ A \subseteq \partial D$}
\]
then conditionally on $\sigma\{u^-\}$ the process $u^+$ is a DW($D,\gamma$) process started
at $\gamma^+ u^-_{[0,\infty)}$ and with exit measure $u^{\partial D,+}$. 
Define measures $I$ on $[0,\infty) \times D$ and $I^{\partial D}$ on $[0,\infty) \times \partial D$
by 
\begin{eqnarray*}
I(h^{(2)}) & = & I^-(h^{(2)}) + I^+(h^{(2)}) = \int^{\infty}_0 u^-_s(h^{(2)}_s) ds + \int^{\infty}_0 \int_{D \times \R} h^{(2)}_{s+r}(x)
\tilde{u}_s(dx,dr) \, ds, \\
I^{\partial D}(h^{(3)}) & = & I^{\partial D,-}(h^{(3)}) + I^{\partial D,+}(h^{(3)})
= \int^{\infty}_0 d_s u^{\partial D}_s(h^{(3)}_s) + \int^{\infty}_0 \int_{\partial(D \times \R)} 
h^{(3)}_{s+r}(x) d_s \tilde{u}_s(dx,dr)
\end{eqnarray*}
for bounded measurable $h^{(2)}:[0,\infty) \times D \to \R$ and $h^{(3)}:[0,\infty) \times \partial D \to \R$.
We claim that $(I(dt,dx), I^{\partial D}(dt,dx))$ has the same law as
$(U_t(dx)dt, d_t U_t^{\partial D}(dx))$ under $Q_{\mu}^{D,\gamma}$. This allows us to define the the required process
$(u,u^{\partial D})$ using $(I, I^{\partial D})$ and the required splitting follows from the definitions of
$I$ and $u^+$. 

To establish the claim on the law of $(I,I^{\partial D})$ we fix smooth 
non-negative $h^{(2)}_s(x),h^{(3)}_s(x)$ that 
vanish for $s \geq t$. It is enough to check that 
$ E[ \exp(-I(h^{(2)}) - I^{\partial D}(h^{(3)})) ] = \exp(-\mu(\phi_t)) $
where $(\phi_s: s \in [0,t])$ is the solution to the log-Laplace equation (\ref{LFD1+}). We start by calculating the conditional
expectation $E[ \exp(-I(h^{(2)}) - I^{\partial D}(h^{(3)})) | \sigma\{u^-\}]$. Note the almost sure limits
\[
I^+(h^{(2)}) = \lim_{N \to \infty} \sum_{k=0}^{N-1} \int^{\infty}_0 \int_{D \times [\frac{k}{N}, \frac{k+1}{N})}
h^{(2)}_{s+\frac{k}{N}}(x) \tilde{u}_s(dx,dr) ds =
\lim_{N \to \infty} \sum_{k=0}^{N-1} \int^{\infty}_0 X_{K,s}(h^{(2)}_{s+\frac{k}{N}}) ds 
\]
and 
\[
I^{\partial D,+}(h^{(3)}) = \lim_{N \to \infty} \sum_{k=0}^{N-1} \int^{\infty}_0 \int_{\partial(D \times \R)}
\I_{r \in [\frac{k}{N}, \frac{k+1}{N})} h^{(3)}_{s+\frac{k}{N}}(x) d_s \tilde{u}_s(dx,dr) =
\lim_{N \to \infty} \sum_{k=0}^{N-1} \int^{\infty}_0 d_s X_{K,s}^{\partial D} (h^{(2)}_{s+\frac{k}{N}})
\]
where $(X_{k,t}:t \geq 0)$, for $k=0,1,\ldots$, are defined by $X_{k,t}(A) = \tilde{u}_t(A \times [\frac{k}{N}, \frac{k+1}{N}))$ and are,
conditionally on $\sigma\{u^-\}$, independent DW($D,\gamma$) processes with initial conditions
$X_{k,0} = \gamma^+ u^-_{[\frac{k}{N}, \frac{k+1}{N}]}$. Using the Laplace functional (\ref{LFD1}) of 
$(X_k,X^{\partial D}_k)$ we find
\begin{eqnarray*}
E[ e^{-I^+(h^{(2)}) - I^{\partial D,+}(h^{(3)})} | \sigma\{u^-\}] & = & 
\lim_{N \to \infty} \exp(-\gamma^+ \sum_{k=0}^{N-1} u^-_{[\frac{k}{N}, \frac{k+1}{N}]}(\phi_{t-\frac{k}{N}})) \\
& = & e^{-\gamma^+ \int^t_0 u^-_s(\phi_{t-s}) ds}.
\end{eqnarray*}
Then we complete the expectation as
\[
E[ e^{-I(h^{(2)}) - I^{\partial D}(h^{(3)})} ] = 
E[ e^{-\int^t_0 u^-_s(h^{(2)} + \gamma^+ \phi_{t-s}) ds - \int^t_0 d_s u^{\partial D,-}_s(h^{(3)}_s)} ] 
\]
which can be calculated again using the Laplace functional of $(u^-,u^{\partial D,-})$. It is straightforward 
to check the required log-Laplace function is solved by $\phi$ and this proves the claim. 
\qed

\noindent
\textbf{Remark.}
A more natural splitting for the final example 
is to let $(u^-_k,u^+_k)$ solve, for each $k$, the system
\begin{eqnarray*}
\partial_t u^- & = &  \Delta u^- -\gamma^- u^- + \sqrt{u^-} \, \dot{W}^-, \\ 
\partial_t u^+ & = &  \Delta u^+ + \gamma^+ u^- -\gamma u^+ + \sqrt{u^+} \, \dot{W}^+, 
\end{eqnarray*}
with $u^+_0=0$ and with orthogonal martingale terms. 
Then $u^-_k$ is a DW($D,\gamma-$) process and
$u_k=u_k^- +  u_k^+$ is a DW($D,\gamma$) process.
The first stage of the construction produces a process
$(u^-,v^-)$ with law $Q^{N,D,\beta,\gamma^-}_{\mu,f}$.
The second stage produces a process $(u^+,v^+)$ which, 
conditional on $\mathcal{G}^-$, has the law
of the approximation process with parameters $\beta,\gamma$, but with
an extra immigration term $\gamma^+ u^-$. This leads to a corresponding 
decomposition for solutions with immigration to (\ref{spde}), but we shall not 
make any use of it. 
\subsection{Completion of the proof of the comparison results} \label{s5.3}
This section contains the proofs of the decomposition lemmas \ref{comparison1}, \ref{comparison2}, 
\ref{comparison3} and \ref{sMp}, using 
the comparison results for the approximations 
from section \ref{s5.2} with the 
convergence of the approximations established in section \ref{s5.1}. 

We start with a simple estimate controlling the death time $\tau= \inf\{t: U_t = 0\}$ 
for the approximation processes
on a bounded domain uniformly in $N$. The change of measure argument
given at the beginning of section \ref{s3.3} can be used to show the limiting death time
has a tail $Q^{D,\beta,\gamma}_{\mu,f}[\tau \geq t] \leq C(\beta,p,|D|, \mu) t^{-p}$ for any
$p<1$. This argument does not apply directly to the approximations since they
are not absolutely continuous with respect to a DW process. However, the argument below 
follows the same main idea. It is combined with a first moment argument that leads to 
a suboptimal bound, but which is sufficient for our needs.
\begin{Lemma} \label{simpledeath}
There exists $C(\beta)<\infty$ so that for all $t \geq e$ and $N \geq 1$
\[
Q^{N,D,\beta,\gamma}_{\mu,f} \left[ \tau \geq t \right] \leq C(\beta) 
\left((\mu(1))^{1/2} + \langle f,1 \rangle \right) \ln(t)^{-1}. 
\]
\end{Lemma}

\textbf{Proof.}
Take an approximation process $u$ with law $Q^{N,D,\beta,\gamma}_{\mu,f}$.
Taking $\phi=1$ in the martingale problem (\ref{AMGPD1}) and 
$\lambda_s = -2/s$, so that $\dot{\lambda}_s = -\lambda_s^2/2$, we 
find for fixed $t>0$ and $s \in [0,t)$
\begin{eqnarray*}
d e^{-\lambda_{t-s} u_s(1)} &=&
e^{-\lambda_{t-s} u_s(1)} \left( u_s(1) (\dot{\lambda}_{t-s} + \frac12 \lambda_{t-s}^2 - \gamma \lambda_{t-s}) ds
+ \lambda_{t-s} du^{\partial D}_s(1) - \beta \lambda_{t-s} u_s(v_s) ds  \right) \\
&& + \sum_{s \leq t} D_s(e^{-\lambda_{t-s} u_s(1)}) + \lambda_{t-s} e^{-\lambda_{t-s} u_{s-}(1)} D_s u_s(1)  
+ \mbox{martingale increments} \\
& \geq &  - \beta \lambda_{t-s} e^{-\lambda_{t-s} u_s(1)}
u_s(v_s) ds + \mbox{martingale increments}.
\end{eqnarray*}
Let $m(dx,ds)$ be the continuous martingale measure constructed from the martingales
$m_t(\phi)$ in (\ref{AMGPD1}) and, as in the change of measure
arguments, set $M_t=\int^t_0 \int v_s(x) m(dx,ds)$. Then Ito's lemma, 
using the cross variation,
\begin{eqnarray*}
d \left[ \mathcal{E}_s(- \beta M), e^{-\lambda_{t-s} u_s(1)} \right]_s
& = &  \beta \lambda_{t-s} e^{-\lambda_{t-s} u_s(1)} \mathcal{E}_s(- \beta M) 
d \left[ M_s, m_s(1)\right]_s  \\
& = & \beta \lambda_{t-s} e^{-\lambda_{t-s} u_s(1)} \mathcal{E}_s(- \beta M) u_s(v_s) ds
\end{eqnarray*}
shows that $ s \to \mathcal{E}_s(- \beta M) (1- \exp(-\lambda_{t-s} u_s(1))$ is a non-negative
supermartingale on $s \in [0,t)$. Taking expectations at $s_n \uparrow t$, we find
\begin{eqnarray*}
E \left[ \mathcal{E}_t(- \beta M) \I(\tau >t) \right] 
& \leq & \lim_{n \to \infty} 
E \left[ \mathcal{E}_{s_n} (- \beta M) (1-e^{-\lambda_{t-s_n} u_{s_n}(1)}) \right] \\ 
& \leq &  \left(1- e^{-\lambda_t \mu(1)} \right)
\leq 2\mu(1)/t.
\end{eqnarray*}
Let $\sigma_K = \inf\{s: \int^s_0 u_s(v_s) ds \geq K\}$. Then, using Cauchy-Schwarz,  
and $v_s \leq f \leq 1$,
\begin{eqnarray}
\left( P \left[ \tau > t, \sigma_K > t \right] \right)^2 & \leq &
E \left[ \I( \tau > t) \mathcal{E}_t(-\beta M) \right] 
E \left[ \I( \sigma_K > t) \mathcal{E}^{-1}_t(-\beta M) \right] \nonumber \\
& \leq & \frac{2\mu(1)}{t}
E \left[ \I( \sigma_K > t) \mathcal{E}_t(+\beta M)  
e^{ \beta^2\int^t_0 u_s(v^2_s) ds} \right]  \leq  \frac{2\mu(1)}{t}
e^{ \beta^2 K}. \label{death101}
\end{eqnarray}
To control $\sigma_K$ we use a crude first moment bound. 
From (\ref{aveqn}) we obtain
\[
E[ (v_t,1)]   =  \langle f,1 \rangle - \int^t_0 \int_D E [ v_s(x) u_s(\hat{\psi}_x) ] dx ds
=  \langle f,1\rangle  - \int^t_0 E [ u_s(v_s) ] ds
\]
(using the fact that $v$ is constant on the partition sets $D_j$).
Then $P[\sigma_K \leq t] \leq \langle f,1 \rangle /K$ by Markov's inequality. Combining this with
(\ref{death101}) and the choice $K= c \ln t$, for small $c= c(\beta)$, leads to the desired bound.
\qed

We may characterize the total occupation and exit measures via the Laplace
functional, defined by 
\[
\Phi_{\phi,\psi}(\nu_1,\nu_2) = e^{-\nu_1(\phi) - \nu_2(\psi)}
\]
for $\nu_1 \in \mathcal{M}(D)$, $\nu_2 \in \mathcal{M}(\partial D)$, 
$0 \leq \phi \in C^0_b(D)$ and $0 \leq \psi \in C^0_b(\partial D)$. 
Fix $\mu \in \mathcal{M}(D)$ and $f \in \mathcal{L}^1(D)$, for a bounded domain $D$.
Choose approximations processes with initial nutrient levels $1 \geq f^{(N)} \to f$ in $L^1$ and
$\mu^{(N)} \to \mu$ weakly. Theorem \ref{convergence} (ii) implies that 
\begin{eqnarray}
\lim_{N \to \infty} Q^{N,D,\beta,\gamma}_{\mu^{(N)},f^{(N)}}
\left[ \Phi_{\phi,\psi} (U_{[0,\infty)},U^{\partial D}_{\infty}) \right] 
&=&  \lim_{N \to \infty} \lim_{T \to \infty} Q^{N,D,\beta,\gamma}_{\mu^{(N)},f^{(N)}} 
\left[ \Phi_{\phi,\psi} (U_{[0,T]},U^{\partial D}_{T}) \right] \nonumber \\
&=& \lim_{T \to \infty} \lim_{N \to \infty} Q^{N,D,\beta,\gamma}_{\mu^{(N)},f^{(N)}}
\left[ \Phi_{\phi,\psi} (U_{[0,T]},U^{\partial D}_{T}) \right] \nonumber \\
& = & \lim_{T \to \infty} Q^{D,\beta,\gamma}_{\mu,f} 
\left[ \Phi_{\phi,\psi} (U_{[0,T]},U^{\partial D}_{T}) \right] \nonumber \\
& = & Q^{D,\beta,\gamma}_{\mu,f} 
\left[ \Phi_{\phi,\psi} (U_{[0,\infty)},U^{\partial D}_{\infty}) \right]. \label{conv12}
\end{eqnarray}
The interchange of limits is justified by the above lemma, which gives uniform (in $N$) control
on $Q^{N,D,\beta,\gamma}_{\mu^{(N)},f^{(N)}}[U_T=0]$, 
and the fact that $U_{[0,\infty)} = U_{[0,T]}$ and $U^{\partial D}_{\infty} = U^{\partial D}_{T}$
on the set $\{U_t =0\}$. 

\textbf{Proof of Lemma \ref{comparison3}.}
Returning to Remark $0$ in section \ref{s5.1} we see that, taking
approximations $f^{(N)}, g^{(N)} \leq 1$ converging to $f,g$ as above,
and noting that $g^{(N)}(x)dx \to g(x) dx$ in $\mathcal{M}(D)$, 
\[
Q^{N,D,\beta,\gamma}_{\mu,f^{(N)}+ g^{(N)}}
\left[ \Phi_{\phi,\psi} (U_{[0,\infty)},U^{\partial D}_{\infty}) \right]
\geq
Q^{N,D,\beta,\gamma}_{\mu + \beta g^{(N)} dx ,f^{(N)}}
\left[ \Phi_{\phi,\psi} (U_{[0,\infty)},U^{\partial D}_{\infty}) \right].
\]
Passing to the limit as above we find that
\[
Q^{D,\beta,\gamma}_{\mu,f+ g}
\left[ \Phi_{\phi,\psi} (U_{[0,\infty)},U^{\partial D}_{\infty}) \right]
\geq
Q^{D,\beta,\gamma}_{\mu + \beta g \, dx ,f}
\left[ \Phi_{\phi,\psi} (U_{[0,\infty)},U^{\partial D}_{\infty}) \right].
\]
This inequality between the Laplace functionals implies the desired
stochastic domination. 

\textbf{Proof of Lemma \ref{comparison1} part (i).} 
Consider the decomposition $u^-,u^+$ given in section
\ref{s5.2} example 1, with initial conditions $\mu= \mu^- + \mu^-$ and $f^{(N)} \to f$ 
in $L^1(D)$. Then pathwise
\begin{equation} \label{5.3.11}
 \Phi_{\phi,\psi}(u_{[0,\infty)}, u^{\partial D}_{\infty} ) 
= 
\Phi_{\phi,\psi} (u^-_{[0,\infty)} + u^+_{[0,\infty)}, u^{\partial D,-}_{\infty} 
+ u^{\partial D,+}_{\infty}).
\end{equation}
As above, the expectation of the left hand side converges, as $N \to \infty$, to 
$Q^{D,\beta,\gamma}_{\mu,f} \left[ \Phi_{\phi,\psi}(U_{[0,\infty)},U^{\partial D}_{\infty}) \right]$.
Using the law of $(u^+,v^+)$ conditional on $\mathcal{G}^-$, the limit of the expectation of 
the right hand side can be written as
\begin{eqnarray*}
&& \hspace{-.4in}  \lim_{N \to \infty} 
E \left[ \Phi_{\phi,\psi} (u^-_{[0,\infty)}, u^{\partial D,-}_{\infty} ) 
\, Q_{\mu^+, v^-_{\infty}}^{N,D,\beta,\gamma} 
\left[ \Phi(U_{[0,\infty)}, U^{\partial D}_{\infty}) \right] \right] \\
&=& \lim_{N \to \infty} \lim_{T \to \infty}
E \left[ \Phi_{\phi,\psi}(u^-_{[0,T]}, u^{\partial D,-}_{T}) 
\, Q_{\mu^+, v^-_{T}}^{N,D,\beta,\gamma} 
\left[ \Phi_{\phi,\psi}(U_{[0,\infty)}, U^{\partial D}_{\infty}) \right] \right] \\
& = & \lim_{T \to \infty} Q_{\mu^-,f}^{D,\beta,\gamma} \left[ \Phi_{\phi,\psi} (U_{[0,T]}, U^{\partial D}_{T}) 
\, Q_{\mu^+, f e^{-U(0,T)}}^{D,\beta,\gamma} 
\left[ \Phi_{\phi,\psi}(U_{[0,\infty)}, U^{\partial D}_{\infty}) \right] \right] \\
&=&  Q_{\mu^-,f}^{D,\beta,\gamma} 
\left[ \Phi_{\phi,\psi}(U_{[0,\infty)}, U^{\partial D}_{\infty}) 
\, Q_{\mu^-, f e^{-U(0,\infty)}}^{D,\beta,\gamma} \left[ \Phi_{\phi,\psi}( U_{[0,\infty)}, 
U^{\partial D}_{\infty})  \right] \right]. 
\end{eqnarray*}
The interchange in limits follows as above, using $v^-_T=v^-_{\infty}$ on
$\{u^-_T = 0\}$. The $N \to \infty$ limit holds using
the convergence of $(u^-_{[0,T]},u^{\partial D,-},v^-_T)$ given by
Theorem \ref{convergence} (ii) and (\ref{conv12}) above.  
The final identity between Laplace functionals implies the desired result. 

\textbf{Proof of Lemma \ref{comparison2}.} 
The proofs of both parts of this lemma follow closely the argument 
used for the proof of Lemma \ref{comparison1} part (i) above. 
For example, for part (i) take the decomposition $u^-,u^+$ given in section
\ref{s5.2} example 4, so that (\ref{5.3.11}) holds, and with
initial conditions $\mu$ and $f^{(N)} \to f$ in $L^1$.
Then argue as above, noting that  
\begin{eqnarray*}
&& \hspace{-.4in}  \lim_{N \to \infty}
E \left[ \Phi_{\phi,\psi}(u^-_{[0,T]}, u^{\partial D,-}_{T}) 
Q_{\beta^+(f-v^-_{T})dx, v^-_{T}}^{N,D,\beta,\gamma}  
\left[ \Phi_{\phi,\psi}(U_{[0,\infty)}, U^{\partial D}_{\infty}) \right] \right] \\
& = & Q_{\mu,f}^{D,\beta^-,\gamma}  
\left[ \Phi_{\phi,\psi} (U_{[0,T]}, U^{\partial D}_{T}) 
\, Q_{\beta^+ f(1-e^{-U(0,T)}) dx, f e^{-U(0,T)}}^{D,\beta,\gamma} 
 \left[ \Phi_{\phi,\psi}(U_{[0,\infty)}, U^{\partial D}_{\infty}) \right] \right] 
\end{eqnarray*}
using the fact that if $v^-_T$ converges in $L^1(D)$ then
$\beta^+(f-v^-_{T})dx$ converges in $\mathcal{M}(D)$.
Part (ii) is entirely similar starting with example 5 in section
\ref{s5.2}. 

\textbf{Proof of Lemma \ref{sMp}.} We can follow the previous argument closely, but 
we write this out since we have to be a little careful about the two 
domains. To indicate that we use the 
joint Laplace functional acting on measures on $D^+$ and $\partial D^+$, we denote it by
$\Phi^+_{\phi,\psi}$. Take the decomposition $u^-,u^+$ given in section
\ref{s5.2} example 3, so that 
$ \Phi^+_{\phi,\psi}(u_{[0,\infty)}, u^{\partial D^+}_{\infty} ) 
= 
\Phi^+_{\phi,\psi} (u^-_{[0,\infty)} + u^+_{[0,\infty)}, u^{\partial D^-,-}_{\infty} |_{ \partial D^+}+  u^{\partial D^+,+}_{\infty})$.
Taking expectations, the left hand side converges, as $N \to \infty$, to 
$Q^{D^+,\beta,\gamma}_{\mu,f} 
\left[ \Phi^+_{\phi,\psi}(U_{[0,\infty)},U^{\partial D^+}_{\infty}) \right]$.
For this proof only, set $\mu^- = \mu \I(D^-), \, \mu^+ = \mu \I(D^+ \setminus D^-)$, 
$f^-= f \I(D^-), \, f^+ = f \I(D^+ \setminus D^-)$ and
$f^{(N),-}= f^{(N)} \I(D^-), \, f^{(N),+} = f^{(N)} \I(D^+ \setminus D^-)$
The limit of the right hand side can be evaluated as
\begin{eqnarray*} 
&& \hspace{-.4in}  \lim_{N \to \infty} 
E \left[ \Phi^+_{\phi,\psi} (u^-_{[0,\infty)}, u^{\partial D^-,-}_{\infty} |_{\partial D^+}) 
\, Q_{\mu^+ + u^{\partial D^-,-}_{\infty} |_{D^+}, 
f^{(N),+} + v^-_{\infty} \I(D^-)}^{N,D^+,\beta,\gamma} 
\left[ \Phi^+(U_{[0,\infty)}, U^{\partial D^+}_{\infty}) \right] \right] \\
&=& \lim_{N \to \infty} \lim_{T \to \infty}
E \left[ \Phi^+_{\phi,\psi}(u^-_{[0,T]}, u^{\partial D^-,-}_{T} |_{\partial D^+}) 
\, Q_{\mu^+ + u^{\partial D^-,-}_{T}|_{D^+}, f^{(N),+} + v^-_{T} \I(D^-)}^{N,D^+,\beta,\gamma} 
\left[ \Phi^+_{\phi,\psi}(U_{[0,\infty)}, U^{\partial D^+}_{\infty}) \right] \right] \\
& = & \lim_{T \to \infty} Q_{\mu^-,f^-}^{D^-,\beta,\gamma}  
\left[ \Phi^+_{\phi,\psi} (U_{[0,T]}, U^{\partial D^-}_T |_{\partial D^+}) 
\, Q_{\mu^+ + U^{\partial D^-}_{T} |_{D^+}, 
f^+ + f e^{-U(0,T)} \I(D^-)}^{D^+,\beta,\gamma} 
\left[ \Phi^+_{\phi,\psi}(U_{[0,\infty)}, U^{\partial D^+}_{\infty}) \right] \right] \\
&=&  Q_{\mu^-,f^-}^{D^-,\beta,\gamma}  
\left[ \Phi^+_{\phi,\psi}(U_{[0,\infty)}, U^{\partial D^-}_{\infty}|_{\partial D^+}) 
\, Q_{\mu^+ + U^{\partial D^-}_{\infty}, 
f^+ + f e^{-U(0,\infty)} \I(D^-)}^{D^+,\beta,\gamma} 
 \left[ \Phi^+_{\phi,\psi}( U_{[0,\infty)}, U^{\partial D^+}_{\infty})  \right] \right]. 
\end{eqnarray*}
The limits above follow as in the previous examples, once we have established that
$ u^{\partial D^-,-}_T|_{\partial D^+} $ converges in $\mathcal{M}(\partial D^+)$ and
$ u^{\partial D^-,-}_T|_{D^+} $ converges in $\mathcal{M}(D^+)$.  Both of these 
follow from the convergence of $ u^{\partial D^-,-}_T $ in $\mathcal{M}(\partial D^-)$
under the hypothesis (\ref{nullset}), since this hypothesis ensures that the limit 
law of $U^{\partial D^-}_T$ under $Q_{\mu^-,f^-}^{D^-,\beta,\gamma}$  
does not charge the discontinuity set 
$S=(\partial D^- \cap \partial D^+) \cap \overline{(\partial D^- \setminus \partial D^+)}$
(check the first moment of $U^{\partial D^-}_T(S)$ is zero). 

\textbf{Proof of Lemma \ref{comparison1} part (ii).} 
To allow this case to follow from the same argument as the earlier examples
we need a slight improvement in Theorem \ref{convergence} part (ii). 
Taking the two stage construction from section
\ref{s5.2} example 3, we need to consider two parts of the approximate density, namely
\[
\overline{v}^-_t = \sum_{1 \leq k \leq L} \psi_k \I(k \not \in \mathcal{S}^-_t) \quad \mbox{and}
\quad \overline{v}^+_t = \sum_{L < k \leq K} \psi_k \I(k \not \in \mathcal{S}^-_t).
\]
Now consider a sequence of models indexed by $N$ in which
$f^{(N),\pm} \to f^{\pm}$ in $L^1(D)$ and where $f=f^- + f^+$. 
Then Theorem \ref{convergence} part (ii) can be extended to show that the 
law of the quadruple
$(u^{(N)}_{[0,T]}, u^{\partial D,(N)}_T, \overline{v}^{(N),-}_T, \overline{v}^{(N),+}_T)$ on
$\mathcal{M}(D) \times \mathcal{M}(\partial D) \times 
L^1(D) \times L^1(D)$ converges to the law of $(U_{[0,T]},U^{\partial D}_T,f^- \exp(-U(0,T)),
f^+ \exp(-U(0,T)))$ under $Q^{D,\beta,\gamma}_{\mu,f}$. 
This can be shown as in Theorem \ref{convergence}
once one notes that the analogue of (\ref{aveqn}) holds for 
$\overline{v}^{\pm}$, namely that there are martingales $m^{v,\pm}_t$ such that
\[
\overline{v}^{\pm}_t(x) 
= f(x) - \int^t_0 \overline{v}^{\pm}_s(x)
u_s(\hat{\psi}_x)ds + m^{v,\pm}_t(x).
\]
The initial conditions for
stage two of the 
construction can be expressed in terms of $\overline{v}^{\pm}$ as
$u^+_0 = \beta(f^+ - \overline{v}^+_{\infty}) dx $ and $v^+_0 = \overline{v}^-_{\infty} 
+ \overline{v}^+_{\infty}$
and the passage to the limit then follows the same lines as previous examples.
\section{Life} \label{s6}
\subsection{Embedded oriented percolation processes} \label{s6.1}
In the proofs of possible life, the main step is to construct 
a discrete one-dimensional oriented percolation (OP) 
$(\omega(j,k):(j,k) \in \mathcal{L})$ where
$\mathcal{L} = \{(j,k): j,k \in \Z, \;j \geq 1, \; \mbox{$j+k$ is even}\}$,
so that life occurs if the percolation is supercritical. 
The variables $\omega(j,k)$ will be defined in terms of the exit measures for
solutions to (\ref{spde}) on a series of boxes.  
We will choose a parameter $L>0$ controlling the length scale. 
In the grid of points $x_{j,k}^L=(3jL, 2kL)$ in $d=2$ or $x_{j,k}^L=(3jL, 2kL,0)$ in 
$d=3$, the $x_1$ direction will play the usual role of time, and the $x_2$ direction
the role of space for the comparison OP process. 
For this argument we therefore need to be in dimension $d \geq 2$.

We will freeze the mass of the solution as it exits the boxes
$D(n) = (-3nL,3nL)^d$ for $n =1,2,\ldots$ 
We define exit measures
$u^{n,\partial D(n)}_{\infty}$ inductively. Choose an initial
condition $\mu$ supported in $D(1)$ and let $(u^{1},v^{1})$ be a solution
on $D(1)$ started from $(\mu,1)$. Suppose that $(u^{j},v^{j})$ have been defined for
$j=1,\ldots,n$. Then, conditional on $ \sigma \{ u^{j}: j=1,\ldots,n \}$
we let $(u^{n+1},v^{n+1})$ 
be a solution on $D(n+1)$ started from 
\[
u^{n+1}_0 = u^{n,\partial D(n)}_{\infty}, \quad \mbox{and} \quad
v^{n+1}_0(x) = 
\left\{ \begin{array}{ll}
v^{n}_{\infty}(x) & x \in D(n), \\
1 & x \in D(n+1) \setminus D(n).
\end{array} \right.
\]
Repeated use of the spatial Markov property in Lemma \ref{sMp}
shows that $u^{n,\partial D(n)}_{\infty}$
has the same law as the exit measure $U^{\partial D(n)}_{\infty}$ under the law
$Q^{D(n),\beta,\gamma}_{\mu,1}$.

In the hyperplane $x_1=0$ we define the box $I_L= \{x: |x| \leq L, x_1=0\}$. 
Define, for $j,k \in \mathcal{L}$,  
$$
\tilde{\omega}(j,k) = \left\{
\begin{array}{ll} 
1 & \mbox{if $u^{j,\partial D(j)}_{\infty} \left( x_{j,k}^L+I_L \right) > M$,} \\
0 & \mbox{otherwise.}
\end{array} \right.
$$
Set $\tilde{\omega}(0,k)=0$ if $k \neq 0$ and $\tilde{\omega}(0,0) = I(\mu(I_L) \geq M)$.
Define for $(j,k) \in \mathcal{L}$, 
$$
\omega(j,k) = \left\{
\begin{array}{ll} 
1 & \mbox{if $\tilde{\omega}(j-1,k-1)=\tilde{\omega}(j-1,k+1)=0$,} \\
\tilde{\omega}(j,k) & \mbox{otherwise.}
\end{array} \right.
$$
We aim to show, for suitable choice of the length and mass scales $L$ and $M$, that 
$(\omega(j,k):(j,k) \in \mathcal{L})$ is a $3$-dependent oriented site percolation
with density at least $1-\epsilon$. We recall the definition of this from
Durrett \cite{durrett}: whenever $(j_n,k_n) \in \mathcal{L}$, for $1 \leq n \leq N$, satisfy
$|j_n-j_m| \geq 3$ or $|k_n-k_m| \geq 3$ for all $n \neq m$, then
\begin{equation} \label{OPhypothesis}
P[ \mbox{$\omega(j_n,k_n)=0$ for $1 \leq n \leq N$}] \leq \epsilon^N.
\end{equation}
We write $(0,0) \to (j,k)$ when there exist a sequence of points
$0=k_0,k_{1}, \ldots, k_{j-1},k_{j}=k$ so that
so that $|k_m - k_{m-1}|=1$ for
$0 < m \leq j$ and satisfying $\omega(m,k_m)=1$ for $1 \leq m \leq j$.
Then let $\mathcal{C}_0 = \{ (j,k): (0,0) \to (j,k) \}$ be the cluster of sites
connected to the origin. This definition is sufficient for 
the key property of percolation to hold: Theorem
4.1 of \cite{durrett} shows that when $\epsilon \leq \epsilon_0= 6^{-196}$ then
$P[|\mathcal{C}_0| < \infty] \leq 1/20 $,
where $|\mathcal{C}_0|$ is the cardinality of $\mathcal{C}_0$.
In particular, if the initial condition is such that
$\tilde{\omega}(0,0) =1$, that is $\mu(I_L) \geq M$, then 
\[
Q^{D(n),\beta,\gamma}_{\mu,1}[U^{\partial D(n)}_{\infty} \neq 0] =
P[ u^{n,\partial D(n)}_{\infty} \neq 0] \geq P[|\mathcal{C}_0| = \infty]
\]
is bounded away from zero, and by Lemmas 
\ref{deathcharacterization1} and \ref{deathcharacterization2}, possible life occurs.

It remains to check the hypothesis (\ref{OPhypothesis})
of being a $3$-dependent oriented site percolation. We do this inductively. 
By conditioning on $\mathcal{G}_M = \sigma\{ u^{j}:j=1,\ldots,M \}$ (and 
setting $\mathcal{G}_0$ to be the trivial sigma field) it suffices to
show an estimate of the form
\[
P[ \mbox{$\omega(M+1,k_n)=0$ for $1 \leq n \leq N_0$} \, |  \, \mathcal{G}_{M}] \leq \epsilon^{N_0}.
\]
whenever $|k_n-k_m| \geq 3$ for all $n \neq m$. 
If $\tilde{\omega}(M,k_n-1)=\tilde{\omega}(M,k_n+1)=0$ for some
$n$ then the conditional expectation is zero. So 
we may also restrict to the set where, for each $n$, there exists $\tilde{k}_n
\in \{k_n-1,k_n+1\}$ for which $\tilde{\omega}(M,\tilde{k}_n) =1$. 
We let $D(M,n)$ be the box $x^L_{M,\tilde{k}_n}+D_{3L}$. Note
that the boxes $D(M,n)$ for $n=1,\ldots,N_0$ are disjoint and contained in $D(M+1)$. 

We now apply the spatial Markov property Lemma \ref{sMp}
to the pairs of domains $D^- = D(M,n) \subseteq D(M+1) = D^+$, for 
$n=1,\ldots,N_0$ in succession. Since the domains $D(M,n)$ are disjoint, this amounts to 
constructing, conditional on $\sigma \{u^{j}:j=1,\ldots,M \}$, independent solutions 
$(u^{M+1,n},v^{M+1,n})$ on $D(M,n)$ with initial condition 
\[
u^{M+1,n}_0 = u^{M,\partial D(M)}_{\infty} \I(D(M,n)), \quad \mbox{and} \quad
v^{M+1,n}_0(x) = 
\left\{ \begin{array}{ll}
v^{M}_{\infty}(x) & x \in D(M,n) \cap D(M), \\
1 & x \in D(M,n) \setminus D(M),
\end{array} \right.
\]
for $n=1,\ldots,N_0$. Then an exit measure 
$\tilde{u}^{M+1,\partial D(M+1)}_{\infty}$ can be constructed by running one further 
solution, conditionally on these $N_0$ processes, that starts with the combined exit measures 
of $(u^{M+1,n}:n=1,\ldots,N_0)$ on their domains and any part of $u^{M,\partial D(M)}_{\infty}$ that is
unused, and run until the exit from $D(M+1)$. This leads to
a construction of an exit measure $\tilde{u}^{M+1,\partial D(M+1)}_{\infty}$ that, by the spatial Markov property, 
is equal in law to $u^{M+1,\partial D(M+1)}_{\infty}$, but that is pathwise 
larger than the sum of 
the exit measures $ \sum_{n=1}^{N_0} u^{M+1,n,\partial D(M,n)}_{\infty}$ restricted to $\partial D(M+1)$.
This in turn implies, since the solutions 
$(u^{M+1,n},v^{M+1,n})$ are conditionally independent, that 
\[
P[ \mbox{$\omega(M+1,k_n)=0$ for $1 \leq n \leq N_0$} \, | \, \mathcal{G}_{M}] \leq 
\prod_{n=1}^{N_0} P [ u^{M+1,n,\partial D(M,n)}_{\infty} (x_{M+1,k_n}^L+I_L) \leq M \; | \; \mathcal{G}_{M}] 
\]
By translational invariance, and the monotonicity in the initial measure, 
the condition (\ref{OPhypothesis}) will then follow from an estimate of the form
\begin{eqnarray} 
&& \hspace{-.3in} Q^{D_{3L},\beta,\gamma}_{\mu,f}
[U^{\partial D_{3L}}_{\infty}(x_{1,\pm 1}^L+I_L) \leq M] \leq \epsilon \nonumber \\
&& \hspace{.4in} \mbox{whenever $f \geq I(x_1 \geq 0)$ and $\mu$ is supported on 
$I_L$ with $\mu(1) \geq M$.} \label{OPhypothesis4}
\end{eqnarray} 
By monotonicity in $f$ it is enough to consider just the case $f=I(x_1 \geq 0)$.
We shall check this block estimate in the next two subsections for two different regions of parameter 
values $(\beta,\gamma)$.
\subsection{Proof that $\liminf_{\beta \to 0} \Psi(\beta)/\beta^2 > 0$ in dimension $d = 3$} \label{s6.2}
After applying the scaling in Lemma \ref{scaling}
with the choices $a=c=L \beta^{-1}, b = L^2 \beta^{-2}, e=1$, solutions to (\ref{spde}) with the parameter 
values $(\beta, L^{-2} \beta^2)$ in $d=3$ become solutions to the rescaled equation
(\ref{life1spde}) studied in the following lemma. Undoing this scaling, the lemma shows that the
block estimate (\ref{OPhypothesis4}) holds for solutions to (\ref{spde}) when $ \beta \in (0,1]$ and 
with the choices 
\[
\mbox{$\gamma=L_0^{-2}(\epsilon) \beta^2$, $L=L_0(\epsilon)\beta^{-1}$, and
$M=L_0^2(\epsilon) \beta^{-2}$.}
\] 
Choosing $\epsilon = \epsilon_0$, the value from in section \ref{s6.1} that ensures 
the OP is supercritical, the OP comparison 
implies possible life for the parameter values $(\beta, L_0^{-2}(\epsilon_0) \beta^2)$, and thus 
$\liminf \Psi(\beta)/\beta^2 \geq L_0^{-2}(\epsilon_0)$ 
as $\beta \to 0$.
\begin{Lemma} \label{life1keylemma}
Consider solutions, in $d=3$, on the domain $D_3 = (-3,3)^3$ to the equation
\begin{equation}  \label{life1spde}
\left\{ 
\begin{array}{l}
\partial_t u = \Delta u + L^2 \beta^{-1} uv - u +  \sqrt{u}\, \dot{W},  \\ 
\partial_t v = - L \beta^{-1} uv, \quad v_0(x) = \I(x \in H),
\end{array}
\right.
\end{equation}
where $H = \{x: x_1 \geq 0\}$. 
For any $\epsilon>0$ there exists $L_0=L_0(\epsilon) \geq 1$ so that, whenever 
$(u,v)$ is a solution to (\ref{life1spde}) with parameters
$L_0$ and $\beta \in (0,1]$, and with an initial condition $\mu$ supported in 
$I_1 = \{x:x_1=0,|x| \leq 1\}$ and satisfying $\mu(1) \geq 1$,
\[ 
P \left[ u^{\partial D_3}_{\infty}(x_{1,\pm 1}+I_1) \leq 1 \right] \leq \epsilon.
\]
where $x_{1,1}=(3,2,0)$ and $x_{-1,1}=(-3,2,0)$.
\end{Lemma}
\textbf{Remark.} The intuition for the lemma comes from 
taking $L$ large, ignoring for a moment the diffusion, killing and noise terms and reducing
to the ODE system $\partial_t x = L^2 \beta^{-1} xy, \; \partial_t y = -L \beta^{-1} xy$, where, 
for $L$ large, the solution converges fast to $x_{\infty} = x_0 + L y_0$.
Thus the nutrient leads to a large build-up of mass, which in turn leads 
to a high probability that the exit measure is large. 
The corresponding scaling in dimension $d=2$ fails to produce this large build up of mass.

\noindent
\textbf{Proof.} 
We first construct a DW process $u^-$ on $D_3$ with law $Q^{D_3,1}_{\mu}$. 
Then, conditional on $\sigma \{u^-\}$, we construct a solution
$(u^+,v^+)$ to (\ref{life1spde}) on $D_3$ with initial conditions
\[
u^+_0 = L \I(x \in H) (1- e^{-L \beta^{-1} u^-(0,\infty,x)}) dx 
\quad \mbox{and} \quad
v^+_0= \I(x \in H) e^{- L \beta^{-1} u^-(0,\infty,x)}. 
\]
Lemma \ref{comparison2}, after a suitable linear scaling, 
shows that $u^{\partial D_3,-}_{\infty} + u^{\partial D_3,+}_{\infty}$ has the same law as
the total exit measure $u^{\partial D_3}_{\infty}$ of a solution to 
(\ref{life1spde}) started at $\mu$.  
Lemma \ref{comparison2} also implies that the exit measure 
$u^{\partial D_3,+}_{\infty}$ is stochastically
larger than the exit measure $U^{\partial D_3}_{\infty}$ where $U$ has the law
$Q^{D_3,1}_{u^+_0}$ of a DW process, conditionally on $\sigma\{u-\}$. 
Choose smooth $h: \partial D_3 \to [0,1]$ so that $\{h>0\} = x_{1,1}+I_1$
(the case of $x_{1,-1}$ is symmetric). Then
\[
E \left[ e^{- u^{\partial D_3}_{\infty}(x_{1,1}+I_1)} \right]
\leq  E \left[ e^{- u^{\partial D_3}_{\infty}(h)} \right] \\
 \leq  E \left[ Q^{D_3,1}_{u^+_0} [ e^{-U^{\partial D_3}_{\infty}(h) } ] \right] \\
 =  E \left[ e^{- u^{+}_0(w)} \right] 
\]
where, using the DW Laplace functional from (\ref{LFD2}), 
$ \Delta w = \frac12 w^2 + w $ on $D_3$ and $w = h$ on $\partial D_3$. 
Let $c_0 = \inf \{w(x): x \in D_2 \cap H \}$. Then a comparison argument shows that 
$c_0 \in (0,\infty)$. Using the initial condition of $u^+_0$ we continue, for 
$ \beta \leq 1$ and $ L \geq 1$,
\begin{eqnarray*}
E \left[ e^{- u^{\partial D_3}_{\infty}(x_{1, 1}+I_1)} \right]  
& \leq & Q_{\mu}^{D_3,1} \left[ \exp \left(-  L \int_{H \cap D_3 } 
\left(1- e^{- L \beta^{-1} U(0,\infty,x)} \right) w(x) dx \right) \right] \\
& \leq &  Q_{\mu}^{D_3,1} \left[ \exp \left(- L c_0 \int_{H \cap D_2} 
\left(1- e^{- U(0,\infty,x)} \right) dx \right) \right]  =:  p(L,\mu).
\end{eqnarray*}
We will show that $p(L,\mu)$ converges to zero, as $L \to \infty$, 
uniformly over $\mu$ as stated in the lemma. Then a Chebychev argument finishes 
the proof.

It is not difficult to see that if $\mu(1) \neq 0$ and $\mu$ is 
supported on $I_1$ then $U_{[0,1]}(D_2 \cap H)>0$, $Q^{D_3,1}_{\mu}$ almost surely. It is possible to 
argue this directly from the Laplace functional, but it follows also from
results in the literature. Indeed the projection of the process onto
$\{x_2=\ldots=x_d=0\}$ is a DW process started at $c \delta_0$, so one can use
absolute continuity and the 
strong law of the logarithm for the closed support shown in Tribe (\cite{tribe91}).
Since $x \to U(0,\infty,x)$ is lower semi-continuous, the set
$\{x \in D_2 \cap H: U(0,\infty,x) >0 \}$ is open and almost surely 
non-empty by the above.
This implies that $p(L,\mu) \to 0$ as $L \to \infty$
for each non-zero $\mu$. 
Lemma \ref{comparison1} shows that $\mu \to p(L,\mu)$ is non-increasing, so we may restrict
to the set $\mathcal{M}_1(I)$ of measures supported in $I$ and with total mass one.
Note that $\mathcal{M}_1(I)$ is compact in the weak topology.  
Since $L \to p(L,\mu)$ is also non-increasing, the required uniform convergence 
follows by a Dini argument if we can show that $\mu \to p(L,\mu)$ is upper semi-continuous.
This will be true since $p(L,\mu)$ is the decreasing limit of $p(L,\mu,\epsilon)$ as $\epsilon \downarrow 0$, 
where
\[
p(L,\mu,\epsilon) = 
Q_{\mu,1}^D \left[ \exp \left(-  L c_0 \int_{D_2 \cap H} 
\left(1- e^{- U(\epsilon,\epsilon^{-1},x)} \right) dx \right) \right]
\] 
Increment estimates (similar to those in (\ref{densitytightness}) 
show that the laws $Q^{D}_{\mu,1}(U(\epsilon,\epsilon^{-1}) \in df)$ 
on $C(\overline{D_2},\R)$ are tight, uniformly over $\mu \in \mathcal{M}_1(I)$. 
This in turn implies that the map $\mu \to Q^{D}_{\mu,1}(U(\epsilon,\epsilon^{-1}) \in df)$ and 
hence also $\mu \to p(L,\mu,\epsilon)$ is continuous, which finishes the proof. \qed

\subsection{Proof that $\lim_{\beta \to \infty} \beta^{-1} \Psi(\beta) =1 $. } \label{s6.3}
We fix $\kappa \in (0,1)$ and show there exists $\beta_0= \beta_0(\kappa)$ so that
possible life occurs whenever $\beta \geq \beta_0$ and $\gamma = \kappa \beta$.
Since $\Psi(\beta) \leq \beta$ this implies the desired limit. 
\begin{Lemma} \label{largebetalifelemma1}
Fix $\kappa \in (0,1)$ and smooth 
$\eta: \overline{D_3} \to [0,1] $ compactly supported inside $D_3 \cap \{x:x_1>0\}$ and so that 
$\eta = 1$ on  $D_2 \cap \{x:x_1 \geq 1\}$.
Consider, in dimensions $d \geq 2$, solutions on the domain $D_3=(-3,3)^d$ to the equation
for a DW process given by 
\begin{equation}  \label{life10spde}
\partial_t u = \Delta u + \theta \left( \frac{1+\kappa}{2} \eta - \kappa \right) u  + \sqrt{\sigma u}\, \dot{W}.
\end{equation}
Given $\epsilon > 0$ we may choose 
$\theta_0(\epsilon),\sigma_0(\epsilon) > 0$ so that 
for any solution $u$ to (\ref{life10spde}) with parameters $\sigma_0$ and $\theta_0$, and 
initial condition $\mu \in \mathcal{M}_1(I_1)$,
\[ 
P \left[ u^{\partial D_3}_{1}(x_{1, \pm 1}+I_1) \geq 1 \right]  \geq 1 - \epsilon  
\]
where $x_{1,\pm 1}=(3, \pm 2,\ldots)$. 
\end{Lemma}

\noindent \textbf{Proof.}
Choose $h: \partial D_3 \to [0,1]$ so that $\{h>0\} = x_{1,1}+I_1$ (the case $x_{1,-1}$ being symmetric).
Using the Laplace functional from (\ref{LFD1}) and scaling, we have
$E[\exp(-u^{\partial D_3}_{1}(h))] = \exp(-\mu(\phi_1))$ where
\begin{equation} \label{lifetemp10}
\left\{ \begin{array}{l}
\partial_t \phi = \Delta \phi + \theta \left( \frac{1+\kappa}{2} \eta - \kappa \right) \phi 
- \frac{\sigma}{2} \phi^2  \quad \mbox{on $D_3$}, \\
\quad \mbox{$\phi = h$ on $\partial D_3$ for $t \geq 0, \quad$ and $\phi_0 = 0$ on $D_3$.} 
\end{array} \right.
\end{equation}
As $\sigma \to 0$ the solution $\phi_1$ converges to $\overline{\phi_1}$ where $\overline{\phi}$ is 
the solution of (\ref{lifetemp10}) with $\sigma =0$. Since $\phi$ is decreasing in the parameter
$\sigma$, the convergence is uniform over $x \in I_1$. 

We claim that $\overline{\phi}_1(x) \to \infty$
as $\theta \to \infty$ and that the convergence is uniform over $x \in I_1$. 
One way to see this is via the Feynman-Kac representation for 
$\overline{\phi}$, namely 
\[
\overline{\phi}_1(x) = E_x \left[ h(X_{\tau}) \I( \tau \leq 1) e^{\theta \int^{\tau}_0 
\left( \frac{1+\kappa}{2} \eta(X_s) - \kappa \right) ds} \right]
\] 
where $X$ is a Brownian motion (with generator $\Delta$) and $\tau= \inf \{t: X_t \in \partial D_3\}$.
Note that if $X_s \in D_2 \cap \{x:x_1 \geq 1\}$ then $(\frac{1+\kappa}{2} \eta(X_s) - \kappa)  \geq (1-\kappa)/2$.
Now we use the support theorem for Wiener measure. There is an open set of paths $O_x \subseteq C([0,1],\R^d)$ 
where $f \in O_x$ satisfy $f_0=x$, $f$ exits 
$D_3$ before time $1$ and $\int^{\tau}_0 (\frac{1+\kappa}{2} \eta(f_s) - \kappa) ds > (1-\kappa)/4$. 
Then $\overline{\phi}_1(x) \geq \exp( \theta (1-\kappa)/4) E_x [ h(X_{\tau}) \I(X_{\cdot} \in O_x) ]$.
Moroever $O_x$ can be chosen so that $E_x[h(X_{\tau}) \I(X_{\cdot} \in O_x)]$ is bounded below uniformly over $x \in I_1$,
establishing the claim. Then we may choose $\theta_0(\epsilon)$ and then $\sigma_0(\epsilon)$ so that
\begin{eqnarray*}
P \left[u^{\partial D_3}_1(x_{1,1}+I_1) \leq 1 \right] & \leq & P \left[u^{\partial D_3}_1(h) \leq 1 \right] \\
& \leq & e \, e^{-\mu(\phi_1)} \quad \mbox{by Markov's inequality} \\
& \leq & e \, e^{-\mu(\overline{\phi}_1)} + \frac{\epsilon}{2} \quad \mbox{by the suitable choice of $\sigma_0(\epsilon)$}  \\
& \leq & \epsilon \quad \mbox{by the suitable choice of $\theta_0(\epsilon)$}
\end{eqnarray*}
which completes the proof. \qed

\begin{Lemma} \label{largebetalifelemma2}
Consider, in dimensions $d = 2$ or $3$, solutions on the domain $D_3$ to the equation
\begin{equation}  \label{life20spde}
\left\{ 
\begin{array}{l}
\partial_t u = \Delta u +  \theta u \left(v - \kappa \right) + \sqrt{\sigma u}\, \dot{W},  \\ 
\partial_t v = - \delta u v, \quad v_0 = \I(x_1 \geq 0).
\end{array}
\right.
\end{equation} 
Given $\epsilon > 0$ we may choose 
$\theta_1(\epsilon), \sigma_1(\epsilon),\delta_1(\epsilon) > 0$ so that 
for any solution $(u,v)$ to (\ref{life20spde}), with parameter values
$\theta_1$, $\sigma_1$ and $\delta \leq \delta_1$, and with initial condition $\mu$ 
supported on $I_1$ satisfying $\mu(1) \geq 1$, 
\[ 
P \left[ u^{\partial D_3}_{\infty}(x_{1, \pm 1}+I_1) \geq 1 \right]  \geq 1 - \epsilon  
\]
where $x_{1,\pm 1}=(3,\pm 2,\ldots)$.  
\end{Lemma}

\noindent \textbf{Proof.}
By the monotonicity in $\mu$ given in Lemma \ref{comparison1}, 
it is enough to consider $\mu \in \mathcal{M}_1(I_1)$.
Fix $\theta_1 = \theta_0(\epsilon/2)$ and $\sigma_1 = \sigma_0(\epsilon/2)$ from the previous lemma. 
Let $\eta$ be as in the previous lemma and let $K$ be the closed support of $\eta$.  
We consider a solution to (\ref{life20spde}) on path space and set
\[
\tau_{\delta} = \inf \{t: \mbox{$\delta U(0,t,x) \geq (1-\kappa)/2$ for some $x \in K$} \}. 
\]
Note that $V(t,x) = \exp(-\delta U(0,t,x)) \geq 1 - (1-\kappa)/2 $ for $t \leq \tau_{\delta}$ and $x \in K$. 
This implies that 
\[
V(t,x)-\kappa \geq \frac{1+\kappa}{2} \eta(x) - \kappa \quad \mbox{for $t \leq \tau_{\delta}$ and $x \in D_3$. }
\]
Informally, upto the time $\tau_{\delta}$ the process dominates a solution to (\ref{life10spde}).
More carefully, we may expand the process $\exp(-u_s(\phi_{1-s}) - u^{\partial D_3}_s(\phi_s))$ over the interval
$s \in [0,\tau_{\delta} \wedge 1]$, with $\phi$ as in the previous lemma, to obtain an upper bound on the  
the Laplace functional $E[\exp(-u^{\partial D_3}_1(h)) \I( \tau_{\delta} \geq 1)]$. 
This shows that
$P[u^{\partial D_3}_{\infty}(x_{1, \pm 1}+I_1) \leq 1, \tau_{\delta} \geq 1] \leq \epsilon/2$ for all 
$\mu \in \mathcal{M}_1(I_1)$ (and for any choice of $\delta>0$). 

It remains to estimate $P[\tau_{\delta} <1]$. Informally, we can control this by choosing $\delta$ small and the fact
that the process is dominated by a DW($D_3, -\theta_1$) process. 
Rather than use the natural pathwise comparison, which is messy (although possible) to establish in this non-Lipschitz setting,
we exploit a change of measure argument. Let $Q_{\mu}$ be the law of $(u, u^{\partial D_3})$ on path space
$\Omega_{D_3}$ of the DW process satisfying $\partial_t u = \Delta u + \sqrt{\sigma_1 u} \, \dot{W}$ starting at
$\mu$. Set
\[
M_t = \sigma_1^{-1} \theta_1 \int^t_0 \int \left( e^{\delta  U(0,s,x)} - \kappa \right) Z(dx,ds)
\]
and let $\mathcal{E}_t(M)$ be the associated stochastic exponential. Arguing as in section \ref{s3.2}, under the the measure
$dP|_{\mathcal{U}_t} = \mathcal{E}_t(M) Q_{\mu}$, the process $U$ is a solution to (\ref{life20spde}). 
Then
\begin{equation}
P[ \tau_{\delta} < 1] = Q_{\mu} \left[ \mathcal{E}_{\tau_{\delta} \wedge 1}(M) 
\I(\tau_{\delta} <1) \right]  \leq  \left( Q_{\mu} \left[ \mathcal{E}_{\tau_{\delta} \wedge 1}(pM) 
e^{((p-1)/2)[M]_{\tau_{\delta} \wedge 1}} \right]\right)^{1/p} 
\left( Q_{\mu} [ \tau_{\delta} <1]\right)^{1/q} \label{temp246}
\end{equation}
for a dual pair $p^{-1} + q^{-1} =1$.
Note that $[M]_t \leq \sigma_1^{-2} \theta_1^2 U_{[0,t]}(1)$ for $t \leq \tau_{\delta}$. 
Reversing the change of measure, under the law $\mathcal{E}(pM) dQ_{\mu}$ the process has a mass creation at most at rate
$p \theta_1$ up to time $\tau_{\delta}$. For small enough $\theta>0$ the exponential moments
$Q^{D_3, -p \theta_1}_{\mu}[\exp(\theta U_{[0,1]}(1))]$ are finite, and bounded uniformly over $\mu \in \mathcal{M}_1(I_1)$.
Thus we may choose $p=p(\epsilon)>1$ so that the first term on the right hand side of (\ref{temp246}) is 
bounded uniformly over $\mu \in \mathcal{M}_1(I_1)$, and it remains to estimate $Q_{\mu}[\tau_{\delta} <1]$. 
Note that 
\[
Q_{\mu}[\tau_{\delta} <1] \leq Q_{\mu} \left[ \min\{1, 2 \delta (1-\kappa)^{-1} \sup_{x \in K} U(0,1,x) \} \right] := p(\delta,\mu).
\]
The continuity of the occupation density $U(0,1,x)$ for $x \in K$ (the finite spead of support ensures that
the solution does not hit $K$ immediately) implies that $\sup_{x \in K} U(0,1,x)$ is finite, and so that
$p(\delta,\mu) \downarrow 0$ as $ \delta \downarrow 0$. But $p(\mu,\delta)$ is the deceasing limit of $p(\delta,\mu,\epsilon)$ as $\epsilon \downarrow 0$ where
\[
p(\delta,\mu,\epsilon)= Q_{\mu} \left[ \min\{1, 2 \delta (1-\kappa)^{-1} \sup_{x \in K} U(\epsilon,1,x) \} \right].
\]
Arguing as in the proof of Lemma \ref{life1keylemma}, the functions $p(\delta,\mu,\epsilon)$ are continuous in $\mu$ and hence $p(\delta,\mu)$
is upper semicontinuous, ensuring that the convergence $p(\delta,\mu) \to 0$ is uniform over
$\mathcal{M}_1(I_1)$. This allows us to choose $\delta_1$ to ensure that 
$P[ \tau_{\delta} < 1] \leq \epsilon/2$ for all $\mu \in \mathcal{M}_1(I_1)$ and completes the proof. \qed

Now chooose $\epsilon= \epsilon_0$ the value that ensure the OP process is supercritical. 
Using the scaling lemma \ref{scaling} for a solution to (\ref{spde}) with the choices $b= \theta_1(\epsilon_0) \beta^{-1}$, $a=b^2$, 
$c^d=ab\sigma_1(\epsilon_0)$ and $e=1$,
the scaled solution $\tilde{u}$ solves (\ref{life20spde}) with the parameter values $\theta_1$, $\sigma_1$ and 
$\delta = \theta_1 \beta^{-1}$. Choosing $\beta$ sufficiently large that $\delta \leq \delta_1(\epsilon_0)$ we see that the block
estimate (\ref{OPhypothesis4}) holds, ensuring possible life. 

\noindent \textbf{Remark.} The correct large $\beta$ asymptotics can be seen from the following formal 
argument.
Choose $\gamma = \beta - \theta \beta^{2/(6-d)}$ for some $\theta >0$. Applying the scaling in 
Lemma \ref{scaling} with the choices $a= \beta^{-(d-2)/(6-d)}$, $b=\beta^{-2/(6-d)}$, $c= \beta^{-1/(6-d)}$
and $e=1$ we find the scaled equation solves
\[
\left\{ 
\begin{array}{l}
\partial_t u =  \Delta u + \beta^{(4-d)/(6-d)} u \left(v - (1-\theta \beta^{-(4-d)/(6-d)}) \right) + \sqrt{u} \, \dot{W}, \\ 
\partial_t v = - \beta^{-(4-d)/(6-d)} uv, \quad v_0 = 1.
\end{array}
\right.
\]
Letting $v = 1- \beta^{-(4-d)/(6-d)} \hat{v}$ we may rewrite this equation as
$\partial_t u =  \Delta u + \theta u - u \hat{v} + \sqrt{u} \, \dot{W}$ and 
$\partial_t \hat{v} = u(1- \beta^{-(4-d)/(6-d)} \hat{v})$. This suggests that as $\beta \to \infty$ the system 
converges to the one parameter system
\begin{equation} 
\left\{ 
\begin{array}{l}
\partial_t u =  \Delta u + \theta u - uv + \sqrt{u} \, \dot{W}, \\ 
\partial_t v = u, \quad v_0 = 0.
\end{array}
\right.
\label{6.310}
\end{equation}
There is monotonicity of total occupation and exit measures in $\theta$ and one expects that there is a critical 
$\theta_c$ so that solutions to (\ref{6.310}) die for $\theta < \theta_c$ and may live for $\theta > \theta_c$. 
It is therefore reasonable to conjecture that $\lim_{\beta \to \infty} \beta^{-2/(6-d)} (\beta - \Psi(\beta)) = \theta_c$.
There is no obvious line of proof for this conjecture, since 
it is not obvious that a life or death block estimate will hold for $\theta$ close to $\theta_c$. 
However such block estimates should hold for sufficiently large and small $\theta$
leading to upper and lower asymptotics of the same order. However note that the block construction 
we have presented in section \ref{s6.2} uses initial nutrient condition $f = \I(x_1 \geq 0)$, and the above 
intuition does not apply as such. One needs a more careful block construction where the 
nutrient level is controlled, and this in turn seems to require a stronger version of the spatial Markov property.
We will present these details in a subsequent paper.  

We remark that a similar one parameter scaling limit should hold for the small $\beta$ asymptotic. 
Indeed letting $\beta \to 0$ in (\ref{life1spde}) suggest that the small $\beta$ behaviour 
when $\gamma = L^{-2} \beta^2$ should approximate the one parameter system
\[
\partial_t u = \Delta u + L \delta_{\sigma(x)} -u + \sqrt{u} \, \dot{W} \quad \mbox{where $\sigma(x) = \inf\{t: U(0,t,x)>0\}$.}
\]
The analogous conjecture is that that this equation makes sense, that it has has a critical value $L_c$, and 
that the limit $\Psi(\beta)/\beta^2$ exists and equals $L_c$. 
\section{Death} \label{s7}
\subsection{Death in dimension $d=1$} \label{s7.1}
We will show that death is certain when $\gamma=0,\beta \geq 0$ and $d=1$. 
\begin{Lemma} \label{d=1deathlemma}
There are closed bounded intervals $I(\mu) \subseteq (\-\infty,\infty)$, indexed by $\mu \in \mathcal{M}(\R)$, 
and probabilities $p_{\beta} >0$, for $\beta \geq 0$, so that when $\mu$ is compactly supported
\[
Q^{\R,\beta,0}_{\mu,1} [U_{[0,\infty)}(I(\mu)^c)=0] \geq p_{\beta}.
\]
\end{Lemma}
By arguing as in the proof of Lemma \ref{deathcharacterization2}, one may check that
\begin{equation}
Q^{\R,\beta,0}_{\mu,f} [U_{[0,\infty)}(I^c)=0] = Q^{(-M,M),\beta,0}_{\mu,f|_{(-M,M)}} [U_{[0,\infty)}(I^c)=0]
\label{oddremark}
\end{equation}
whenever $I$ is a closed interval supporting $\mu$ and satisfying $I \subseteq (-M,M)$. 
We define stopping times iteratively as follows. Let
$\tau_1 = \inf\{t: U_t(I(\mu)^c)>0\}$ (with $\inf\{\emptyset\}=+\infty$). Supposing
$\tau_n <\infty$ we set
\[
\tau_{n+1} = \inf\{t \geq \tau_n: U_t(I(U_{\tau_n})^c) >0\}.
\]
If $\tau_n =\infty$ then we let $\tau_{n+1}=\infty$. 
Set $p(\mu,f) = Q_{\mu,f}^{\R,\beta,0}[ U_{[0,\infty)} (I(\mu)^c)>0]$. The decomposition in Lemma \ref{comparison1}
and (\ref{oddremark}) show that $f \to p(\mu,f)$ is non-decreasing. By the lemma above,
$Q_{\mu,1}^{\R,\beta,0}[\tau_1 < \infty] = p(\mu,1) \leq 1-p_{\beta}$.
By the strong Markov property at $\tau_n$ we have on the set $\{\tau_n <\infty\}$
\begin{eqnarray*}
Q_{\mu,1}^{\R,\beta,0} [\tau_{n+1} < \infty| \mathcal{U}_{\tau_n}] 
& = & p(U_{\tau_n}, \exp(-U(0,\tau_n))) \\
& \leq & p(U_{\tau_n},1) \leq 1-p_{\beta}.
\end{eqnarray*}
Then iterating the lemma shows that
\[
Q^{\R,\beta,0}_{\mu,1}[\mbox{$U_{[0,\infty)}$ is not compactly supported}] \leq 
Q^{\R,\beta,0}_{\mu,1}[\tau_n < \infty] \leq (1-p_{\beta})^n 
\]
completing the proof of death
by the characterization in Lemma \ref{deathcharacterization1}. 

\noindent 
\textbf{Proof of Lemma \ref{d=1deathlemma}.} 
Choose a closed bounded interval 
$\tilde{I} = \tilde{I}(\mu) \subseteq (-\infty,\infty)$ so that
$Q^{\R}_{\mu,0}[U_{[0,\infty)}(\tilde{I}^c) =0 ] \geq 1/2$. This is possible since 
paths die out and are compactly supported. Without loss we may assume $\tilde{I}$ has length at least $2$. 
Now set $I = \tilde{I} + [-1,1]$. Choose $M$ so that $I \subseteq (-M,M)$. 
By (\ref{oddremark}) we must show
$Q^{(-M,M),\beta,1}_{\mu,0}[U_{[0,\infty)}(I^c) =0] \geq p_{\beta}$.  
Construct processes $(u^-,u^+,v^+)$ so that 
$u^-$ is a DW($(-M,M),0$) process started at $\mu$, and conditional on $\sigma\{u^-\}$ the process
$(u^+,v^+)$ is a solution to (\ref{spde}) on $(-M,M)$ with parameters $(\beta,0)$ and initial conditions
$ u^+_0 = \beta (1-e^{-u^-(0,\infty,x)}) dx $ and  $v^+_0 = e^{-u^-(0,\infty)}$. 
Then Lemma \ref{comparison2} (i) shows that $u^-_{[0,\infty)} + u^+_{[0,\infty)}$ has the same law as
$U_{[0,\infty)}$ under $Q^{(-M,M),\beta,0}_{\mu,1}$. 
On the set $\{u^-_{[0,\infty)}(\tilde{I}^c) =0\}$ we have
$u^+_0 \leq \beta \I(x \in \tilde{I}) dx$. When this occurs, Lemma \ref{comparison3} shows that
the total occupation measure
$u^+_{[0,\infty)}$ becomes stochastically larger if we replace the initial conditions $(u^+_0,v^+_0)$ by
$(\beta \I(x \in \tilde{I}) dx, \I((-M,M) \setminus \tilde{I}))$. 
Therefore on $\{u^-_{[0,\infty)}(\tilde{I}^c) =0\}$,
\begin{eqnarray*}
P[u^+_{[0,\infty)}(I^c(\mu)) =0 | \, \sigma\{u^-\}] &=&  
Q^{(-M,M),\beta,0}_{u^+_0,v^+_0} [U_{[0,\infty)}(I^c) =0] \\
& \geq &  Q^{(-M,M),\beta,0}_{\beta \I(x \in \tilde{I}) dx, \I((-M,M) \setminus \tilde{I})} 
[U_{[0,\infty)}(I^c) =0] \\
& = & Q^{\R,\beta,0}_{\beta \I(x \in \tilde{I}) dx, \I(\tilde{I}^c)} 
[U_{[0,\infty)}(I^c) =0] \quad \mbox{ by (\ref{oddremark}) } \\
& \geq & Q^{\R, - \beta \I(\tilde{I}^c)}_{\beta \I(x \in \tilde{I}) dx} [U_{[0,\infty)}(I^c) =0]  \\
& \geq & \exp(-\beta \int_{\tilde{I}} w(x) dx) \quad \mbox{using (\ref{globaltsupport2})}
\end{eqnarray*}
provided that $w$ solves 
\begin{equation} \label{d=1comp}
\mbox{$\Delta w \leq \frac12 w^2 - \eta w$ on $\mbox{int}(I)$ and $w(x) \uparrow \infty$ as $x \to \partial I$} 
\end{equation}
with smooth $\eta \geq 0$ satisfying $\eta \geq \beta $ on $I \setminus \tilde{I}$. The proof will be 
complete if we find such a function $w$ and show that there is an 
upper bound on $\int_{\tilde{I}} w(x) dx$ that is independent of the 
length of the interval $I$. 

Choose smooth non-increasing $\overline{\eta}:[0,\infty) \to [0,1]$ so that $\overline{\eta}(x) = 1$ for $ x \in [0,1]$ and 
$\overline{\eta}(x) = 0$ for $ x \geq 2$.
There is a unique non-increasing $\overline{w}(x) \geq 0$, for $x \in (0,\infty)$ that solves
\[
\mbox{$\Delta \overline{w} = \frac12 \overline{w}^2 - 2 \beta \overline{\eta} \, \overline{w}$ on $(0,\infty), \quad$  
$\overline{w}(x) \to \infty$ as $x \to 0$ and $\overline{w}(x) \to 0$ as $x \to \infty$}.
\]
Moreover the decay as $x \to \infty$ ensures that $ \int_1^{\infty} \overline{w}(x) dx < \infty$. 
Suppose $I(\mu) = [a,b]$ and $\tilde{I}=[a+1,b-1]$, where $b \geq a+4$. Then $w(x) = 
\overline{w}(x-a) + \overline{w}(b-x)$ satisfies the required properties (\ref{d=1comp})
with the choice $\eta(x) = \beta \overline{\eta}(x-a) + \beta \overline{\eta}(b-x)$. \qed

\noindent \textbf{Remark.} The extinction estimates for a DW($\R,\gamma$) can be used to show that
$Q^{\R,\gamma}_{\delta_0}[ U_{[0,\infty)}([-L,L]^c)>0] \sim \exp(-\gamma^{1/2} L)$ for large $L$. 
This in turn yields the finiteness of small positive exponential moments of 
$\int (1- \exp(-2U(0,\infty,x)) dx$ which implies, via Novikov's criterion, that the exponential martingale
$t \to \mathcal{E}_t(M^{\beta,0,1})$ is uniformly integrable for small $\beta$ and hence certain death.
We omit the details since this simple method does not seem to be applicable for large $\beta$ in $d=1$ 
or any $\beta>0$ in dimension $d>1$. 
\subsection{A block estimate for death} \label{s7.2}

By Lemmas \ref{deathcharacterization1} and \ref{deathcharacterization2}, to establish
certain death for some parameter values
$(\beta,\gamma)$ it is sufficient to choose a non-zero compactly supported $\mu$ and show that
\begin{equation} \label{deathsufficient}
\lim_{n \to \infty} Q^{\R^d, \beta,\gamma}_{\mu,1} 
\left[ U^{\partial D_n}_{\infty} = 0 \right] = 1.
\end{equation} 
It can be simpler to verify the following block estimate, which will
imply (\ref{deathsufficient}) via a simple iteration reminiscent of the 
more direct branching process comparison used in \cite{mueller+tribe}.
\begin{Lemma} \label{deathsufficiencylemma}
Suppose that there exists $L,M>0$ so that for all $\mu$ supported in $\overline{D}_L=[-L,L]^d$
and satisfying $\mu(1) \leq M$ the following bounds hold:
\begin{equation}
Q^{D_{3L},\beta,\gamma}_{\mu,1} \left[ U^{\partial D_{3L}}_{\infty} \neq 0\right] < \epsilon_0 \quad
\mbox{and} \quad 
Q^{D_{3L},\beta,\gamma}_{\mu,1} \left[ U^{\partial D_{3L}}_{\infty}(1) \right] < \epsilon_0 M,
\label{deathblockestimates}
\end{equation}
where $\epsilon_0 = 1/(4 \cdot 3^d)$. Then (\ref{deathsufficient}) holds for $\mu = M \delta_0$  
and hence certain death holds for $(\beta,\gamma)$. 
\end{Lemma}

\noindent \textbf{Proof.} 
We want to break a measure $\mu \in \mathcal{M}(\R^d)$ into submeasures supported on 
blocks of size $L$ and of mass at most $M$. We consider
translates $x+[-L,L)^d$ for $x \in 2 L \Z^d$ and break 
the restriction of $\mu$ on $x+[-L,L)^d$ into at most $\mbox{Int}[M^{-1} \mu(x+[-L,L)^d)]$ 
parts, each of measure at most $M$ 
(here $\mbox{Int}[z]$ is the smallest integer greater than or equal to $z$). 
We require only that this splitting is done in a
measurable manner. 
Let $N_{L,M}(\mu)$ be the number of pieces into which $\mu$ is thus decomposed,
that is 
\[
N_{L,M}(\mu) = \sum_{x \in L\Z^d} \mbox{Int} \left[ M^{-1}\mu(x+[-L,L)^d) \right].
\]
We now inductively define a sequence $(X_n:n \geq 0)$ of random measures. 
Let $X_0= M \delta_0$.
Let $X_1=u^{(1),\partial D_{3L}}_{\infty}$ for a solution $(u^{(1)},v^{(1)})$
to (\ref{spde}) on $D_{3L}$ with initial conditions $u^{(1)}_0 = X_0$ and $v^{(1)}_0 = 1$.  
The spatial Markov property (Lemma \ref{sMp}) and the monotonicity of exit measures with respect to
the initial nutrient level (Lemma \ref{comparison1} (ii)), imply that for $R\geq 3L$
\[
 \hspace{-.3in} Q^{D_R,\beta,\gamma}_{X_0,1} \left[ U^{\partial D_R}_{\infty} =0 \right] 
\geq
E \left[  Q^{D_R,\beta,\gamma}_{X_1,1} [ U^{\partial D_R}_{\infty} =0 ] \right].
\]
Then inductively, provided that $X_{n-1} \neq 0$, 
 we choose one of the submeasures $\tilde{X}_{n-1}$ that form the decomposition
of $X_{n-1}$ as in the first paragraph (in some measurable manner), and denote 
its supporting cube as $x_{n-1} + [-L,L)^d$. Let $D(n) = x_{n-1} + (-3L,3L)^d$. Conditionally independent
of $(u^{(j)},v^{(j)}:j=1,\ldots,n-1)$ we run a solution $(u^{(n)},v^{(n)})$ 
to (\ref{spde}) on $D(n)$ with initial
conditions $u^{(n)}_0 = \tilde{X}_{n-1}$ and $v^{(n)}_0 = 1$. We then set
\begin{equation}
X_n = X_{n-1} - \tilde{X}_{n-1} + u^{(n), \partial D(n)}_{\infty}.
\label{d=1eqn1}
\end{equation}
When $X_{n-1}=0$ we set $X_n=0$.

We apply the decomposition Lemma \ref{comparison2} (i) to the domain $D(n)$ to see that the 
total exit measure on $\partial D(n)$ started from $(X_{n-1}|_{D(n)},1)$ can be constructed in two parts
as $u^{(n),\partial D(n)}_{\infty} + \tilde{u}^{(n),\partial D(n)}_{\infty}$ where
$\tilde{u}^{(n)}$ is a solution on $D(n)$ started at $X_{n-1}|_{D(N)} - \tilde{X}_{n-1}$ and initial 
nutrient level $\exp(-u^{(n)}(0,\infty))$. The monotonicity of exit measures in the nutrient level
shows that $\tilde{u}^{(n),\partial D(n)}_{\infty}$ is stochastically smaller than the total exit measure
of a solution on $D(n)$ started at $(X_n,1)$. 
Noting that $X_n$ is certainly supported inside $\overline{D}_{(2n+1)L}$, the spatial Markov property 
then gives, for $R \geq (2n+1)L$,
\begin{equation} \label{temp104}
Q^{D_R,\beta,\gamma}_{X_0,1} \left[ U^{\partial D_R}_{\infty} =0 \right] 
\geq
E \left[  Q^{D_R,\beta,\gamma}_{X_n,1} [ U^{\partial D_R}_{\infty} =0 ] \right].
\end{equation}
Using the simple inequality $\mbox{Int}[z] \leq z + \I(z>0)$, the hypotheses of this
lemma imply, when $X_n \neq 0$, that 
\[
E[N_{L,M}(u^{(n), \partial D(n)}_{\infty})] \leq 3^d (\epsilon_0 + \epsilon_0) \leq 1/2.
\]
Note from (\ref{d=1eqn1}) that $N_{L,M}(X_n) \leq N_{L,M}(X_{n-1}) - 1 + N_{L,M}(u^{(n),\partial D(n)}_{\infty})$
when $X_{n-1} \neq 0$. Then $n \to N_{L,M}(X_n)$ is a $\N$ valued  supermartingale
whose only possible limit point is zero, which implies that 
$X_n=0$ for large $n$, almost surely. Then (\ref{temp104}) implies the sufficient condition
(\ref{deathsufficient}) for certain death. \qed

As an illustration of this block condition, we shall use it to show 
for any $\beta > 0$ there exists $\gamma < \beta$ so that certain death occurs
for $(\beta,\gamma)$. This implies that $\Phi(\beta) <\beta$. The following lemma 
will imply that the block estimates will hold when $\gamma = \beta$, and a perturbation 
argument will ensure they hold $\gamma$ sufficiently close to $\beta$. 
\begin{Lemma} \label{death2keylemma}
Consider solutions in the domain $D_3=(-3,3)^d$ to the equation
\begin{equation}  \label{death1spde2}
\left\{ 
\begin{array}{l}
\partial_t u = b^{-(4-d)/d} \Delta u + b \beta u(v - 1) + \sqrt{u}\, \dot{W},  \\ 
\partial_t v = -uv, \quad v_0 = 1.
\end{array}
\right.
\end{equation}
Then, given $\epsilon>0$, there exists $b_0 = b_0(\beta,d,\epsilon)>0$ 
so that for $b \geq b_0$, and whenever $u_0=\mu$ is supported in
$\overline{D_1}$ and satisfies $\mu(1) \leq 1$, 
\[ 
E \left[ u^{\partial D_3}_{\infty}(1) \right] < \epsilon
\quad \mbox{and} \quad 
P \left[ u^{\partial D_3}_{\infty} \neq 0 \right] < \epsilon.
\]
\end{Lemma}

\noindent 
\textbf{Remark.} Note that (\ref{death1spde2}) comes from (\ref{spde}) via the 
scaling lemma with the choices $e=1, a=b, c=b^{2/d}$. Undoing this scaling shows 
that the block estimates (\ref{deathblockestimates}) hold when $\gamma = \beta$ and 
$L=b^{2/d}, M=b$ with $b = b_0(\beta,d,\epsilon_0)$.

\noindent \textbf{Proof.} 
We first consider solution to (\ref{death1spde2}) on the smaller domain $D_2$ and estimate the 
expected total exit measure. 
Using Lemma \ref{comparison3} we can and shall, at the expense of obtaining a stochastically 
larger exit measure, change the initial condition to $\mu + \mathbf{I}(D_2)(dx)$
and the initial nutrient to $v_0 = (1- (b \beta)^{-1}) \I(D_2)$. With these changes 
the reaction term $b \beta u(v-1)$ is at most $-u$. 
Let $\phi_{b,r}(x)$ solve $b^{-(4-d)/d} \Delta \phi_{b,r} = \phi_{b,r}$ on the domain 
$D_r$ with boundary conditions $\phi_{b,r}=1$ on $\partial D_r$. Calculus shows for $r>2$ that 
$u_t(\phi_{b,r}) + u^{\partial D_2}_t(\phi_{b,r})$ is a non-negative supermartingale. Taking expectations and 
$ r \downarrow 2$ we find  
$E[u^{\partial D_2}_{\infty}(1)] \leq \int_{D_2} \phi_{b,2}(x) (\mu(dx) + dx)$.
As $ b \to \infty$, $\phi_{b,2}(x)$ 
decreases to zero as $b \to \infty$ for $x \in D_2$ (argue, for example, using the probabilistic representation) 
and therefore uniformly for $x \in \overline{D_1}$. 
This implies that $E[u^{\partial D_2}_{\infty}(1)] \to 0$
as $b \to \infty$, uniformly over $\mu$ as in the statement of the lemma.

For a solution to (\ref{death1spde2}) on $D_3$ we apply the spatial Markov property Lemma \ref{sMp}
with the subdomain $D_2$. Using also extinction estimates as in (\ref{globaltsupport2})
this shows that, when 
$\mu$ supported in $\overline{D_1}$, 
\[
E[u^{\partial D_3}_{\infty}(1)] \leq E[ u^{\partial D_2}_{\infty}(1)] \quad \mbox{and} \quad
P[u^{\partial D_3}_{\infty}(1) >0] \leq E[ 1 - e^{-u^{\partial D_2}_{\infty}(w)}] \leq 
E[u^{\partial D_2}_{\infty}(w)]
\]
provided that $b^{-(4-d)/d} \Delta w \leq w^2/2$ on $D_3$ and 
$\inf\{w(x): x \in \partial D_r\} \to \infty$ as $r \uparrow 3$. 
The test function $w(x) = 12 \sum_{i=1}^d (x_i+3)^{-2} + (3-x_i)^{-2}$ satisfies these
requirements provided $b \geq 1$ and also $\sup \{w(x): x \in \partial D_2 \} < \infty$. The lemma
now follows from our control on the expected exit mass $E[u^{\partial D_2}_{\infty}(1)]$. \qed

The proof that $\Psi(\beta) < \beta$ will follow from Lemma \ref{deathsufficiencylemma}
once we have shown that on a box $D$, the block estimates
$ Q^{D,\beta,\gamma}_{\mu,1} [ U^{\partial D}_{\infty} \neq 0 ]$
and $ Q^{D,\beta,\gamma}_{\mu,1} [ U^{\partial D}_{\infty}(1)]$
are continuous as $\gamma \uparrow \beta$, uniformly over $\mu$ supported in a certain strict sub-box 
and with a certain bounded total mass.
This follows by a change of measure argument. Indeed the derivative 
$dQ^{D,\beta,\gamma}_{\mu,1}/dQ^{D,\beta,\beta}_{\mu,1}$ on $\mathcal{U}_t$ is given by
the exponential martingale $\mathcal{E}_t((\beta-\gamma)M(1))$ where $M_t(1)$ is the martingale part of $U_t(1)$. 
Note that $[M(1)]_t = U_{[0,t]}(1)$. The uniform integrability of
$\mathcal{E}_t((\beta-\gamma)M(1))$, when $\beta - \gamma$ is small, follows from the finiteness of
the exponential moments 
\[
Q^{D,\beta,\beta}_{\mu,1} \left[\exp(\lambda U_{[0,\infty)}(1)) \right] 
\leq Q^{D,0}_{\mu} \left[\exp(\lambda U_{[0,\infty)}(1)) \right] =
e^{ \mu(\phi^{(\lambda)})} < \infty
\]
where $\phi^{(\lambda)}$ solves, for small enough $\lambda>0$, $\Delta \phi^{(\lambda)} = -(1/2) (\phi^{(\lambda)})^2 - \lambda$
on $D$ and $\phi^{(\lambda)} = 0 $ on $\partial D$. Then 
\begin{eqnarray*}
 \left| Q^{D,\beta,\gamma}_{\mu,1} \left[ U^{\partial D}_{\infty}(1) \right] -
Q^{D,\beta,\gamma}_{\mu,1} \left[ U^{\partial D}_{\infty}(1) \right] \right|^2 
& \leq & \left( Q^{D,\beta,\beta}_{\mu,1} \left[ U^{\partial D}_{\infty}(1) 
(\mathcal{E}_{\infty}((\beta-\gamma) M(1))-1) \right] \right)^2 \\
& \leq &  Q^{D,\beta,\beta}_{\mu,1} \left[ (U^{\partial D}_{\infty}(1))^2 \right] 
 Q^{D,\beta,\beta}_{\mu,1} \left[ (\mathcal{E}_{\infty}((\beta-\gamma) M(1))-1)^2 \right].
\end{eqnarray*}
Then $Q^{D,\beta,\beta}_{\mu,1} [ (U^{\partial D}_{\infty}(1))^2 ] \leq Q^{D,0}_{\mu} [ (U^{\partial D}_{\infty}(1))^2 ]$
can be bounded in terms of the total mass $\mu(1)$ while, when $\beta-\gamma$ is small enough, 
\begin{eqnarray*}
&& \hspace{-.4in} Q^{D,\beta,\beta}_{\mu,1} \left[ (\mathcal{E}_{\infty}( (\beta-\gamma) M(1))-1)^2 \right] \\
& = & Q^{D,\beta,\beta}_{\mu,1} \left[ e^{2(\beta-\gamma) M_{\infty}(1) - (\beta-\gamma)^2 U_{[0,\infty)}(1)} \right] - 1 \\
& \leq & \left( Q^{D,\beta,\beta}_{\mu,1} \left[ \mathcal{E}_{\infty}(4(\beta-\gamma)M(1)) \right] \right)^{1/2}
\left( Q^{D,\beta,\beta}_{\mu,1} \left[ e^{6(\beta-\gamma)^2 U_{[0,\infty)}(1)} \right] \right)^{1/2} - 1 \\
& \leq & \left( Q^{D,0}_{\mu} \left[ e^{6(\beta-\gamma)^2 U_{[0,\infty)}(1)} \right] \right)^{1/2} - 1 \\
& = & \left(\exp(\mu( \phi^{(\lambda)}) \right)^{1/2} - 1
\end{eqnarray*}
where $\lambda= 6 (\beta- \gamma)^2$. Since $\phi^{(\lambda)}$ decreases to $0$ as $\lambda \downarrow 0$, this 
shows the required continuity in $\gamma$. The argument for 
$ Q^{D,\beta,\gamma}_{\mu,1} [ U^{\partial D}_{\infty} \neq 0 ]$ is similar. 
\subsection{Proof that $\Psi(\beta) \leq c_1 \beta^2$ in $d=3$} \label{s7.3}
We sketch a rather simple argument for death, based on the decomposition in 
Lemma \ref{comparison2} (i). This decomposition suggests we may construct the total occupation measure for a 
solution to (\ref{spde}) in two parts: $u^-_{[0,\infty)}$ from a process $u^-$ with
law $Q^{\R^d,\gamma}_{\mu}$, and $u^+_{[0,\infty)}$ from a solution $u^+$ to (\ref{spde}) 
that conditional on $\sigma\{u^-\}$ has initial condition $u^+_0(dx) = \beta (1- \exp(-u^-(0,\infty,x))dx$. 
The expected initial mass
$E[u^+_0(1)]$ can be exactly calculated via the Laplace functional of $u^-(0,\infty,x)$ and shown to 
be bounded by $c_0 \beta \gamma^{-1/2} \mu(1)$. This suggests, when $c_0 \beta \gamma^{-1/2}<1$, that this 
decomposition can be iterated leading to a convergent geometric series and a finite total occupation measure.  

We now formalise this argument, working on finite domains and with exit measures. 
\begin{Lemma} \label{death3keylemma}
When $c_0 \beta \gamma^{-1/2} <1$ we have 
\[
Q^{D_L,\beta,\gamma}_{\mu,1} \left[ U^{\partial D_L}_{\infty}(1) \right] 
\leq \mu(\phi^{(L)}) + \mu(1) \sum_{k=1}^{\infty} (c_0 \beta \gamma^{-1/2})^k
\]
where $\phi$ solves $\Delta \phi^{(L)} = \gamma \phi^{(L)}$ on $D_L$ and $\phi^{(L)} =1$ on $\partial D_L$, 
and $c_0 = \int_{\R^3} \psi(x) dx$ where $(\psi(x): \R^3  \setminus \{0\})$ is the maximal solution to
\begin{equation}
\mbox{$\Delta \psi = \frac12 \psi^2 + \psi$ on $\R^3 \setminus \{0\}$, $\quad \psi(x) \to 0$ as $|x| \to \infty \;$ and
$\; \psi(x) \to \infty$ as $|x| \to 0$.}
\label{psieqn}
\end{equation}
\end{Lemma}

\noindent \textbf{Remark.} The maximal solution $\psi$ can be constructed as follows: non-negative 
solutions $\psi^{(\epsilon)}$ to $\Delta \psi^{(\epsilon)} = (1/2) (\psi^{(\epsilon)})^2 + \psi^{(\epsilon)}$
on $\{x: 0< \epsilon < |x| < \infty\}$, with boundary conditions
$\psi^{(\epsilon)}(x) \to 0$ as $|x| \to \infty$ and
$\psi^{(\epsilon)}(x) \to \infty$ as $|x| \downarrow \epsilon$ exist and are unique (take limits through 
increasing boundary conditions, using the upper bound $C(|x|-\epsilon)^{-2}$, for existence 
and use the maximum principle for uniqueness.)
Moreover $\psi^{(\epsilon)}$ decrease as $\epsilon \downarrow 0$ to function $\psi$ defined on $\R^3 \setminus \{0\}$. 
Interior regularity estimates show that that $\psi$ solves (\ref{psieqn}). The boundary condition for $\psi^{(\epsilon)}$ 
at $|x|=\epsilon$ imply, by a comparison argument, that $\psi$ is the maximal solution. Finally, another comparison shows 
that $0 \leq \psi \leq 4 |x|^{-2}$ and that $\psi$ has exponential decay at infinity so that 
$c_0 = \int \psi(x) dx < \infty$.

\noindent \textbf{Proof.} $L,\beta,\gamma$ are fixed throughout. For convenience in this proof we
let $F(\mu,f) = Q^{D_L,\beta,\gamma}_{\mu,f} [ U^{\partial D_L}_{\infty}(1) ]$ and let 
$\Xi(\mu,d\nu)$ be the the measurable kernel given by the 
law of $\beta(1-\exp(-U(0,\infty,x))dx$ on $\mathcal{M}(D_L)$ under
$Q^{D_L,\gamma}_{\mu}$. Then the decomposition in Lemma \ref{comparison2} (i), monotonicity and
the exact formula for exit measures of superprocesses (derived, say, from (\ref{LFD2})) imply  
\begin{eqnarray}
F(\mu,1) & = & Q_{\mu}^{D_L,\gamma} \left[ U^{\partial D_L}_{\infty}(1) \right] 
+ Q^{D_L,\gamma}_{\mu} \left[ F(\beta(1-e^{-U(0,\infty,x)}dx, e^{-U(0,\infty,x)}) \right] \nonumber \\
& \leq & Q_{\mu}^{D_L,\gamma} \left[ U^{\partial D_L}_{\infty}(1) \right] 
+ Q^{D_L,\gamma}_{\mu} \left[ F(\beta(1-e^{-U(0,\infty,x)}dx, 1) \right] \nonumber \\
&=& \mu(\phi^{(L)}) + \int F(\nu,1) \, \Xi(\mu,d\nu). \label{Fiterate}
\end{eqnarray}
Define the range by $\mathcal{R} = \cup_{\delta >0} \, \overline{ \cup_{t \geq \delta} \, \mbox{support}(U_t)}$, where $\mbox{support}(\mu)$ is 
the closed support of $\mu$. Then
\begin{eqnarray*}
\int \nu(1) \, \Xi(\mu,d \nu) & = & \beta \int_{D_L}  Q^{D_L,\gamma}_{\mu} \left[ 1- e^{-U(0,\infty,x)} \right] \, dx \\
& \leq & \beta \int_{D_L}  Q^{D_L,\gamma}_{\mu} \left[ \{x\} \in \mathcal{R} \right] \, dx
= \beta \int_{D_L} \mu(\tilde{\psi}^{(\gamma,x)}) dx 
\end{eqnarray*}
where $\tilde{\psi} = \tilde{\psi}^{(\gamma,x)}$ defined for $y \in \overline{D_L} \setminus \{x\}$ is the maximal solution 
to 
\[ \mbox{$\Delta \tilde{\psi} = \frac12 \tilde{\psi}^2 + \gamma \tilde{\psi}$ on $D_L \setminus \{x\}$, 
$\quad \tilde{\psi} = 0$ on $\partial D_L \;$ and $ \; \tilde{\psi}(y) \to \infty$ as $|y-x| \to 0$.} 
\]
The final equality follows using mostly the arguments in \cite{perkins1} Theorem III.5.7 (a) 
(and constructing $\tilde{\psi}^{(\gamma,x)}$ as in the above remark). 
A comparison argument shows that $\tilde{\psi}^{(\gamma,x)}(y)
\leq \gamma \psi(\gamma^{1/2}(y-x))$ and this leads to the bound
$\int \nu(1) \, \Xi(\mu,d \nu) \leq c_0 \beta \gamma^{-1/2} \mu(1)$.  
Note that $\phi^{(L)} \leq 1$. We can now interate the inequality (\ref{Fiterate}) 
and the lemma will follow once we have show that the remainder term converges to zero, and this 
will follow from the fact that $F(\mu,1) \to 0$ when $\mu(1) \to 0$. To see this, note that
\begin{eqnarray*}
F(\mu,1) & \leq & Q^{D_L,\beta,\gamma}_{\mu,1} \left[ U^{\partial D_L}_{\infty}(1) \I(\tau >t)  \right]
+ Q^{D_L,\beta,\gamma}_{\mu,1} \left[ U^{\partial D_L}_t(1) \right] \\
& \leq & \left(  Q^{D_L,\beta,\gamma}_{\mu,1} \left[ (U^{\partial D_L}_{\infty}(1))^2  \right]
 Q^{D_L,\beta,\gamma}_{\mu,1} \left[ \tau > t \right] \right)^{1/2} 
+ Q^{D_L,\gamma-\beta}_{\mu} \left[ U^{\partial D_L}_t(1) \right] \\
& \leq & \left(  Q^{D_L,\gamma}_{\mu + \beta \I(D_L) dx} \left[ (U^{\partial D_L}_{\infty}(1))^2  \right]
 Q^{D_L,\beta,\gamma}_{\mu,1} \left[ \tau > t \right] \right)^{1/2} 
+ \mu(\theta_t)
\end{eqnarray*}
where $\partial \theta = (\beta-\gamma) \theta$ on $[0,t] \times D_L$, $\theta_0 =0$
and $\theta =1$ on $[0,t] \times \partial D_L$. Note the term 
$ Q^{D_L,\gamma}_{\mu + \beta \I(D_L) dx} \left[ (U^{\partial D_L}_{\infty}(1))^2  \right]$
is bounded by $C(L,\beta,\gamma)<\infty$ for $\mu(1) \leq 1$. The term
$ Q^{D_L,\beta,\gamma}_{\mu,1} \left[ \tau > t \right]$ converges to zero as $t \to \infty$, uniformly over 
$\mu(1) \leq 1$, for instance by the estimate in Lemma \ref{simpledeath}. These together allow us to see
that $F(\mu,1) \to 0$ when $\mu(1) \to 0$ completing the proof. 
\qed

To complete the proof that $\Psi(\beta) \leq c_1 \beta^2$ we 
will apply the lemma above to a sequence of domains $D(n) := D_{n \gamma^{-1/2}}$, and show that the expected 
exit measures decay geometrically. 
Scaling shows that the solution $\phi^{(\gamma^{-1/2})}$ from the above lemma
satisfies $\phi^{(\gamma^{-1/2})}(x) = \hat{\phi}(\gamma^{1/2}x)$ where
$ \Delta \hat{\phi} = \hat{\phi}$ on $D_1$ and $\hat{\phi}=1$ on $\partial D_1$. 
Moroever a comparison argument shows that $\phi^{((n+1) \gamma^{-1/2})}(x) \leq 
\phi^{(\gamma^{-1/2})}(0) = \hat{\phi}(0)<1$ for $x \in \partial D(n)$. The lemma therefore implies that
\begin{equation}
Q^{D(n+1),\beta,\gamma}_{\mu,1} \left[ U^{\partial D(n+1)}_{\infty}(1) \right] 
\leq \mu(1) \left( \hat{\phi}(0) + (c_0 \beta \gamma^{-1/2})(1-c_0 \beta \gamma^{-1/2})^{-1} \right)
\label{temp303} 
\end{equation}
whenever $\mu$ is supported on $\partial D(n)$. 
We may now choose $\gamma = c_1 \beta^2$ with $c_1 < \infty$ large enough that the right hand side
of (\ref{temp303}) is at most $(1/2)(1+\hat{\psi}(0))\mu(1) < \mu(1)$. Iterating shows when 
$\mu$ is supported on $\partial D(1)$ that 
\[
Q^{D(n+1),\beta,c_2 \beta^2}_{\mu,1} \left[ U^{\partial D(n+1)}_{\infty}(1) \right] 
\leq ((1+\hat{\phi}(0))/2)^n \mu(1).
\]
This implies certain death. Indeed the spatial Markov property and the 
extinction estimate (\ref{globaltsupport2}) show that
\[
Q^{D(n),\beta,c_2 \beta^2}_{\mu,1} \left[ U^{\partial D(n)}_{\infty} \neq 0 \right] 
\leq Q^{D(n),\beta,c_2 \beta^2}_{\mu,1} \left[ 1- e^{-U^{\partial D(n-1)}_{\infty} (w)} \right]
\to 0 \quad \mbox{as $n \to \infty$} 
\]
where 
\[
w(x) = \sum_{i=1}^3 \left( 2 \beta + \frac{12}{(x_i+ n \gamma^{-1/2} )^2} \right) 
+ \left( 2 \beta + \frac{12}{(n \gamma^{-1/2}-x_i)^2} \right) 
\]
satisfies $\Delta w \leq (1/2) w^2 - \beta w$ on $D(n)$.  
Certain death follows from Lemmas \ref{deathcharacterization1} and \ref{deathcharacterization2}.

\noindent \textbf{Remark}
The approach above can be applied to the case of $d=2$ and small $\beta$. An analysis of 
the first moment shows that when $\gamma = \exp(-C/\beta)$ the above series construction converges
for suitable $C$ and hence certain death. We do not include the details since the argument is quite crude 
and it seems easier to conjecture that there is death when $\gamma=0$ for small $\beta$. 
Note however that when $\gamma=0$ a first moment argument will not 
show death since the first moment $Q^{D_R,\beta,0}_{\mu,1}[U^{\partial D_R}_{\infty}(1)] \to \infty$ 
as $R \to \infty$. We hope to comment on this in a subsequent paper. 

\end{document}